\documentclass[11pt,letterpaper]{amsart}
\usepackage[
    top=3cm, bottom=2.9cm, inner=3cm, outer=3cm,
    marginparwidth=2cm 
    ]{geometry}

\usepackage{amssymb,latexsym, amsmath, amsxtra, mathrsfs, bm}

\usepackage[dvips]{graphics}
\usepackage[all]{xy}
\usepackage{subcaption,xcolor}
\usepackage[T1]{fontenc}

\usepackage{verbatim}
\usepackage[abs]{overpic}
\usepackage{xurl}
\usepackage{mathtools}
\usepackage{hyperref}

\usepackage{tikz}           
\usepgflibrary{arrows.meta}                 
\usetikzlibrary{decorations.pathreplacing}  
\usepackage{tikz-cd}        

\newcommand*\circledd[1]{\tikz[baseline=(char.base)]{
            \node[shape=circle,draw,inner sep=2pt] (char) {#1};}}

\newcommand{\numbercircle}[5]{%
    \node[
        circle,
        draw=#3,
        fill=#3!#4!white,
        text opacity=1,
        text=black,
        line width=1pt,
        minimum size=5mm,
        inner sep=0pt
    ] (#1) at #2 {#5};
}

\DeclareMathAlphabet{\mathbbold}{U}{bbold}{m}{n}

\allowdisplaybreaks[4]

\theoremstyle{plain}
        \newtheorem{theorem}{Theorem}[section]
        \newtheorem*{theorem*}{Theorem}
        \newtheorem*{conj*}{Conjecture}
        \newtheorem{lemma}[theorem]{Lemma}
        \newtheorem{prop}[theorem]{Proposition}
        \newtheorem{cor}[theorem]{Corollary}

\theoremstyle{definition}
        \newtheorem{definition}[theorem]{Definition}

\theoremstyle{remark}
        \newtheorem*{remark}{Remark}

\numberwithin{equation}{section}
\numberwithin{theorem}{section}
\numberwithin{table}{section}
\numberwithin{figure}{section}

\renewcommand{\le}{\leqslant}
\renewcommand{\leq}{\leqslant}
\renewcommand{\ge}{\geqslant}
\renewcommand{\geq}{\geqslant}

\newcommand{\var}  {\operatorname{var}}

\newcommand{\id} {\operatorname{id}}

\newcommand{\supp}{\operatorname{supp}}
\newcommand{\ve}{\varepsilon}
\newcommand{\R}{\mathbb{R}}

\newcommand{\N}{\mathbb{N}}

\providecommand{\abs}[1]{\lvert#1\rvert}
\providecommand{\Absbig}[1]{\bigl\lvert#1\bigr\rvert}
\providecommand{\Absbigg}[1]{\biggl\lvert#1\biggr\rvert}

\providecommand{\norm}[1]{\|#1\|}

\renewcommand{\:}{\colon}

\newcommand{\marg}{\operatorname{marg}}

\newcommand{\LIP}{\operatorname{LIP}}

\newcommand{\MMM}{\mathcal{M}}

\newcommand{\vertiii}[1]{{\left\vert\kern-0.25ex\left\vert\kern-0.25ex\left\vert #1
    \right\vert\kern-0.25ex\right\vert\kern-0.25ex\right\vert}}

\newcommand{\wt}[1]{\widetilde{#1}}

\renewcommand{\=}{\coloneqq}

\newcommand{\Mmax}{\MMM_{\operatorname{max}}}

\newcommand{\ouv}{\operatorname{sep}}

\newcommand{\Lip}{\operatorname{Lip}}

\newcommand{\Lock}{\operatorname{Lock}}

\newcommand{\conc}{\operatorname{conc}}
\newcommand{\cou}{\operatorname{cou}}

\newcommand{\bF}{\bm{F}}

\newcommand{\bQ}{\bm{Q}}

\newcommand{\bV}{\bm{V}}

\newcommand{\ba}{\bm{a}}
\newcommand{\bb}{\bm{b}}
\newcommand{\bc}{\bm{c}}

\newcommand{\be}{\bm{e}}
\newcommand{\bg}{\bm{g}}

\newcommand{\bp}{\bm{p}}

\newcommand{\bs}{\bm{s}}

\newcommand{\bw}{\bm{w}}
\newcommand{\bx}{\bm{x}}

\newcommand{\bz}{\bm{z}}

\newcommand{\boo}{\mathfrak{o}}

\newcommand{\PO}{\mathcal{O}}
\newcommand{\Path}{\mathcal{P}}

\newcommand{\cC}{\mathcal{C}}

\newcommand{\cF}{\mathcal{F}}

\newcommand{\cH}{\mathcal{H}}
\newcommand{\cI}{\mathcal{I}}

\newcommand{\cM}{\mathcal{M}}

\newcommand{\cO}{\mathcal{O}}
\newcommand{\cP}{\mathcal{P}}

\newcommand{\cV}{\mathcal{V}}

\newcommand{\fC}{\mathfrak{C}}

\newcommand{\sE}{\mathscr{E}}

\newcommand{\sN}{\mathscr{N}}

\newcommand{\sP}{\mathscr{P}}

\newcommand{\sU}{\mathscr{U}}

\newcommand{\tg}{\widetilde{g}}

\newcommand{\black}[1]{\color{black}}

\newcommand{\mpe}{\beta}
\newcommand{\ext}{\operatorname{ext}}
\newcommand{\smallo}{o}

\begin{document}

\title[Prevalence of periodicity of maximizing measures]{Prevalence of the periodicity of maximizing measures}

\makeatletter
\renewcommand{\@setauthors}{%
  \begingroup
  \trivlist
  \centering\normalsize \@topsep30\p@\relax
  \advance\@topsep by -\baselineskip
  \item\relax
  \authors
  \vspace{14pt}
  \endtrivlist
  \endgroup
}
\makeatother

\author[Zelai~Hao, Illya~Koval, and Zhiqiang~Li]{%
  \begin{tabular}{c}Zelai~Hao\\[-3pt]\Small Peking University
  \end{tabular}\;
  \begin{tabular}{c}Illya~Koval\\[-3pt]\Small Institute of Science and Technology Austria
  \end{tabular}\;
  \begin{tabular}{c}Zhiqiang~Li\\[-3pt]\Small Peking University  \end{tabular}
}



\begin{abstract}
For a continuous map $T\: X\rightarrow X$ on a compact metric space $(X,d)$, we say that a function $f\: X \rightarrow \R$ has the property $\sP_T$ if its time averages along forward orbits of $T$ are maximized at a periodic orbit. In this article, we establish the Hunt--Ott conjecture (cf.~\cite{HO96a,HO96b}), that is, the property $\sP_T$ is prevalent (in the sense of Hunt--Sauer--Yorke) in the space of Lipschitz functions with respect to a standard metric for subshifts of finite type. This result is a strengthening of \cite[Theorem~A]{BZ16} and \cite[Theorem~1.1]{DLZ24}, confirms the prediction mentioned in the ICM proceeding contribution of Bochi (\cite[Section~1]{Boc18}) suggested by experimental evidence. 
\end{abstract}

\maketitle

\tableofcontents

\section{Introduction}
In this article, we prove an infinite-dimensional version of the Hunt--Ott conjecture for Lipschitz potentials in the context of symbolic dynamics. 

In the late 1990s,  Hunt, Ott, and Yuan conjectured that for a typical chaotic system $T$, a typical smooth or Lipschitz potential $f$ admits a maximizing measure supported on a periodic orbit, where the typicality on the potential $f$ was both in the sense of probability (in terms of the Lebesgue measure of the parameter(s) for the potential $f$ in a finite-dimensional parameter space) (\cite[Conjecture]{HO96a} and \cite[Conjecture~2]{HO96b}) and in the sense of topology (topological genericity) (\cite[Conjecture~1.1]{YH99}). While considerable progress has been made on the topological direction since \cite{Co16}, much less is known on the probabilistic side. Here we establish a prevalence result for Lipschitz potentials with respect to standard metrics on subshifts of finite type. 

\begin{theorem}	\label{thm:main}
	Let $(\Sigma,\sigma)$ be a one-sided subshift of finite type equipped with some standard metric $d$. Then the set of Lipschitz potentials satisfying the locking property is prevalent in the space $\Lip(\Sigma,d)$ of Lipschitz potentials. 
\end{theorem}

Here a potential function $f\in\Lip(\Sigma,d)$ is said to satisfy the \emph{locking property} if $f$ admits a unique maximizing measure $\mu$, which is supported on a periodic orbit, and there exists a nonempty neighbourhood $U$ of $f$ in $\Lip(\Sigma,d)$ such that for each $g\in U$, $\mu$ is the unique maximizing measure of $g$.

The notion of \emph{prevalence}, introduced by \cite{HSY92}, is an infinite-dimensional analog of the notion of ``Lebesgue almost every''. More precisely, let $V$ be a topological vector space over $\R$ equipped with a complete metric. A Borel set $S \subseteq V$ is called \emph{shy} if there is a compact Borel subset $\cH \subseteq V$ and a probability measure $\nu$ supported on $\cH$ such that for every $f \in V$ we have
\begin{equation*}
	\nu(\{g \in \cH:f+g \in S\}) =0.
\end{equation*}
A set $S\subseteq V$ is called \emph{prevalent} if its complement is contained in a shy set.

\subsection{Setting and terminology}

For a dynamical system on a compact metric space $(X,d)$ with a continuous map $T$, let $\MMM (X, T)$ denote the set of invariant Borel probability measures of $(X,T)$. Recall the \emph{ergodic supremum}

\begin{equation*}\label{e:ergodicmax}
  \mpe(X, T, f)\=\sup\limits_{\mu\in\MMM(X,T)}\int\! f \, \mathrm{d}\mu
  =\max\limits_{\mu\in\MMM(X,T)}\int\! f \, \mathrm{d}\mu.
\end{equation*}

Let $\Mmax(X,T,f)$ be the set of maximizing measures. If there is no confusion on the space $X$ or the map $T$, we may omit them in the notations $\mpe(X, T, f)$ and $\Mmax(X, T, f)$.

Denote by $\Lip(X, d)$ the space of functions with finite \emph{minimal Lipschitz constant}
\begin{equation*}
	\LIP_d (f) \= \sup\{ \abs{f(x) - f(y)} / d(x,y)  :  x,\,y \in X,\, x \neq y \}.
\end{equation*}
We equip this space with the Lipschitz norm $\norm{\cdot}_{\Lip}$, given by $\norm{  f  }_{\Lip} \= \norm{ f }_{\infty} + \LIP_d(  f )$, i.e., the sum of the sup-norm and the minimal Lipschitz constant $\LIP_d(  f )$ of $f$.

For a typical $f$ in $\Lip(X, d)$, how do the measures in $\Mmax(f)$ look like? There is an easy way to construct simple invariant measures: for every periodic orbit $\boo$ of $(X, T)$ one can take a unique invariant measure $\mu_{\boo}$ supported on this orbit. Then, one have
$    \mpe(\boo, f) = \int\! f\, \mathrm{d}\mu_{\boo}$.
Our goal is to show that for a \emph{typical} $f$, the maximizing measure $\mu_{\max}$ is unique and that it is supported on a periodic orbit. In particular, we have $\mpe(X, f)=\mpe(\boo, f)$ for some periodic orbit $\boo$. 
However, $f$ being typical can mean different things. In this article, more specifically, we will consider a \emph{prevalent} Lipschitz function $f$.

A function $f \in \Lip(X, d)$ is said to satisfy the \emph{periodic optimization property} (abbr.\ PO property) if there is a unique maximizing measure of $f$ and this measure is supported on a periodic orbit. Moreover, $f$ is said to have the \emph{locking property} if $f$ has PO property and $\Mmax(g) = \Mmax(f)$ for all functions $g\in \Lip(X, d)$ sufficiently close to $f$ in the Lipschitz norm. This means that if $f$ satisfies the locking property then the maximizing measure is locked to be a unique periodic measure under a small perturbation. The set of Lipschitz functions $f \in \Lip(X, d)$ with the locking property is denoted by $\Lock_{d} (X,\sigma)$. The locking property clearly depends on the metric $d$.

Consider a one-sided subshift of finite type $(\Sigma_{\bF}, \sigma )$ on a set $\Sigma_{\bF} \subseteq \bQ^{\N} = \{0,\,1,\, \cdots,\, q-1\}^{\N}$ of sequences on a $q$-symbol alphabet that do not contain the set of forbidden words $\bF$. Here the \emph{left-shift operator} $\sigma  \: \Sigma_{\bF} \rightarrow \Sigma_{\bF}$ is given by
\begin{equation*}
\sigma  ( \{w_i\}_{i\in\N} ) = \{w_{i+1}\}_{i\in\N}  \quad \text{ for } \{w_i\}_{i\in\N} \in \Sigma_{\bF}.
\end{equation*}

The space $\Sigma_{\bF}$ admits the following metric. For a strictly decreasing sequence $\ba = \{ a_n \}_{n\in\N}$ of positive numbers converging to zero, define a function $d_{\ba} \: \Sigma_{\bF} \times \Sigma_{\bF} \rightarrow \R$ by assigning, for each pair of distinct $\bw= \{w_n\}_{n\in\N}, \, \bz=\{z_n\}_{n\in\N} \in \Sigma_{\bF}$, 
\begin{equation}   \label{e:def_d_a}
d_{\ba} ( \bw, \bz ) \= a_n,
\end{equation}
where $n$ is the smallest positive integer with the property that $w_n \neq z_n$, and by setting $d_{\ba}  ( \bz, \bz ) \= 0$ for each $\bz \in \Sigma_{\bF}$. It is easy to see that $d_{\ba}$ is a metric on $\Sigma_{\bF}$. One obtains a \emph{standard metric} on $\Sigma_{\bF}$ by setting $a_n \= \theta^{n}$ for some $\theta \in (0, 1)$.

\subsection{Previous results}

A remarkable phenomenon, suggested by numerical experiments of Hunt and Ott in the 1990s \cite{HO96a,HO96b}, is that for a chaotic map $T$, a typical potential $\phi$ within a sufficiently regular class, such as the class of smooth or H\"{o}lder potentials, has a unique maximizing measure supported on some $T$-periodic orbit $\boo$, i.e., the $\phi$-average $\frac{1}{\abs{\boo}}\sum_{x\in\boo}\phi(x)$ is strictly larger than the $\phi$-average $\int\!\phi\,\mathrm{d}\mu$ with respect to any $T$-invariant Borel probability measure $\mu$ not supported on $\boo$. In \cite{HO96a, HO96b}, Hunt and Ott investigated finite-dimensional parametrized families of potentials and conjectured that the set of parameter values for which maximizing measures are periodic has \emph{full Lebesgue measure}; subsequently Yuan and Hunt \cite{YH99}  formulated a \emph{topological} analog of the Hunt--Ott Conjecture. These two conjectures are among the central problems in \emph{ergodic optimization}, the study of maximizing measures; for comprehensive surveys, see \cite{Boc18,Je06,Je19}. Despite fruitful results stimulated by the Yuan--Hunt Conjecture (see e.g.\ \cite{Bou01,Bou08,BQ07,Co16,CLT01,GSZ25,HHJL25,HLMXZ25,HJXZ26,HJL25,JLS26,LZ25,Mo08}), the Hunt--Ott Conjecture has been less thoroughly investigated. 

Motivated by the hope for an analog of the notions of ``Lebesgue almost every'' and ``Lebesgue measure zero'' in infinite-dimensional spaces such as various Banach spaces of functions studied in dynamics (cf.~\cite{Ko54}), Hunt, Sauer, and Yorke introduced the notion of \emph{prevalence} in \cite{HSY92}, which has since proven to play a central role in understanding typical behaviors in dynamics from a probabilistic (measure-theoretic) perspective. For comprehensive surveys on prevalence, see Hunt--Kaloshin \cite{HK10} and Ott--Yorke \cite{OY05}. In finite-dimensional Euclidean spaces, the notion of prevalence coincides with the notion of ``Lebesgue almost every''. Recent all general progress towards the Hunt--Ott Conjecture (cf.~\cite{BZ16,DLZ24,GSZ25}) has come from results showing that for a fixed chaotic map $T$ and a suitable infinite-dimensional  Banach space $\cV$ of regular potentials, the set of potentials in $\cV$ having a unique periodic maximizing measure is prevalent in $\cV$. More precisely, for the one-sided full shift on two symbols, Bochi and Zhang \cite{BZ16} and Ding, Li, and Zhang \cite{DLZ24} established prevalence results in spaces of potentials with rapidly decaying variations, with the latter work allowing a substantially slower rate of decay. Both works concern potentials whose variations decay faster than those of Lipschitz potentials with respect to the standard symbolic metrics. Subsequently, Gao, Shen, and Zhang proved prevalence results for analytic circle expanding maps and smooth potentials in \cite{GSZ25}, and have expressed in multiple occations that their method does not apply to symbolic settings. Our approach appears applicable to settings far beyond the symbolic case. Following the tradition in \cite{BZ16, DLZ24}, we have chosen to present it in the symbolic setting to illustrate the ideas, leaving its presentation in other settings and the development of our ideas for further investigation in ergodic optimization beyond the current phenomenon to a separate upcoming work.

\subsection{Idea of the proof}

The main focus of the proof is the analysis of periodic orbits $\boo \in\PO_*$. For a periodic orbit $\boo \in \PO_*$ two quantities of it would be crucial throughout the proof. The first of these quantities is the \emph{(minimal) period} $\abs{\boo}$ of the orbit. The second property is the \emph{gap} $\Delta(\boo)$ of the orbit, defined by
\begin{equation}    \label{eq:gap}
    \Delta (\boo) \= \begin{cases}
        \min \{ d_{\ba}(\bw, \bz): \bw, \, \bz \in \boo,\, \bw \ne \bz\} &\text{ if }\abs{\boo} > 1,\\
        a_1&\text{ if }\abs{\boo} = 1.
    \end{cases}
\end{equation}
Orbits with a small gap almost visit the same point twice during the period, see Figure~\ref{fig:4examples}. 
Here the idea using the gap of periodic orbits comes from the use of de Bruijn--Good graphs in \cite{BZ16} and \cite{DLZ24}. Note that a periodic orbit $\boo$ satisfies $\Delta(\boo) > a_n$ if and only if it corresponds to a cycle in the $n$-th level de Bruijn--Good graph.
Roughly speaking, periods and gaps control the structure of the whole proof, meaning that orbits of different periods and gaps are eliminated in different mostly independent sections. 

\begin{figure}
\begin{tikzpicture}
    \pgfmathdeclarefunction{f}{1}{\pgfmathparse{-1.45+4*#1^(-0.25)}}
    \pgfmathdeclarefunction{g}{1}{\pgfmathparse{-1.0+4*#1^(-0.25)}}


    \fill[red!20]
        plot[domain=0.5:10, samples = 50, smooth] (\x,{f(\x)}) --
        plot[domain=10:0.5, samples = 50, smooth] (\x,{g(\x)}) -- cycle;

    \fill[green!20]
        plot[domain=0.5:7.5, samples = 50, smooth] (\x,{f(\x)}) --
        plot[domain=7.5:0.5, samples = 50, smooth] (\x,{(f(\x)*f(7.5))^(0.5)}) -- cycle;

     \fill[blue!20]
        plot[domain=0.5:3.5, samples = 50, smooth] (\x,{(f(\x)*f(7.5))^(0.5)}) --
        plot[domain=3.5:2.5, samples = 50, smooth] (\x,{(f(\x)*f(7.5))^(0.5)*(\x-2.5)}) -- (0.5, 0) -- cycle;

    \fill[yellow!20]
        plot[domain=7.5:10, samples = 50, smooth] (\x,{f(\x)}) -- (10, 0) --
        plot[domain=2.5:3.5, samples = 50, smooth] (\x,{(f(\x)*f(7.5))^(0.5)*(\x-2.5)}) -- plot[domain=3.5:7.5, samples = 50, smooth] (\x,{(f(\x)*f(7.5))^(0.5)}) -- cycle;

    \draw[black, very thick] (0.5,0) -- (10,0);
    \draw[black, very thick] (0.5,0) -- (0.5,{g(0.5)});
    \draw[black, very thick] (10,0) -- (10,{g(10)});

    \draw[black] plot[domain=0.5:10, samples = 50, smooth] (\x,{f(\x)});
    \draw[black] plot[domain=0.5:7.5, samples = 50, smooth] (\x,{(f(\x)*f(7.5))^(0.5)});
    \draw[black] plot[domain=2.5:3.5, samples = 50, smooth] (\x,{(f(\x)*f(7.5))^(0.5)*(\x-2.5)});
    \draw[black, very thick] plot[domain=0.5:10, samples = 50, smooth] (\x,{g(\x)});

    \draw[->] (10,0) -- (11,0) node[below] {$\Delta(\boo)$};
    \draw[->] (0,0) -- (0,4) node[left] {$\abs{\boo}$};
    \draw[black] (0,0) -- (0.5, 0) node[below] {$\theta^n$};

    \node[] at (2.6,1.9) {\textbf{E}};
    \node[] at (1.0,2.05) {\textbf{A}};
    \node[] at (6,0.55) {\textbf{C}};
    \node[] at (1.5,0.75) {\textbf{B}};
\end{tikzpicture}
\caption{Four groups of periodic orbits ordered by their periods $\abs{\boo}$ and gaps $\Delta(\boo)$.}
\label{fig:4groups}
\end{figure}

What does it mean to ``eliminate'' an orbit? Take some finite set of orbits, for example those with $\Delta(\boo) > a_n$ for some large $n$. As invariant measures supported on periodic orbits are dense, the maximum of $\mpe( \boo, f)$ on our finite set should approximate $\mpe(\Sigma_{\bF},f)$ fairly well. Moreover, to avoid the locking property, some other orbits in this set should be competitive with the maximum, otherwise the maximal $\widehat{\boo}$ will simply support the maximizing measure. Our goal would be to understand this maximal $\widehat{\boo}$, see Lemma~\ref{lemma.betastar}. Thus, by eliminating an orbit $\boo$ we will mean proving that with high probability either the locking property holds or $\boo$ is not the maximal one.

Orbits will be eliminated in $4$ groups, depending on their period and gap. This will be done in Sections~\ref{sec:thetareduction}, \ref{sec:A}, \ref{sec:B}, and~\ref{sec:C}. The visual representation of those groups is given in Figure~\ref{fig:4groups} with orbit examples given in Figure~\ref{fig:4examples}. Overall, the division into groups is not very strict and border orbits can be eliminated in both ways. One may use this flexibility to extend the result of Theorem~\ref{thm:main} to other systems. By eliminating the orbits, we will use a bootstrapping argument to establish, with high probability, successively stronger properties of the maximal $\widehat{\boo}$ (see Theorems~\ref{thm:maim.redu} and~\ref{thm:prop}) until the sufficient condition for the locking property given in Lemma~\ref{lemma.betastar} is satisfied. Now, we discuss these $4$ groups in more detail.

\subsubsection{Group $E$ -- long orbits -- Section~\ref{sec:thetareduction}} 

In this section, we reduce Theorem~\ref{thm:main} to the case when the topological entropy of $\Sigma_{\bF}$ is small, see Theorem~\ref{thm:maim.redu} for more details. The main tool for this is the \emph{prevalence of zero entropy} of maximizing measures (Theorem~\ref{prop:pre.zero.entropy}).  As a side effect, this reduction removes all the periodic orbits $\boo$ with $\abs{\boo} > \theta^{-n/3}$ and $\Delta(\boo) > a_n$ from consideration because by the pigeonhole principle orbits with bounded gap cannot be so long for small entropy. So, in summary, this argument deals with exponentially long orbits. 

To prove the prevalence of zero-entropy of maximizing measures, we mainly combined the ideas in \cite{Mo08} in proving genericity of every maximizing measure having zero entropy, the Hilbert brick construction in \cite{BZ16}, and the probability estimation in \cite{DLZ24}. More precisely, following \cite{Mo08}, we use a combination of the Ma\~n\'e lemma and the Bressaud--Quas closing lemma to prove that for sufficient large $N \in \N$, a perturbation of scale $\smallo \bigl( n^{1/ \smallo(1)} \bigr)$ in the $\frac{\log n}{\smallo(1)}$ level of the Hilbert brick $\cH[\bb]$ can make the upper bound of entropy of maximizing measures arbitrarily small. Based on this, we use an analogous idea in \cite[Proposition~5.1]{DLZ24} which we call $\delta$-argument to prove that 
\begin{equation*}
    \mathbb{P} \Bigl( \Bigl\{ g \in \cH[\bb] : \sup_{\mu\in \Mmax(f+g)} h_{\mu}(\sigma) < \ve \Bigr\} \Bigr) <\ve
\end{equation*}
for all $\ve>0$ and $f\in \Lip(\Sigma_{\bF}, d_{\ba})$, yielding the prevalence of zero-entropy of maximizing measures. The $\delta$-argument is introduced in Subsection~\ref{subsec:darg} and is used for other groups as well. 

We use the compactness of the space of compact subsets of $\Sigma_{\bF}$ to divide the space of Lipschitz functions with zero entropy into countably many subsets $\{L_n\}_{n\in \N}$ such that for functions in each subset $L_n$, every maximizing measure is supported on a low entropy subshift $Y_n$ (see Corollary~\ref{cor:sep.lowentropy}). For functions $f$ in each subset $L_n$, we reduce the ergodic optimization studying to that restricted to $(Y_n, \sigma)$, which satisfies the assumptions in Theorem~\ref{thm:maim.redu}. Theorem~\ref{thm:maim.redu} shows that the subset of functions without the locking property in each $L_n$ is shy. This and the prevalence of zero-entropy of maximizing measures give Theorem~\ref{thm:main}.

\subsubsection{Groups A--C}
Our next improvement over \cite{DLZ24} concerns how we measure the competitiveness of two orbits. Whereas \cite{DLZ24} simply takes the difference of their ergodic suprema, we normalize this difference by the closeness factor; see \eqref{eq:betastardef}. This helps, because if two orbits are close to each other (say differ in one symbol), their ergodic suprema would be close, but they would not necessarily be competitive. This idea will be used for Groups A--C.

\begin{figure}
\begin{tikzpicture}
    \pgfmathsetseed{12346}
    \foreach \i in {1,...,20} {
        \pgfmathsetmacro{\x}{2*rnd-3}
        \pgfmathsetmacro{\y}{2*rnd-0.5}
        \pgfmathtruncatemacro{\l}{mod(\i,3)}
        \pgfmathsetmacro{\o}{15*\l}

        \numbercircle{R\i}{(\x,\y)}{red}{\o}{\l}
    }
    \draw[->, thick] (R1) -- (R13);
    \draw[->, thick] (R13) -- (R14);
    \draw[->, thick] (R14) -- (R16);
    \draw[->, thick] (R12) to[out=-30, in=80] (R6);
    \draw[->, thick] (R6) to[out=-50, in=30] (R15);
    \draw[->, thick] (R15) to[out=-120, in=5] (R8);
    \draw[->, thick] (R18) -- (R10);
    \draw[->, thick] (R17) -- (R3);
    \draw[->, thick] (R5) to[out=30, in=135] (R12);
    \draw[->, thick] (R8) to[out=-135, in=-90] (R20);
    
    \numbercircle{Y1}{(0,0)}{orange!70!green}{15}{1}
    \numbercircle{Y2}{(1,0)}{orange!70!green}{15}{1}
    \numbercircle{Y3}{(0, 1)}{orange!70!green}{0}{0}
    \numbercircle{Y4}{(1,1)}{orange!70!green}{0}{0}
    \draw[->, thick] (Y3) -- (Y4);
    \draw[->, thick] (Y4) -- (Y2);
    \draw[->, thick] (Y2) -- (Y1);
    \draw[->, thick] (Y1) -- (Y3);

    \numbercircle{G1}{(2.9,0.5)}{green!60!black}{0}{0}    
    \numbercircle{G2}{(2.1,0.5)}{green!60!black}{15}{1}
    \numbercircle{G3}{(4.5,0.5)}{green!60!black}{0}{0}    
    \numbercircle{G4}{(5.3,0.5)}{green!60!black}{15}{1}
    \numbercircle{G5}{(3,1.2)}{green!60!black}{0}{0}    
    \numbercircle{G6}{(3,-0.2)}{green!60!black}{30}{2}
    \numbercircle{G7}{(2.7,1.9)}{green!60!black}{15}{1}    
    \numbercircle{G8}{(2.7,-0.9)}{green!60!black}{30}{2}
    \numbercircle{G9}{(4.2,1.2)}{green!60!black}{30}{2}    
    \numbercircle{G10}{(4.2,-0.2)}{green!60!black}{15}{1}
    \numbercircle{G11}{(4.7,1.9)}{green!60!black}{30}{2}    
    \numbercircle{G12}{(4.7,-0.9)}{green!60!black}{0}{0}
    \draw[->, thick] (G1) -- (G5);
    \draw[->, thick] (G6) -- (G1);
    \draw[->, thick] (G2) -- (G7);
    \draw[->, thick] (G8) -- (G2);
    \draw[->, thick] (G3) -- (G10);
    \draw[->, thick] (G9) -- (G3);
    \draw[->, thick] (G4) -- (G12);
    \draw[->, thick] (G11) -- (G4);
    \draw[->, thick] (G5) to node[above=-1pt,  sloped, pos=0.6] {\scriptsize$5$} (G11);
    \draw[->, thick] (G5) to[out=-30, in=-150] node[below=-1pt,  sloped] {\scriptsize$1$} (G9);
    \draw[->, thick] (G7) to[out=30, in=150] node[below=-1pt,  sloped] {\scriptsize$3$} (G11);
    \draw[->, thick] (G7) to node[above=-1pt,  sloped, pos=0.4] {\scriptsize$7$} (G9);
    \draw[->, thick] (G10) to[out=150, in=30] node[above=-1pt,  sloped, pos=0.5] {\scriptsize$8$} (G6);
    \draw[->, thick] (G10) to node[below=-1pt,  sloped, pos=0.6] {\scriptsize$2$} (G8);
    \draw[->, thick] (G12) to node[below=-1pt,  sloped, pos=0.4] {\scriptsize$4$} (G6);
    \draw[->, thick] (G12) to [out=-150, in=-30] node[above=-1pt,  sloped] {\scriptsize$6$} (G8);

    \numbercircle{B1}{(6.7,0.5)}{blue}{15}{1}
    \numbercircle{B6}{(7.3,0.5)}{blue}{15}{1}
    \numbercircle{B2}{(6.2,1.2)}{blue}{0}{0}
    \numbercircle{B5}{(7.8,1.2)}{blue}{30}{2}
    \numbercircle{B3}{(6.6,1.9)}{blue}{30}{2}
    \numbercircle{B4}{(7.4,1.9)}{blue}{0}{0}
    \numbercircle{B10}{(6.2,-0.2)}{blue}{15}{1}
    \numbercircle{B7}{(7.8,-0.2)}{blue}{0}{0}
    \numbercircle{B9}{(6.6,-0.9)}{blue}{0}{0}
    \numbercircle{B8}{(7.4,-0.9)}{blue}{30}{2}
    \draw[->, thick] (B1) -- (B2);
    \draw[->, thick] (B2) -- (B3);
    \draw[->, thick] (B3) -- (B4);
    \draw[->, thick] (B4) -- (B5);
    \draw[->, thick] (B5) -- (B6);
    \draw[->, thick] (B6) -- (B7);
    \draw[->, thick] (B7) -- (B8);
    \draw[->, thick] (B8) -- (B9);
    \draw[->, thick] (B9) -- (B10);
    \draw[->, thick] (B10) -- (B1);

    \node[] at (-2,-1.7) {\textbf{E}};
    \node[] at (0.5,-1.7) {\textbf{C}};
    \node[] at (3.7,-1.7) {\textbf{A}};
    \node[] at (7.0,-1.7) {\textbf{B}};
\end{tikzpicture}
\caption{Examples of periodic orbits in $4$ groups.}
\label{fig:4examples}
\end{figure}
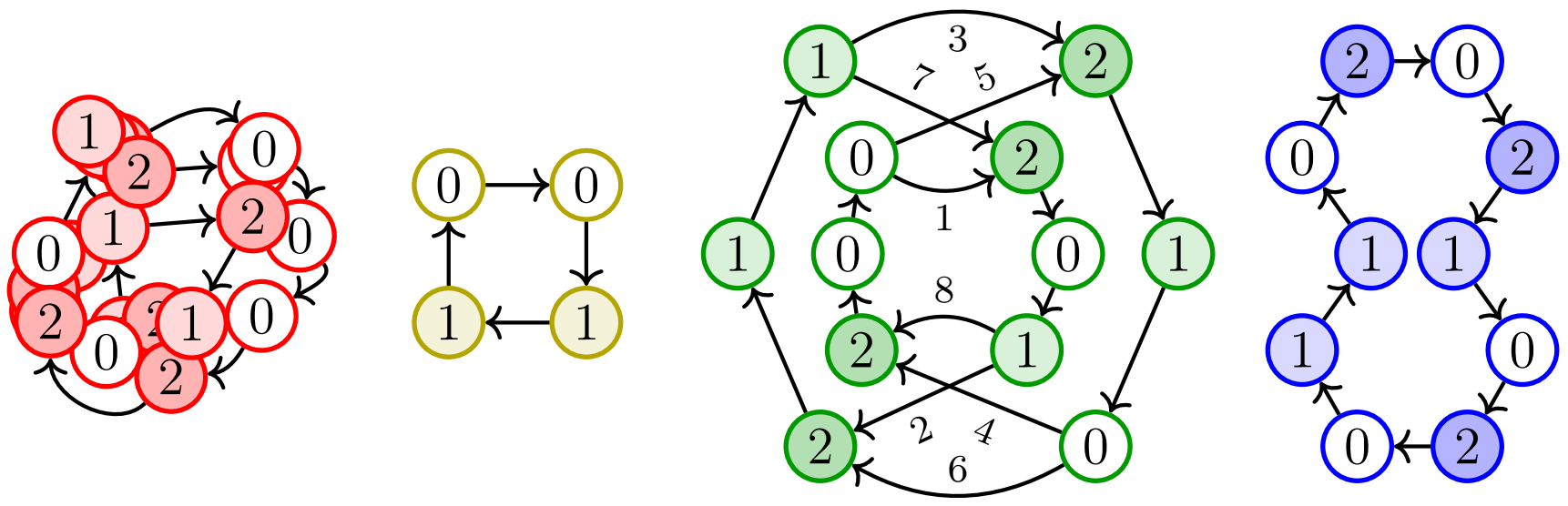

\subsubsection{Group $A$ -- medium length orbits -- Section~\ref{sec:A}} An example of a Group A orbit is given in Figure~\ref{fig:4examples}. The numbered edges should be traversed in the increasing order. So, the orbit has period $24$ and goes through paths $(200)$, $(201)$, $(210)$, and $(211)$ twice during the period. The main idea of the section is to show that if this orbit is maximal, then other competitive orbits should visit these paths the same proportion of times (have the same concentration). So, if some other orbit of period $24$ visits $(201)$ thrice, then it would not be competitive. Moreover, the same principle applies to the sections of the maximal orbit. If one section of the periodic orbit has different concentrations, compared to some other section, then this orbit is likely not maximal. Using combinatorial estimates we prove that such uniform-concentration orbits cannot be that long, thus eliminating all Group A orbits. Group A consists of orbits $\boo$ with $\Delta(\boo) > a_n$ and
\begin{equation*}
    \exp \bigl( n^{1/10}\bigr) \le \abs{\boo} \le \theta^{-n/3}.
\end{equation*}
The uniform-concentration principle holds due to the similar $\delta$-argument as for Group~E. If one half of the orbit frequently visits $(200)$ and the second -- $(211)$, then the ergodic supremum of the first part will be related to the value of $f$ on words containing $(200)$, opposed to $(211)$ for the second part. With high probability, these values will not be close to one another. Say $(200)$-words will have higher value. Then, if one closes the first half of the orbit into its own orbit, one would get an orbit with a larger ergodic supremum than that of the original. Thus orbits without uniform concentration are likely not maximal. A more precise argument is given in Section~\ref{sec:A}; it also uses a $\delta$-argument and requires small entropy reduction of Section~\ref{sec:thetareduction}.

\subsubsection{Group $B$ -- short orbits with a small gap -- Section~\ref{sec:B}}
When an orbit has a small gap, it can be split into two or more shorter orbits that on average give a similar ergodic average compared to the whole orbit, see Lemma~\ref{lemma.gapcompose}. For example, Group B orbit on Figure~\ref{fig:4examples} repeats $(1020)$ twice, so it can be split into $(10202)^{\N}$ and $(10201)^{\N}$. So, if such an orbit is maximal, it can be cut in half to obtain an almost-maximal orbit. We divide the set of periodic orbits with $\Delta(\boo) >a_n$ into many levels, with the $i$-th level consisting of orbits with $a_{k_i-1} < \Delta(\boo) \le a_{k_i}$ for $k_i \= n-is$ where $1\le i \le R$, and the $(R+1)$-th level consisting of the remaining large gap orbits.

For an orbit in level $i$, we use Lemma~\ref{lemma.gapcompose} to cut it up and obtain ``suborbits'' in level $\ge i+1$. For the maximal orbit $\widehat{\boo}$ among those with $\Delta(\boo) > a_n$, we recursively cut $\widehat{\boo}$ until we get an orbit in level $R+1$. On one hand, with Group~A eliminated, short length of $\widehat{\boo}$ must lead to a jump in level, i.e., on one step the level is $i$ and on the next step the level is $\ge i+2$. On the other hand, using modified concentration arguments from Section~\ref{sec:A}, we prove with high probability that at a jump of level, the orbit at this step has the same concentration as $\widehat{\boo}$, leading to a contradiction. In conclusion, we eliminate the case $\Delta(\widehat{\boo}) \le a_{n-s}$. This idea of recursively cutting the orbits follows from Proposition~2.1 of \cite{HLMXZ25}, see also \cite{Boc19}. The concentration arguments are modified as follows: at each length of words $k$, we use a $l^1$-norm ball as a perturbation space instead of a $l^\infty$-norm ball used in Sections~\ref{sec:thetareduction} and~\ref{sec:A}. See the remark below Corollary~\ref{cor:concB} for more details.

\subsubsection{Group $C$ -- short orbits with a large gap -- Section~\ref{sec:C}} This group in particular contains some very simple orbits, like $(01)^{\N}$ or $(0011)^{\N}$. If an orbit like $(01)^{\N}$ is the maximal one in our finite set, that means that it gives the best approximation for $\mpe(\Sigma_{\bF},f)$ in this set. But if $\mpe\bigl((01)^{\N}, f\bigr) < \mpe(\Sigma_{\bF},f)$, then it cannot be the best approximation for arbitrary large $n$. Thus we should either have the locking property with $\mpe\bigl((01)^{\N}, f\bigr)= \mpe(\Sigma_{\bF},f)$ or expect $\mpe\bigl((01)^{\N}, f\bigr)$ to be away from the maximal, making these orbits quite volatile.

To control these volatile orbits we basically use the $\delta$-argument developed in \cite{DLZ24} with a similar modification to that in Section~\ref{sec:B}: using $l^1$ norm instead of $l^\infty$ norm. With Groups~E, A, and B eliminated, the probability estimation improves a lot  compared with that in \cite{DLZ24},  ultimately yielding the locking property with large probability.

\bigskip

Regarding other sections, Section~\ref{sec:prelim} introduces notions that are used throughout the proof, while Section~\ref{sec:advantage} introduces normalized ergodic supremum and gives the locking criterion in terms of periodic orbits. Ideas from it are used to eliminate Groups A, B, and C. Finally, Section~\ref{sec:split} ties together arguments from Sections~\ref{sec:A}--\ref{sec:C} to prove Theorem~\ref{thm:main} for small entropy.

\subsection*{Acknowledgments}
Z.~H.\ and Z.~L.\ were partially supported by Beijing Natural Science Foundation (JQ25001) and National Natural Science Foundation of China (12471083). 

I.~K.\ expresses deep gratitude to Peking University for the hospitality during his month-long visit in April 2026, during which the essential progress on this collaboration was made, following the collaborative work earlier. He is also thankful to Chongqing Normal University, Jilin University, Soochow University, and Peking University for their warm hospitality during his subsequent visits in May and June 2026, when he made further progress on this project together with Z.~H.\ and Z.~L.
The authors are grateful to Vadim Kaloshin for introducing I.~K.\ to the other authors and for his encouragement throughout the project. 

\subsection*{AI statement}
AI systems were used in this research as detailed below:

During the initial exploratory phase of this research, some of the authors utilized OpenAI’s ChatGPT~5.5 for brainstorming, by providing AI our own potential approaches and preliminary findings such as Proposition~\ref{prop:pre.zero.entropy}. While we did not get any useful ideas from these human--AI discussions, regrettably, this process led to an unintended confidentiality leak: five days after the new model GPT~6.0 was released in September 2026, we noticed that our Proposition~\ref{prop:pre.zero.entropy} had been inadvertently exposed and reproduced in the model's outputs from an account that had never been given this proposition.\footnote{Similar concerns have also been raised by mathematicians in connection with OpenAI, e.g., on finite-time blow-up for the forced Navier--Stokes equations  (\url{https://cims.nyu.edu/~tristanb/statement.pdf}), non-sofic groups (\url{https://terrytao.wordpress.com/2026/09/11/on-the-existence-of-non-sofic-groups/}), and the abundance conjecture (\url{https://arxiv.org/pdf/2610.35913v1}, Remark~0.1).} 
Though this proposition is only an elementary component in our proof, we were concerned by this apparent leak and terminated all use of AI-assisted tools for the remainder of this project.

Following discussions with the other authors between April and August 2026 on this project, one of the authors, who worked independently of AI throughout this project (except for literature searching), developed the central ideas and the overall structure of the main proof without drawing on the exploratory AI discussions described above. AI served only as an indirect aid on two isolated technical choices in the implementation of this approach. First, the geometric sequence construction of a sequence $\{n_i\}$, used to help estimate the length of the orbit in Subsection~\ref{subs:counting}, was inspired by discussions with AI on some related combinatorial counting problem. Second, the choice of the sequence $\{k_i\}$, which is used only in the half-cutting argument in Proposition~\ref{prop:halfcut}, was inspired by an AI-assisted summary and reformulation of Bochi's interpretation in \cite{Boc19}.

\section{Preliminaries}
\label{sec:prelim}

The set of positive integers is denoted by $\N$, and the set of nonnegative integers by $\N_0 \= \N \cup \{0\}$. The cardinality of a set $S$ is denoted by $\abs{S}$. For a finite or infinite sequence $\bx =\{ x_i\}_{i\in I}$ ($I$ is a finite or countable set) of real numbers, we denote its $l^\infty$-norm by $\norm{\bx}_{\infty} \= \sup_{i\in I}\abs{x_i}$, its $l^1$-norm by $\norm{\bx}_{1} \= \sum_{i\in I} \abs{x_i}$, and its $l^2$-norm by $\norm{\bx}_{2} \= \bigl(\sum_{i\in I} x_i^2 \bigr)^{1/2}$. Throughout the proof, we will assume that $n_{\theta} \in \N$ is a sufficiently large number, in particular all forbidden words $\bw \in \bF$ have $\abs{\bw} < n_{\theta}/2$.

\subsection{Subshifts and subshifts of finite type}

All shifts and subshifts in this article are assumed to be one-sided.
For each $n \in \N$, the set $\Sigma_{n} \= \bQ^{n}$ of \emph{words of length $n$} (or \emph{$n$-words}) in the alphabet $\bQ$ is the set of $q$-nary sequences of length $n$. Write $\Sigma_{0} \=\{\emptyset \}$ and call $\emptyset$ the word of length $0$. Define $\abs{\bw} \=n$ for each word of length $n\in\N_0$.  Denote the set of \emph{finite words} by $\Sigma_* \= \bigcup_{n=0 }^{+\infty}\Sigma_{n}$ and the set of words of length at most $k$, for $k\in\N_0$, by $\Sigma_{\leq k} \= \bigcup_{n = 0}^k \Sigma_{n}$. Note that $\Sigma_*$ includes $\emptyset$, i.e., the word of length $0$.

The set $\Sigma^+ \= \bQ^{\N}$ of \emph{infinite words} in the alphabet $\bQ$ is the set of one-sided infinite $q$-nary sequences. The system $(\Sigma^+, \sigma)$ is called a \emph{full shift} on $q$ symbols. For a compact subspace $X \subseteq \Sigma^+$, $(X,\sigma)$ (or simply, $X$) is called a \emph{subshift} if $\sigma(X)=X$.

For a subshift $X$ and $n \in \N$, we denote the set of length-$n$ allowed words by 
\begin{equation}\label{eq:def.nallowed}
    \Sigma_n(X) \= \{ \bz_{[1,n]} : \bz\in X \},
\end{equation}
where $\bz_{[1,n]} \in \Sigma_n$ denotes the $n$-prefix of $\bz$, and we write $\Sigma_n^c(X) \= \Sigma_n \smallsetminus \Sigma_n(X)$ for its complement in $\Sigma_n$.  We set $\Sigma_0(X) \= \{ \emptyset\}$.
We denote the \emph{$m$-cylinder cover} of $X$ by
\begin{equation*}
    [\Sigma_m(X)] \= \bigcup_{\bw \in \Sigma_m(X)} [\bw]
\end{equation*}
and define $[\Sigma_m^c(X)]$ similarly.
Moreover, we define $\Sigma_*(X) \= \bigcup_{n=0 }^{+\infty} \Sigma_n (X) $ and $\Sigma_{\le k}(X) \= \bigcup_{n =0}^{k} \Sigma_n (X)$ for each $k\in \N_0$. We often use these notations in the case where $X$ is some periodic orbit $\boo$.

Let $\bF \subseteq \Sigma_*$ be a finite subset of words, called \emph{forbidden words}. We often write $\{w_n\}_{n\in\N} = w_1 w_2 \dots w_n \dots$ for an infinite word $\bw= \{w_n\}_{n\in\N} \in \Sigma^+$, and similarly, write $\{z_n\}_{n=1}^m = z_1 z_2 \dots z_m$ for an $m$-word $\bz= \{z_n\}_{n=1}^m \in \Sigma_m$, $m\in\N$. The \emph{subshift of finite type determined by $\bF$} is defined as
\begin{equation}\label{eq:def.SFT}
    \Sigma_{\bF} \= \{ \bx\in \Sigma^+ : x_i x_{i+1} \dots x_{i+n-1} \notin \bF \text{ for all } i,\, n \in \N \}.
\end{equation}
A subshift $X$ is called a \emph{subshift of finite type} if $X$ is the subshift of finite type determined by some finite subset of words $\bF$. If $X =\Sigma_{\bF}$ is a subshift of finite type, for each $n\in \N_0$, we denote $\Sigma_n(\bF) \= \Sigma_n(X)$, $\Sigma_{\le n}(\bF) \= \Sigma_{\le n}(X)$, and $\Sigma_*(\bF) \= \Sigma_*(X)$ for simplicity. For every $n\in \N$, we denote
\begin{equation}\label{eq:deffk}
	Z_n \= \abs{\Sigma_{n}(\bF)}.
\end{equation}

Let $X \subseteq \Sigma^+$ be a subshift and $n\in \N$. For $\bw \in \Sigma_n(X)$, we define the \emph{$\bw$-cylinder} by $[\bw] \=  \{ \bx \in X : \bx_{[1,n]} = \bw \}$ and call $[\bw]$ a length-$n$ cylinder of $X$. For a function $f \: X \to \R$ and $n\in \N$, we recall the \emph{$n$-th variation} of $f$ given by 
\begin{equation*}
    \var_n(f) \= \max_{ \bw\in \Sigma_n(X)} \Bigl\{ \sup_{\bx\in [\bw]} f(\bx) - \inf_{\bx \in [\bw] }f(\bx) \Bigr\} .
\end{equation*}
Under metric $d_{\ba}$, an immediate bound on $\var_n(f)$ is given by the minimal Lipschitz constant:
\begin{equation}    \label{eq:varlip}
    \var_n(f) \le \LIP_{d_{\ba}}(f)\cdot a_{n+1}.
\end{equation}

\subsection{Hilbert bricks and admissible gauges}

We recall the construction of Hilbert bricks and admissible gauges (cf.~\cite[Section~2]{BZ16} and \cite[Lemma~2.1]{DLZ24}).
In the present setting, we use cylinder indicators in the defining series.
Let $\bb=\{b_{\bw}\}_{\bw \in \Sigma_*(\bF)}$ be a family of nonnegative numbers such that
\begin{equation} \label{eq:def.bn}
    \sum_{n=1}^{+\infty}b_n<+\infty,\quad
    \text{ where } b_n\= \max_{\bw\in \Sigma_n(\bF)} b_{\bw}.
\end{equation}
For a family of real numbers $\be=\{e_{\bw}\}_{\bw\in \Sigma_*(\bF)}$ satisfying $\abs{e_{\bw}}\leq b_{\bw}$ for each
$\bw\in\Sigma_*(\bF)$, the series $\sum_{\bw \in \Sigma_*(\bF)} e_{\bw} \mathbbold{1}_{[\bw]}$ converges uniformly on $\Sigma_{\bF}$ to a continuous function $g_{\be}\:\Sigma_{\bF}\to\R$, since for each $n\in\N$, at most one cylinder of length $n$ contains any given point.
Define a \emph{Hilbert brick} $\cH[\bb]$ to be the set of the following functions (cf.~\cite[Lemma~2.1]{DLZ24}):
\begin{equation*}
    \cH[\bb] \= \{ g_{\be}: \be=\{ e_{\bw} \}_{\bw \in \Sigma_*(\bF)} \text{ with }  \abs{e_{\bw}} \leq b_{\bw} \text{ for each }\bw \in \Sigma_*(\bF)\}.
\end{equation*}

We say that $\bb$ is an \emph{$\ba$-admissible gauge} if $\cH[\bb]$ is compact in $\Lip(\Sigma_{\bF}, d_{\ba})$. By Lemma~\ref{lem:hilbert-brick}, a sufficient condition for a gauge $\bb$ to be $\ba$-admissible is 
\begin{equation}\label{eq:a-admissible}
    \lim_{n\to +\infty} b_n/a_n = \lim_{n\to +\infty} \max_{\bw\in \Sigma_n(\bF)} b_{\bw} /a_n =0,
\end{equation}
For an $\ba$-admissible gauge $\bb$, let $\cI[\bb]\=\prod_{\bw\in\Sigma_*(\bF)}[-b_{\bw},b_{\bw}]$ be equipped with the product topology. For each $\bw \in \Sigma_*(\bF)$, let $\nu_{\bw}$ be a probability measure supported on $[-b_{\bw}, b_{\bw}]$.
The product measure $\overline\nu \= \otimes_{\bw \in \Sigma_*(\bF)} \nu_{\bw}$ defines a probability measure on $\cI[\bb]$.

Define the projection $\Pi\:\cI[\bb]\to\cH[\bb]$ by
\begin{equation}  \label{eq:def_Pi}
    \Pi(\be)
    \=g_{\be} 
    = \sum_{\bw \in \Sigma_*(\bF)} e_{\bw} \mathbbold{1}_{[\bw]}.
\end{equation}
The continuity of $\Pi$ follows from the uniform
convergence of the defining series and the compactness of $\cH[\bb]$ in the Lipschitz norm; see  Lemma~\ref{lem:hilbert-brick}. Hence the pushforward measure $\nu\=\Pi_{*}\overline\nu$
is a Borel probability measure on $\cH[\bb]$.
We call any such measure $\nu$ an \emph{$\ba$-admissible measure}.

For the quantitative estimates throughout the proof, we introduce the following notion of $\log$-admissibility.

\begin{definition} \label{def:logadmis}
An $\ba$-admissible gauge $\bb$ is \emph{$C$-$\log$-$\ba$-admissible} for some $C \in (0, 1)$ if
\begin{equation}    \label{eq:logadmis}
    b_{\bw}  \begin{cases}
        = 0 &\text{ if }\abs{\bw} < n_{\theta},\\
        \le Ca_{\abs{\bw}} / \log \abs{\bw}&\text{ if }\abs{\bw} \ge n_\theta.
    \end{cases}
\end{equation}
\vspace{-\baselineskip}
\end{definition}
We call a Hilbert brick (resp.\ an admissible measure) \emph{$C$-$\log$-$\ba$-admissible} if it is induced by a $C$-$\log$-$\ba$-admissible gauge as above. Note that all $g$ in a $\log$-admissible Hilbert brick $\cH[\bb]$ satisfy $\LIP(g) < 1/4$, provided that $n_{\theta}$ is large enough.

If $\nu'$ and $\nu''$ are respectively $C'$ and $C''$-$\log$-$\ba$-admissible measures (perhaps with different $\bb$), we can consider their convolution $\nu \= \nu' \ast \nu''$. Then, $\nu$ is a $(C' + C'')$-$\log$-$\ba$-admissible measure, see Proposition~\ref{prop:finite-convolution}. Moreover, if the elements of the sequence $\{\nu_i\}_{i\in \N}$ are respectively $C_i$-$\log$-$\ba$-admissible and $C_* \= \sum_{i=1}^{+\infty} C_i <  +\infty$, then 
$    \nu\=\nu_1\ast\nu_2\ast\nu_3\ast\cdots$
is a $C_*$-$\log$-$\ba$-admissible measure, see Proposition~\ref{prop:infinite-convolution}. This is useful for the following: let $\nu_A$ and $\nu_B$ be two $\ba$-admissible measures. Assume that we have independent random variables $g_A\sim \nu_A$ and $g_B \sim \nu_B$, and assume that for every $f_0$ the functions $f_0+g_A$ and $f_0 +g_B$ satisfy properties $A$ and $B$ respectively with high probability. Then, for $\nu \= \nu_A\ast \nu_B$ and $g\sim \nu$ and for all $f_0$, $f_0+g$ satisfies both $A$ and $B$ with high probability.

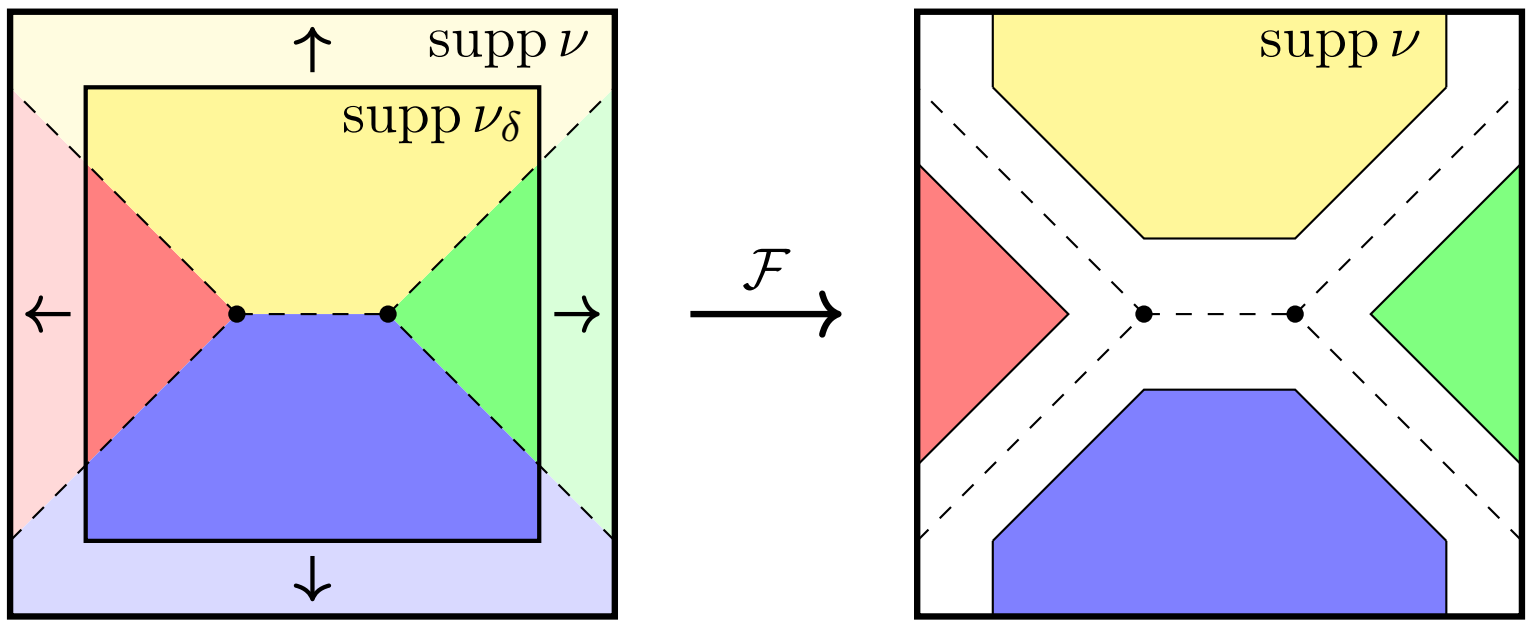
\begin{figure}
\begin{tikzpicture}
    \fill[red!50]
        (0.5, 1) -- (1.5, 2) -- (0.5, 3) --  cycle;
    \fill[red!50]
        (6, 1) -- (7, 2) -- (6, 3) --  cycle;
    \fill[green!50]
        (3.5, 1) -- (2.5, 2) -- (3.5, 3) --  cycle;
    \fill[green!50]
        (10, 1) -- (9, 2) -- (10, 3) --  cycle;
    \fill[yellow!50]
        (3.5, 3.5) -- (3.5, 3) -- (2.5, 2) -- (1.5, 2) -- (0.5, 3) -- (0.5, 3.5) --  cycle;
    \fill[yellow!50]
        (9.5, 4) -- (9.5, 3.5) -- (8.5, 2.5) -- (7.5, 2.5) -- (6.5, 3.5) -- (6.5, 4) --  cycle;
    \fill[blue!50]
        (3.5, 0.5) -- (3.5, 1) -- (2.5, 2) -- (1.5, 2) -- (0.5, 1) -- (0.5, 0.5) --  cycle;
    \fill[blue!50]
        (9.5, 0) -- (9.5, 0.5) -- (8.5, 1.5) -- (7.5, 1.5) -- (6.5, 0.5) -- (6.5, 0) --  cycle;
    \fill[red!15]
        (0, 0.5) -- (0.5, 1) -- (0.5, 3) -- (0, 3.5) -- cycle;
    \fill[green!15]
        (4, 0.5) -- (3.5, 1) -- (3.5, 3) -- (4, 3.5) -- cycle;
    \fill[yellow!15]
        (0, 4) -- (0, 3.5) -- (0.5, 3) -- (0.5, 3.5) -- (3.5, 3.5) -- (3.5, 3) -- (4, 3.5) -- (4, 4) -- cycle;
    \fill[blue!15]
        (0, 0) -- (0, 0.5) -- (0.5, 1) -- (0.5, 0.5) -- (3.5, 0.5) -- (3.5, 1) -- (4, 0.5) -- (4, 0) -- cycle;

    \draw[->, black, thick] (2, 3.6) -- (2, 3.9);
    \draw[->, black, thick] (2, 0.4) -- (2, 0.1);
    \draw[->, black, thick] (3.6, 2) -- (3.9, 2);
    \draw[->, black, thick] (0.4, 2) -- (0.1, 2);

    \draw[->, black, very thick] (4.5, 2) -- (5.5, 2);

    \draw[black, thick] (0.5,0.5) -- (3.5,0.5) -- (3.5,3.5) -- (0.5, 3.5) -- cycle;

    \draw[black, dashed] (0, 0.5) -- (1.5, 2) -- (0, 3.5);
    \draw[black, dashed] (4, 0.5) -- (2.5, 2) -- (4, 3.5);
    \draw[black, dashed] (1.5, 2) -- (2.5, 2);
    \draw[black, dashed] (6, 0.5) -- (7.5, 2) -- (6, 3.5);
    \draw[black, dashed] (10, 0.5) -- (8.5, 2) -- (10, 3.5);
    \draw[black, dashed] (7.5, 2) -- (8.5, 2);
    \draw[black] (6.5, 3.5) -- (6.5, 4);
    \draw[black] (9.5, 3.5) -- (9.5, 4);
    \draw[black] (6.5, 0.5) -- (6.5, 0);
    \draw[black] (9.5, 0.5) -- (9.5, 0);
    \draw[black] (6.5, 3.5) -- (7.5, 2.5) -- (8.5, 2.5) -- (9.5, 3.5);
    \draw[black] (6.5, 0.5) -- (7.5, 1.5) -- (8.5, 1.5) -- (9.5, 0.5);
    \draw[black] (6, 1) -- (7, 2) -- (6, 3);
    \draw[black] (10, 1) -- (9, 2) -- (10, 3);
    
    \draw[black, very thick] (0,0) -- (4,0) -- (4,4) -- (0, 4) -- cycle;
    \draw[black, very thick] (6, 0) -- (10,0) -- (10,4) -- (6, 4) -- cycle;

    \draw[fill=black] (1.5, 2) circle (0.05);
    \draw[fill=black] (2.5, 2) circle (0.05);
    \draw[fill=black] (7.5, 2) circle (0.05);
    \draw[fill=black] (8.5, 2) circle (0.05);
    \node[] at (3.3,3.75) {$\supp \nu$};
    \node[] at (8.8,3.75) {$\supp \nu$};
    \node[] at (2.8,3.25) {$\supp \nu_{\delta}$};
    \node[] at (5,2.3) {$\cF$};
\end{tikzpicture}
\caption{$\delta$-argument applied to $4$ functions $\mpe_1 = -x-1$, $\mpe_2 = y$, $\mpe_3 = x-1$, and $\mpe_4 = -y$. It shows that with high probability in an $(x, y)$-square one of these functions is quantitatively larger than the rest. Here, both $\nu$ and $\nu_\delta$ are uniform measures.}
\label{fig:deltaargument}
\end{figure}

\subsection{$\delta$-argument}
\label{subsec:darg}
Now we describe the $\delta$-argument that helps connect the problem of prevalence to the problem of density. An example of such an argument is shown in Figure~\ref{fig:deltaargument}. Let us say that we have an $\ba$-admissible measure $\nu$ and the goal is to show that for $g \sim \nu$ for all $f$ the function $f +g$ likely satisfies some property~$X$.
To do this, we can find a cut-off measure $\nu_\delta \le \nu$ with $\int\! \, \mathrm{d} \nu_\delta > 1-\delta$ for some small $\delta>0$. Then, if for all $f$ we find a Borel map $\cF\:\supp \nu_{\delta} \rightarrow \supp \nu$ with $\cF^*(\nu_{\delta}) \le \nu$ such that $f+\cF(g)$ always has property $X$, then $f+g$ satisfies $X$ with probability at least $1-\delta$.

Once we know that property~$X$ holds for a dense set of functions, we would find many near-identity $\cF$ so that $f+\cF(g)$ always has $X$. Then, the problem of prevalence reduces to constructing $\cF$ in an injective and measure-respecting way.

\subsection{Truncations of $f$}
\label{subsec:AF}
To state the locking criterion we need to define a truncation of $f$ that makes it constant on each length-$m$ cylinder, for a given $m \in \N_0$. In \cite{BZ16} this is done by averaging $f$ over the cylinder. For a subshift of finite type we use a similar construction for a truncated $A_mf$. Let $\mu$ be a full-supported probability measure on $\Sigma_F$. Then, for every $\bw \in \Sigma_{\bF}$, define
\begin{equation*}
    A_mf (\bw')\= \frac{1}{\mu([\bw_{[1,m]}])}\int_{[\bw_{[1,m]}]}\! f \,\mathrm{d}\mu.
\end{equation*}
In particular, $A_m$ is a linear operator. The truncated $A_mf$ inherits some properties of $f$, also see \cite{BZ16}. A property that follows from \eqref{eq:varlip} is that
\begin{equation}\label{eq:fAmfclose}
    \norm{f-A_mf}_{\infty} \le \var_m(f) \le  \LIP(f)\cdot a_{m+1}. 
\end{equation}
We also have
\begin{equation}    \label{eq:truncnormbound}
    \norm{A_mf}_{\infty} \le \norm{f}_{\infty}, \quad
    \LIP(A_m f) \le \LIP(f), \quad \text{ and }\quad\ A_m A_n f = A_n f 
\end{equation}
if $1\le m \le n$.

\subsection{Periodic orbits and paths}
The shift $\sigma$ has many periodic orbits. Define $\PO_*$ to be the set of all periodic orbits and let $\PO_n$ (resp.\ $\PO_{\le n}$) denote the set of periodic orbits with period equal to (resp.\ no greater than) $n\in\N$. 

Let $\boo \in \PO_*$ be a periodic orbit. Denote its \emph{period} by $\abs{\boo}$ and let $\mu_{\boo}$ be the unique invariant measure supported on it. Then, recall $\mpe(\boo, f) = \int\! f\, \mathrm{d}\mu_{\boo}$ and call it the \emph{score} of $\boo$. Let $\MMM_n \= \{\mu_{\boo}: \boo \in \PO_n\}$ and $\MMM_{\le n} \= \{\mu_{\boo}: \boo \in \PO_{\le n}\}$ be sets of measures supported on periodic orbits. 

We remind that the gap $\Delta(\boo)$ of the orbit was defined in \eqref{eq:gap}. Next, let $n_{\Delta}(\boo) \in \mathbb{N}$ be such that $a_{n_{\Delta}(\,\cdot\,)} = \Delta(\,\cdot\,)$. Let $\cO_{\Delta,n}$ denote the set of periodic orbits $\boo$ with gap $\Delta(\boo) > a_n$.  The gap indicates how composite the orbit $\boo$ is and controls how well $\boo$ can be approximated by the shorter orbits.

\begin{lemma}
    Let $\Sigma_{\bF}$ be a subshift of finite type, $\boo \in \cO_*$, $m\in \N$ with $m  \ge \sup_{\bw \in \bF}\abs{\bw}$, and $f \in \Lip(\Sigma_{\bF}, d_{\ba})$. If $\Delta(\boo) \le a_m$, then there exists a finite sequence of periodic orbits $\{\boo_i\}_{i=1}^s$ with gaps larger than $a_m$ and $\sum_{i=1}^s \abs{\boo_i} = \abs{\boo}$ such that
    \begin{align}
        \Absbigg{ \abs{\boo} \mpe(\boo, f) - \sum_{i=1}^s \abs{\boo_i} \mpe(\boo_i, f) } &\le 2(s-1) \LIP(f) \frac{\theta a_m}{1-\theta} \quad\text{ and }  \label{eq:shadowing1} \\
        \mu_{\boo} ([\bw]) &= \sum_{i=1}^s \frac{\abs{\boo_i}}{ \abs{\boo}} \mu_{\boo_i} ([\bw]) \label{eq:shadowing4}
    \end{align}
    for every $\bw\in \Sigma_m(\bF)$.
    Moreover, if $f= A_m f$, then 
    \begin{equation}\label{eq:shadowing2}
        \abs{\boo} \mpe(\boo, f) = \sum_{i=1}^s \abs{\boo_i} \mpe(\boo_i, f).
    \end{equation}
    \label{lemma.gapcompose}
    \vspace{-\baselineskip}
\end{lemma}

\begin{proof}
    Write $\boo = \{x_i \}_{i=1}^k$ with $x_i = \sigma^{i-1}(x_1)$ and $x_1 = (\alpha_1 \alpha_2 \dots \alpha_k)^\infty$. The indices of $\{x_i\}_{i=1}^k$ and $\{ \alpha_i \}_{i=1}^k$ are taken modulo $k$.
    Let $2\le j \le k$ be the minimal integer such that there exists $1\le i < j$ such that $d_{\ba}(x_i, x_j) \le a_m$. Since $\Delta(\boo) \le a_m$, such $j$ always exists. Now we set $y_1 \= (\alpha_i \dots \alpha_{j-1})^\infty$ and $z_1 \= (\alpha_j \dots \alpha_k \alpha_1 \dots \alpha_{i-1})^\infty$. Let $\boo_1$ and $\boo_1'$ be the orbit of $y_1$ and $z_1$, respectively. Since the minimal period divides any period,
    set $p_1\= (j-i)/\abs{\boo_1}\in \N$ and $q_1\=(k-j+i)/\abs{\boo_1'}\in \N$.
    Note that $d_{\ba}(x_i, x_j) \le a_m$ implies $\alpha_{i+s} = \alpha_{j+s}$ for all $s \in\{ 0,\, 1,\, \ldots,\, m-2\}$. So, every subword of $y_1$ or $z_1$ of length at most $m$ is a subword of $x_1$ and is thus allowed. Hence, $y_1$, $z_1 \in \Sigma_{\bF}$ are well defined.
    Thus, by (\ref{eq:fAmfclose}), we estimate
    \begin{equation}\label{eq:est.shadow}
        \begin{aligned}
            &\abs{ \abs{\boo} \mpe(\boo, f) - p_1\abs{\boo_1} \mpe(\boo_1, f) - q_1 \abs{\boo_1'} \mpe(\boo_1', f) } \\ &\qquad\le \sum_{l=0}^{j-i-1} \Absbig{ f\bigl(\sigma^l(y_1)\bigr) - f(x_{i+l} ) } + \sum_{l=0}^{k+i-j-1} \Absbig{ f\bigl(\sigma^l(z_1)\bigr) - f(x_{j+l} ) } \\
            &\qquad \le  \sum_{l=0}^{j-i-1} \var_{m+l}(f) +   \sum_{l=0}^{k+i-j-1} \var_{m+l}(f) \\
            &\qquad \le \LIP(f) \biggl( \sum_{l=0}^{j-i-1} a_{m+1+l} + \sum_{l=0}^{k+i-j-1} a_{m+1+l} \biggr) \\
            &\qquad \le 2\LIP(f) \frac{\theta a_m}{1-\theta}.
        \end{aligned}
    \end{equation}

    By our construction and the minimality of $j$, $\abs{\boo} = p_1\abs{\boo_1} + q_1\abs{\boo_1'}$ and $\Delta(\boo_1) > a_m$. 
    If $\Delta(\boo_1') \le a_m$, we do the same for $\boo_1'$ and find $\boo_2$ and $\boo_2'$ with $\abs{\boo_1'} = p_2\abs{\boo_2} + q_2\abs{\boo_2'}$, $p_2,\,q_2\in \N$, and $\Delta(\boo_2) > a_m$. We iterate this step until we obtain $\boo_{t-1}'$ such that $\Delta(\boo_{t-1}') > a_m$, where at each step every copy of the orbit being decomposed is replaced by the corresponding numbers of copies of the two resulting orbits. We then relabel the resulting sequence of periodic orbits as $\boo_1,\,\boo_2,\,\ldots,\boo_s$, with $\Delta(\boo_i)>a_m$ for each $i=1,\,2,\,\ldots,\,s$ and $\sum_{i=1}^s\abs{\boo_i}=\abs{\boo}$. Since each individual decomposition increases the number of orbits by at least one, (\ref{eq:shadowing1}) follows from applying (\ref{eq:est.shadow}) at most $s-1$ times. Note that at each step in the process above, the times all orbits at this step visiting $[\bw]$ is not changed for every $\bw \in \Sigma_m(\bF)$. Thus, (\ref{eq:shadowing4}) follows. Moreover, if $f=A_m f$, using the fact that $\var_n(f) =0$ for every integer $n \ge m$ and iterating (\ref{eq:est.shadow}), we obtain (\ref{eq:shadowing2}). This lemma is now proved.
\end{proof}

Let and $\boo \in \PO_*$ be a periodic orbit. For each $m \in \N$, the \emph{$m$-concentration} of $\boo$, $\conc_{m, \boo} \: \Sigma_m \rightarrow[0, 1]$, given by
\begin{equation*}
\conc_{m, \boo}(\bw_m)\= \frac1{\abs{\boo}}\Absbig{\bigl\{0\le i\le \abs{\boo}-1:\sigma^i(\bz)_{[1,m]}=\bw_m\bigr\}},
\end{equation*}
where $\bz$ is any point in $\boo$, measures how frequently the word $\bw_m$ appears in $\boo$. Concentrations are probability measures on $\Sigma_{m}$ and are thus associated to points in $\R^{q^m}$, or restricted to a subshift $X$, probability measures on $\Sigma_m(X)$ are associated to a point in $\R^{\abs{\Sigma_m(X)}}$.

Then, we introduce the paths.
Let $k\in\N_0$ and let $\bw$ and $\bw'$ be two allowed words of length $k$. For an integer $n>k$, a word $\bp\in\Sigma_{n+k}(\bF)$ that begins with $\bw$ and ends with $\bw'$ is called a \emph{path} of length $n$ from $\bw$ to $\bw'$.
Denote by $\Path_{n, \bw, \bw'}$ the set of all such paths and let $\Path_{n, k}$ be the union of these sets over all $\bw$ and $\bw'$ of length $k$. If a word $\bp \in \Sigma_{n+k}(\bF)$ has a subword $\bp' \in \Sigma_{n'+k'}(\bF)$, we call $\bp'\in \Path_{n', k'}$ a subpath of $\bp \in \Path_{n, k}$. 

For a path $\bp \in \Path_{n, k}$, with $k \ge n_{\theta}$, let $\overline\bp$ be an element of $\Sigma_{\bF}$ obtained by appending infinitely many symbols $\bp_{\ext} \in \Sigma_{\bF}$ to the end of $\bp$ to create a valid word in $\Sigma_{\bF}$. We choose $\bp_{\ext}(\bw')$ so that they only depend on the $n_\theta$-suffix of $\bw'$. Note that this is possible to do as this suffix is longer than any element of $\bF$, so the rest of $\bp$ and $\bp_{\ext}$ cannot interact to create a forbidden word.  Then, define the \emph{score} of $\bp$ and the measure corresponding to $\bp$ by
\begin{equation}    \label{eq:betapaths}
    \mpe(\bp, f) \= \frac{1}{n}\sum_{i = 0}^{n-1} f \bigl(\sigma^i(\overline \bp) \bigr) \quad \text{ and 
 } \quad \mu_{\bp} \= \frac{1}{n} \sum_{i=0}^{n-1} \delta_{\sigma^i(\overline{\bp})}.
\end{equation}
For $m \in\N$, define $\Sigma_m(\bp)$, $[\Sigma_m(\bp)]$, and $\conc_{m, \bp}$ using a similar idea to $\PO_*$: take all the words of length $m$ that are the prefix of some $\sigma^i(\overline \bp)$ for $0 \le i \le n-1$ and divide the number of such $i$-s by $n$. Let $\Sigma_m^\ast(\bp)$ and $\conc_{m, \bp}^\ast$ be the same objects, but when we take iterates $0 \le i \le n$ and divide by $n+1$ instead. Let 
\begin{equation*}
    \cou_{m, \bp} \= n \conc_{m,\bp} \quad \text{ and } \quad 
    \cou_{m, \bp}^{\ast} \= (n+1) \conc_{m, \bp}^{\ast}
\end{equation*}
be the undivided concentrations. When $m \le k+1$ all of $\Sigma_m(\bp)$, $\conc_{m, \bp}$, and $\cou_{m, \bp}$ are independent of the choice of $\bp_{\ext}$ (for $*$ objects we need $m\le k$) and thus do not require $k \ge n_{\theta}$ to define.

\subsection{Number of paths with the same concentration}
\begin{prop} \label{prop:mbp}
    Let $n$, $k$ be integers with $1\le k<n$, $\bw, \, \bw'\in \Sigma_k$, and $\bp \in \Path_{n, \bw, \bw'}$. Then 
    \begin{equation*}   
        M(\bp) 
        \= \abs{\{\bp' \in \Path_{n, \bw, \bw'} : \conc_{k+1, \bp'} = \conc_{k+1, \bp} \}}
        \le (n+k)^{\frac{n}{k}+1}.
    \end{equation*}
    \vspace{-\baselineskip}
\end{prop}
The proof is based on \cite{DLS81}. We also use Proposition~4.7 from \cite{KPRW17}, that itself is based on \cite{AB51}. Clearly, it is enough to consider a full shift with $\bF = \emptyset$. First, we define the graph $G_{\bp}$. The vertices of the graph are elements of $\Sigma_k^*(\bp)$. An edge is drawn between $\bw_0$ and $\bw_1$ if $\bw_0\ast = \circledast \bw_1 = \bw_2 \in \Sigma_{k+1}(\bp)$, where $\ast, \, \circledast \in \bQ$. This edge has multiplicity $\cou_{k+1, \bp}(\bw_2)$.  
\begin{lemma}[Proposition~4.7 in \cite{KPRW17}]
    \begin{equation}        \label{eq:detlapl}
        M(\bp) = \det \Delta_{\bw'}(G_{\bp}) \prod_{\bw^\ast \in \Sigma_k^\ast(\bp)} \frac{(\cou_{k, \bp}^*(\bw^\ast)-1)!}{\prod_{\alpha\in \bQ}(\cou_{k+1, \bp}(\bw^\ast \alpha))!},
    \end{equation}
    where $\Delta_{\bw'}(G_{\bp})$ is the graph Laplacian with the row and column indexed by $\bw'$ deleted.
\end{lemma}
\begin{proof}[Proof of Proposition~\ref{prop:mbp}]
    First, by Kirchhoff's theorem the determinant of the Laplacian minor measures the number of spanning oriented trees in $G_{\bp}$, rooted at $\bw'$. The number of rooted trees cannot exceed the product of the non-root vertices outdegrees $\cou_{k, \bp}^*(\bw^\ast)$, thus \eqref{eq:detlapl} leads to
    \begin{equation}        \label{eq:detlapl2}
        M(\bp) \le\prod_{\bw^\ast \in \Sigma_k(\bp)} \frac{(\cou_{k, \bp}(\bw^\ast))!}{\prod_{\alpha\in \bQ}(\cou_{k+1, \bp}(\bw^\ast \alpha))!}.
    \end{equation}
    Note that we changed $\cou^*$ into $\cou$ in the numerator as $\cou^*$ would have counted $\bw'$ one more time, but $\bw'$ is a root, so we do not multiply by its degree. We changed $\Sigma_k^\ast(\bp)$ into $\Sigma_k(\bp)$ due to the same reason.

    Then, in \eqref{eq:detlapl2} we get a product of multinomial coefficients. Consider an infinite periodic word $\bz \in [\bp]$ with $\sigma^{n+k}(\bz) = \bz$ obtained by glueing the ends of $\bp$ together. We will use $\cou_{k, \bz}$ to denote the number of visits of $\bz$ to respective $k$-cylinders along its first $n+k$ iterates, and let $\Sigma_k(\bz)$ be the set of visited $k$-words. We claim that
    \begin{equation}        \label{eq:detlapl3}
        M(\bp) \le\prod_{\bw^\ast \in \Sigma_k(\bz)} \frac{(\cou_{k, \bz}(\bw^\ast))!}{\prod_{\alpha\in \bQ}(\cou_{k+1, \bz}(\bw^\ast \alpha))!}.
    \end{equation}
    We have added $k$ new instances of words to both $\cou_{k}$ and $\cou_{k+1}$. However, every time we have added an instance of $\bw^\ast \alpha$ with $\alpha \in \bQ$ we have also added an instance to $\bw^\ast$, so the multinomial coefficient did not decrease, which proves the claim.
    
    Now, we bound the multinomial coefficients using the following:
    \begin{equation*}
        \binom{m}{s_0, s_1, \cdots, s_{q-1}} \le \exp(m\mathrm{H}(\bs/m)),
    \end{equation*}
    where 
    \begin{equation*}
        \mathrm{H}(\bs/m) \=-\sum_{i=0}^{q-1} \frac{s_i}{m} \cdot \log \frac{s_i}{m}
    \end{equation*}
    is the Shannon entropy with a natural logarithm. Thus,
    \begin{equation}\label{eq:Mpestimate}
        \log M(\bp) \le \sum_{\bw^\ast \in \Sigma_k(\bz)} \cou_{k, \bz}(\bw^\ast)\cdot \mathrm{H}\biggl(\frac{(\cou_{k+1, \bz}(\bw^\ast \alpha))_{\alpha\in \bQ}}{\cou_{k, \bz}(\bw^\ast)}\biggr).
    \end{equation}
    Consider the following set of random variables. Let $i$ be a uniformly random integer in $[0, n+k-1]$, define random variables $\bV_r$, $r = 0,\, 1,\, \ldots,\, k+1$, by
    \begin{equation*}
        \bV_r \= \sigma^i  (\bz)_{[1,r]} \in \Sigma_r.
    \end{equation*}
    Note that $\bV$ with larger $r$ determine the shorter $r$ ones:
    \begin{equation*}
        \bV_{r'} = {\bV_r}_{[1,r']}.
    \end{equation*}
    Thus, the Shannon entropy of $\bV_{k+1}$ can be expressed in terms of conditional entropies:
    \begin{equation*}
        \mathrm{H}(\bV_{k+1}) 
        = \mathrm{H}(\bV_0)+\sum_{r = 0}^k \mathrm H(\bV_{r+1}|\bV_r) 
        = \sum_{r = 0}^k \mathrm H(\bV_{r+1}|\bV_r). 
    \end{equation*}
    Next, note that for each integer $r >1$,
    \begin{equation*}
        \begin{aligned}
        \mathrm{H}(\bV_{r+1}|\bV_r) 
        &= \mathrm{H}\bigl(\sigma^{i+1}(\bz)_{[1,r]}\big|\bV_r\bigr)
        = \mathrm{H}\bigl(\sigma^{i+1}(\bz)_{[1,r]}\big|\sigma^{i+1}(\bz)_{[1,r-1]}, \bV_1\bigr) \\ 
        &\le \mathrm{H}\bigl(\sigma^{i+1}(\bz)_{[1,r]}\big|\sigma^{i+1}(\bz)_{[1,r-1]}\bigr) 
        = \mathrm{H}(\bV_{r}|\bV_{r-1}),
        \end{aligned}
    \end{equation*}
    as $i$ and $i+1$ are identically distributed. Thus,
    \begin{equation}        \label{eq:entropysum}
        \mathrm{H}(\bV_{k+1}) \ge (k+1)\cdot \mathrm{H}(\bV_{k+1}|\bV_k).
    \end{equation}
    Note that this conditional entropy is related to the expression in \eqref{eq:detlapl3}:
    \begin{equation*}
        \begin{aligned}
        \mathrm{H}(\bV_{k+1}|\bV_k) 
        &= \sum_{\bw^\ast \in \Sigma_k} p(\bV_k = \bw^\ast) \cdot \mathrm H(\bV_{k+1}|\bV_k=  \bw^\ast) \\ 
        &= \frac{1}{n+k}\sum_{\bw^\ast \in \Sigma_k(\bz)} \cou_{k,\bz}[\bw^\ast] \cdot \mathrm H(\bV_{k+1}|\bV_k = \bw^\ast).
        \end{aligned}
    \end{equation*}
    In particular, from \eqref{eq:Mpestimate} and \eqref{eq:entropysum},
    \begin{equation*}
        \frac{k+1}{n+k}\cdot \log M(\bp) \le \mathrm{H}(\bV_{k+1}).
    \end{equation*}
    Since random integer $i$ determines $\bV_{k+1}$, we then get
    $        \frac{k+1}{n+k}\cdot \log M(\bp) \le \mathrm{H}(i) = \log (n+k)$.
    Taking an exponent finishes the proof.
\end{proof}

\section{Small entropy reduction}
\label{sec:thetareduction}
For a system $(X, T)$, define
\begin{align*}
\sP\sU(X, T) & \= \{ f \in  C(X) : f \text{ has periodic optimization property} \}.
\end{align*}

To prove Theorem~\ref{thm:main} for all subshifts of finite type, we first reduce the question to the case where the entropy of $(\Sigma_{\bF}, \sigma)$ is small, which is a key condition in Section \ref{sec:B}. More precisely, the aim of this section is to reduce Theorem~\ref{thm:main} to the following theorem:

\begin{theorem}\label{thm:maim.redu}
	There exists a constant $M>0$ such that the following holds for all $q \ge 2$ and $\theta \in (0, 1)$:
    
    Let $(X,\sigma)$ be any subshift of finite type on $q$ symbols, equipped with the metric $d_{\ba}$ given by $\ba=\{a_n\}_{n\in\N} = \{\theta^n\}_{n\in\N}$, that satisfies 
	\begin{equation}    \label{eq:main.redu}
		M h_{\operatorname{top}}(\sigma) + \log \theta <0.
	\end{equation}
	Then $\sP\sU (X,\sigma) \cap \Lip(X,d_{\ba})$ is a prevalent subset of $\Lip(X,d_{\ba})$.
\end{theorem}
We will use Theorem~\ref{thm:maim.redu} to prove Theorem~\ref{thm:main} at the end of this section.

\subsection{Prevalence of zero entropy}

The main tool to reduce the main theorem to Theorem~\ref{thm:maim.redu} is the prevalence of zero-entropy condition. Let $(\Sigma_{\bF}, \sigma)$ be a subshift of finite type of Theorem~\ref{thm:main} on $q$ symbols. We define
\begin{equation*}
	\sE_0 \= \{ f \in \Lip(\Sigma_{\bF}, d_{\ba}) : h_{\mu}(\sigma) =0 \text{ for every } \mu\in \Mmax(\Sigma_{\bF}, \sigma,f) \}.
\end{equation*}

First, we note that:
\begin{lemma}\label{lem:zero.top.entropy}
	If $f \in \sE_0$ and $K$ is the closure of the union of the supports of all maximizing measures of $f$, then $(K, 
    \sigma)$ is a subshift and $h_{\operatorname{top}}(\sigma|_{K})=0$.
\end{lemma}

\begin{proof}
    Clearly, $\sigma(\supp \mu) = \supp \mu$ for every $\mu\in \Mmax(\Sigma_{\bF},\sigma,f)$. Thus, continuity of $\sigma$ gives that $\sigma(K) \subseteq K$. Since $\sigma(K)$ is dense in $K$ and $K$ is compact, we have $\Sigma(K) =K$.

	By the Ma\~n\'e lemma, there exists $u \in \Lip(\Sigma_{\bF}, d_{\ba})$ such that $g \= f + u - u \circ \sigma - \mpe(\Sigma_{\bF},f) \le 0 = \mpe(g)$. Note that $\Mmax(\sigma,f) = \Mmax(\sigma,g)$. Clearly, $g$ is a Lipschitz function. For every $\mu \in \Mmax(\sigma,f)$, we have $g=0$ on $\supp \mu$. So $g = 0$ on $K$. Thus, for every $\mu' \in \cM(\Sigma_{\bF}, \sigma)$, if $\supp \mu' \subseteq K$, then $\mu' \in \Mmax(\sigma, g) = \Mmax(\sigma,f)$, so $h_{\mu'}(\sigma)=0$. By the variational principle, $h_{\operatorname{top}}(\sigma|_{K}) =0$.
\end{proof}

The remaining part of this subsection is dedicated to prove:

\begin{prop}\label{prop:pre.zero.entropy}
	The subset $\sE_0 $ is  prevalent in $\Lip(\Sigma_{\bF}, d_{\ba})$. 
\end{prop}

First, we construct a Hilbert brick to estimate probability. Fix $\theta_* \in (0,\theta)$ and define
\begin{equation*}
	 b_{\bw} \= \theta_*^{\abs{\bw}} 
\end{equation*}
for every $\bw \in \Sigma_*(\bF)$. Clearly, we have $b_{\bw} = \cO (a_{\abs{\bw}})$ as $\abs{\bw} \to +\infty$. By Lemma~\ref{lem:hilbert-brick}, this Hilbert brick $\cH[\bb]$ is compact. Recall $b_n = \theta_*^n$ for $n\in\N$ (cf.\ \eqref{eq:def.bn}). 

Following Morris \cite[p.~6]{Mo08}, we have the following result.

\begin{lemma}  \label{entropty_est}
	Let $\boo$ be a periodic orbit with $\abs{\boo} = n$. Let $\gamma \in (0,1]$, $m\in \N$, and $\mu \in \cM(\Sigma_{\bF}, \sigma)$. If $\mu\bigl( \bigcup_{\bw \in \Sigma_m^c(\boo)} [\bw] \bigr) \le \gamma$, then there exists $0 < \gamma' \le \gamma$ such that
	\begin{equation*}
		h_{\mu}(\sigma)  \le \frac{1-\gamma'}{m} \log n + \gamma' \log q + \frac{2}{m}.
	\end{equation*}
\end{lemma}

\begin{proof}
	It follows from \cite[Lemma~2.1]{Mo08} that if $x_1, \, x_2,\, \ldots,\, x_k$ are nonnegative, then 
	\begin{equation}\label{e:Mo.lem2.1}
		-\sum_{i=1}^{k} x_i \log x_i \le 1+  \sum_{i=1}^{k} x_i  \log k.
	\end{equation}
	Now denoting $\gamma' \= \mu\bigl( \bigcup_{\bw \in \Sigma_m^c(\boo)} [\bw] \bigr)$, we have $0< \gamma' \le \gamma$. By definition, we have 
	\begin{equation}\label{e:gamma'}
		\gamma' = \sum_{\bw \in \Sigma_m^c(\boo)} \mu([\bw]).
	\end{equation}
	Combining (\ref{e:Mo.lem2.1}), (\ref{e:gamma'}), and the fact that $\abs{\Sigma_m(\boo)} \le n$, we obtain that
	\begin{align*}
		h_{\mu}(\sigma) 
        &\le -\frac{1}{m} \sum_{\bw \in \Sigma_m(\boo)} \mu([\bw]) \log \mu([\bw]) - \frac{1}{m} \sum_{\bw \in \Sigma_m^c(\boo)} \mu([\bw]) \log \mu([\bw]) \\
		& \le \frac{1}{m} (1+ (1-\gamma') \log n) + \frac{1}{m} (1+ \gamma' \log(q^m)) 
		 \le \frac{1-\gamma'}{m} \log n + \gamma' \log q + \frac{2}{m}.
	\end{align*}
	This concludes this lemma.
\end{proof}

For each $f\in \Lip(\Sigma_{\bF}, d_{\ba})$ and each $n\in \N$, we define 
\begin{equation*}
	\mpe_n (f) \= \sup_{\boo \in \PO_{\le n}} \mpe(\boo, f).
\end{equation*}
For each $\delta>0$, define 
\begin{equation}\label{eq:def.Hdelta}
	\cH[\bb,\delta] \= \Bigl\{  \sum_{\bw \in \Sigma_*(\bF)} g_{\bw} \mathbbold{1}_{[\bw]} : g_{\bw}\in\R, \, \abs{g_{\bw}} \le b_{\bw} \bigl(1- (2q)^{-\abs{\bw}}\delta \bigr) \text{ for each } \bw \in \Sigma_*(\bF) \Bigr\},
\end{equation}
$\cI[\bb] \= \prod_{\bw\in \Sigma_*(\bF)} [-b_{\bw}, b_{\bw}]$, $\marg_n \= b_n(2q)^{-n}\delta$ for every $n\in \N$, and
\begin{equation*}
    \cI[\bb,\delta] \= \prod_{\bw\in \Sigma_*(\bF)} \bigl[-b_{\bw} \bigl(1- (2q)^{-\abs{\bw}}\delta \bigr), b_{\bw} \bigl(1- (2q)^{-\abs{\bw}}\delta \bigr) \bigr].
\end{equation*}
 Equip $ \cI[\bb]$ with the product probability $\nu_{\times}$ of Lebesgue probability measures $\nu_{\bw}$ on each $[-b_{\bw}, b_{\bw}]$, and equip $\cH[\bb]$ with the probability $\nu$ as the push-forward of $\nu_{\times}$ by the projection $\Pi \: \cI[\bb] \to \cH[\bb]$ defined in \eqref{eq:def_Pi}.
Straightforward calculation yields
\begin{equation}\label{eq:prob.subbrick}
	\nu_{\times} (\cI[\bb,\delta]) = \prod_{i=1}^{+\infty} \bigl(1-(2q)^{-i}\delta \bigr)^{q^i} >  \prod_{i=1}^{+\infty} \bigl(1-2^{-i}\delta \bigr) > 1-\delta.
\end{equation}
Moreover, the measure $\nu_{\times}$ is translation-invariant. More precisely, if $E \subseteq \cI[\bb]$ is a Borel set, $\bg \in \cI[\bb]$, and $E + \bg \subseteq \cI[\bb]$, then $\nu_{\times}(X) = \nu_{\times} (E+\bg)$.

Define
\begin{equation*}
	\sN \= \{ f\in \Lip(\Sigma_{\bF}, d_{\ba}) : \mpe(\boo, f) \neq \mpe(\boo', f) \text{ if } \boo, \, \boo' \in \PO_* \text{ and } \boo \neq \boo' \},
\end{equation*}

We first prove that $\sN$ is ``almost everywhere''.

\begin{lemma}\label{lem:sN.fullprob}
	For every $f\in \Lip(\Sigma_{\bF}, d_{\ba})$, we have $\nu_{\times}(\{ \bg \in \cI[\bb] : f + \Pi(\bg) \in \sN \}) =1$.
\end{lemma}

\begin{proof}
	For $f\in \Lip(\Sigma_{\bF}, d_{\ba})$, define $\sN(f) \= \{\bg \in \cI[\bb] : f + \Pi(\bg) \in \sN \}$. Let $\{ \boo_i \}_{i\in \N}$ be all periodic orbits. For distinct positive integers $i, \, j$, define $\sN_{i,j}(f) \= \{ \bg \in \cI[\bb] : \mpe(\boo_i, f + \Pi(\bg)) \neq \mpe(\boo_j, f+ \Pi(\bg)) \}$. We have $\sN(f) = \bigcap_{1\le i < j} \sN_{i,j}(f)$. So it is sufficient to prove $\nu_{\times}( \sN_{i,j}(f) ) = 1$ for all integers $1\le i < j$.
	
	Now fix integers $1\le i<j$. There exists $\bw' \in \Sigma_*(\bF)$ such that $\conc_{\abs{\bw'}, \boo_i}(\bw') \ne \conc_{\abs{\bw'}, \boo_j}(\bw')$. Given any $\{ g_{\bw}\}_{\bw \in \Sigma_*(\bF) \smallsetminus \{ \bw' \}}$, we have
	\begin{equation*}
		\nu_{\bw'} ( \{ g_{\bw'} : \{ g_{\bw} \}_{\bw \in \Sigma_*(\bF) }  \in \sN_{i,j}(f)  \}  ) =1 .
	\end{equation*}
	Thus, Fubini's theorem gives $\nu_{\times} ( \sN_{i,j}(f) ) =1$. This concludes this lemma.
\end{proof}

\begin{proof}[Proof of Proposition~\ref{prop:pre.zero.entropy}]
	For every $\ve>0$, define
	\begin{equation*}
		\sE_{\ve} \= \{ f\in \Lip(\Sigma_{\bF}, d_{\ba}) : \sup \{ h_{\mu}(\sigma) : \mu\in \Mmax(\sigma,f) \} \le \ve \}.
	\end{equation*}
	Note that it is sufficient to prove that $\sE_{\ve}$ is prevalent for each $\ve>0$. Moreover, it is sufficient to prove that for each $f \in \Lip(\Sigma_{\bF}, d_{\ba})$ and each $\ve>0$,
    \begin{equation}\label{eq:ze.suffcon}
        \nu (  \{  g\in \cH[\bb] : f + g \in \sE_\ve \} ) \ge 1- \ve.
    \end{equation}
	Since $\nu = \Pi_{*}(\nu_{\times})$, (\ref{eq:ze.suffcon}) is equivalent to 
    \begin{equation}\label{eq:ze.suffcon'}
        \nu_{\times} (  \{  \bg \in \cI[\bb] : f + \Pi(\bg) \in \sE_\ve \} ) \ge 1- \ve.
    \end{equation}

	Now fix $f \in \Lip(\Sigma_{\bF}, d_{\ba})$ and $\ve>0$. 
	Define $\delta \= \ve/2$. By (\ref{eq:prob.subbrick}), we have $\nu_{\times}( \cI[\bb,\delta] ) > 1- \ve/2$. For every $f'\in f + \cH[\bb]$, by the Bressuad--Quas closing lemma (see \cite[Theorem~2]{BQ07} and \cite[Lemma~2.2]{HLMXZ25}), we have
	\begin{equation*}
		\lim_{n\to +\infty} (\log(\mpe(f') - \mpe_{n} (f') )) / \log n = -\infty.
	\end{equation*}
	Note that $b_n = \theta_*^n$. Thus, defining $m_n \= \bigl\lceil \frac{3 \log n}{\ve} \bigr\rceil$, and recalling $\marg_{m_n} \= b_{m_n}(2q)^{-m_n}\delta$, for each $f' \in f+ \cH[\bb]$, we have
	\begin{equation}\label{eq:point.conv}
		\lim_{n\to +\infty}  (\mpe(f') - \mpe_n(f')) / \marg_{m_n} = 0.
	\end{equation}
	Set $\gamma \= \frac{\ve}{3 \log q}$ and denote
	\begin{equation}\label{eq:Sn}
		G_n \= \{ \bg\in  \cI[\bb,\delta] :  (\mpe(f+\Pi(\bg)) - \mpe_n(f+ \Pi(\bg))) / \marg_{m_n} < \gamma \}
	\end{equation}
	for each $n\in \N$.
	Since the convergence in (\ref{eq:point.conv}) holds pointwise in $f+\cH[\bb]$ and $\nu_{\times}(\cI[\bb]) $ is finite, it converges almost surely. Thus, $\lim_{n \to +\infty} \nu_{\times}(G_n) = \nu_{\times}(\cI[\bb, \delta]) > 1 - \frac{\ve}{2}$, then there exists $N\in \N$ such that $\nu_{\times}(G_N) > 1-\ve$ and we set $G \= G_N$ and fix $m \= m_N$. Assume without loss of generality that $N> n_{\theta}$. Clearly, $G$ is a Borel set. For each $\boo \in \PO_{\le N}$, denote
	\begin{equation*}
		G_{\boo} \= \{ \bg\in \cI[\bb] : \mpe(\boo, f+ \Pi(\bg)) >\max\{ \mpe(\boo', f+\Pi(\bg)) : \boo' \in \PO_{\le N}, \, \boo'\neq \boo\} \} \cap G.
	\end{equation*}
	Now, we prepare a $\delta$-argument. Clearly, the sets $G_{\boo}$, $\boo \in \PO_{\le N}$, are mutually disjoint, and by Lemma~\ref{lem:sN.fullprob}, we have
	\begin{equation}\label{eq:probeq}
		\sum_{\boo \in \PO_{\le N}} \nu_{\times} (G_{\boo}) = \nu_{\times}(G) > 1-\ve.
	\end{equation}  
	Now fix $\boo \in \PO_{\le N}$ and for $\bg\in G_{\boo} \subseteq \cI[\bb,\delta]$ define
    \begin{equation*}
    \cF_{\boo, \bw}(g_{\bw}) \= \begin{cases}
        g_{\bw}+\marg_{m}&\text{ if } \bw\in \Sigma_{m}(\boo),\\
        g_{\bw}& \text{ otherwise}.
        \end{cases}
    \end{equation*} 
	Let $\cF_{\boo}$ act coordinate-wise on $G_{\boo}$. Since $\bg\in \cI[\bb,\delta]$, by (\ref{eq:def.Hdelta}), $\cF_{\boo} (\bg) \in \cI[\bb]$.
	Note that 
    \begin{equation*}
        \mpe (\boo, f+\Pi ( \cF_{\boo}(\bg))) > \mpe (\boo', f+\Pi ( \cF_{\boo}(\bg)))
    \end{equation*}
    for all other $\boo' \in \PO_{\le N}$. Thus, $\{\cF_{\boo}(G_{\boo}) \}_{\boo \in \PO_{\le N}}$ are mutually disjoint. So, using (\ref{eq:probeq}) and that $\nu_{\times}$ is translation-invariant, we obtain
	\begin{equation}\label{eq:prob.lowentropy}
		\nu_{\times} \Bigl( \bigcup_{\boo \in \PO_{\le N}} \cF_{\boo}(G_{\boo}) \Bigr) 
        = \sum_{\boo \in \PO_{\le N}} \nu_{\times}  (\cF_{\boo}(G_{\boo})   ) 
        = \nu_{\times}(G)
        > 1-\ve.
	\end{equation}
	Now for every $\bg\in \cF_{\boo}(G_{\boo})$, let $\mu \in \Mmax (\sigma,f+ \Pi( \bg) ) $. Assuming $\bg = \cF_{\boo}(\bg')$, we have
    \begin{align*}
        \marg_m\mu\Bigl( \bigcup_{\bw \in \Sigma_m(\boo)} [\bw] \Bigr) 
        &= \mpe(\mu, f+\Pi(\bg)) - \mpe(\mu, f+\Pi(\bg')) \\
        &\ge \mpe(f+ \Pi(\bg)) - \mpe_N(f+\Pi(\bg')) -\gamma \marg_m \\
        &=  \mpe(f+\Pi(\bg))- \mpe_N(f+\Pi(\bg)) + \marg_m - \gamma \marg_m \\
        &\ge (1-\gamma)\marg_m,
    \end{align*}
    where we used \eqref{eq:Sn} since $\bg' \in G$. Thus, we have $\mu\bigl( \bigcup_{\bw \in \Sigma_m^c(\boo)} [\bw] \bigr) \le \gamma$. Recall that $\gamma = \frac{\ve}{3\log q}$, $m = \bigl\lceil \frac{3 \log N}{\ve} \bigr\rceil \ge  \frac{3 \log N}{\ve} $, and $N > n_{\theta}$. So, applying Lemma~\ref{entropty_est}, we obtain 
	\begin{equation*}
		h_{\mu}(\sigma) < \frac{\log N }{m} + \gamma \log q + \frac{2}{m} \le \frac{(\log N +2) \ve}{3 \log N} + \frac{\ve}{3} \le \ve.
	\end{equation*}
	Therefore, $f+ \Pi \bigl( \bigcup_{\boo \in \PO_{\le N}} \cF_{\boo}(G_{\boo}) \bigr) \subseteq \sE_\ve$, so by (\ref{eq:prob.lowentropy}) we obtain (\ref{eq:ze.suffcon'}).
\end{proof}

\subsection{Reduction to small entropy}

In this subsection, we mainly use Proposition~\ref{prop:pre.zero.entropy} and Theorem~\ref{thm:maim.redu} to prove Theorem~\ref{thm:main}.

Let $\fC$ be the space of compact subsets of $\Sigma_{\bF}$ equipped with the Hausdorff distance $d_{\operatorname{H}}$ and define
\begin{equation*}
	\fC_{\operatorname{inv}} \= \{ X \in \fC : \sigma(X)=X \} \quad \text{ and } \quad
	\fC_0 \= \{ X \in \fC_{\operatorname{inv}} :  h_{\operatorname{top}}(\sigma|_{X})=0 \}.
\end{equation*}
Since $\Sigma_{\bF}$ is compact, $\fC$ is compact. So $\fC_0$ is second countable.

\begin{lemma}\label{lem:lowentropy.decomp}
	For each $\ve>0$ and each $X \in \fC_0$, there exists a subshift of finite type $X_{\ve} \in \fC_{\operatorname{inv}}$ containing $X$ with $h_{\operatorname{top}} (\sigma|_{X_{\ve}}) < \ve$ such that every $X' \in \fC_{\operatorname{inv}}$ sufficiently close to $X$ is contained in $X_{\ve}$.
\end{lemma}

\begin{proof}
	For each $n\in \N$, define a subshift of finite type 
	\begin{equation}\label{eq:def.SFTYn}
		Y_n \= \{ \alpha_1 \alpha_2 \ldots \in \Sigma_{\bF} : \alpha_i \alpha_{i+1} \ldots \alpha_{i+n-1} \in \Sigma_n(X) \text{ for every } i\in \N \}.
	\end{equation}
	By definition, we have $X \subseteq Y_n$ and $\lim_{i\to +\infty} Y_i =X$. Since $h_{\operatorname{top}} \: \fC_{\operatorname{inv}} \to [0,+\infty)$ is upper-semicontinuous (see e.g.\ \cite[Proposition~1.3]{FT17}), we have $\lim_{i \to +\infty} h_{\operatorname{top}} (\sigma|_{Y_i}) =0$. Now we let $X_{\ve}$ be some $Y_m$ with topological entropy smaller than $\ve$.
    
    Thus, for every $X' \in \fC_{\operatorname{inv}}$, if $d_{\operatorname{H}} (X,X') < a_m$, then $\Sigma_m(X') = \Sigma_m(X) $, by (\ref{eq:def.SFTYn}), we have $X' \subseteq X_{\ve}$. This lemma is therefore proved.
\end{proof}

\begin{cor}\label{cor:sep.lowentropy}
	For each $\ve>0$, there exists a sequence of subshifts of finite type $\{Y_n\}_{n\in \N}$ with topological entropy smaller than $\ve$ such that every $X \in \fC_0$ is contained in one of them.
\end{cor}

\begin{proof}
	For every $X \in \fC_0$, we obtain a subshift of finite type $X_{\ve}$ from Lemma~\ref{lem:lowentropy.decomp} and a neighborhood $U(X) \subseteq \fC_0$ of $X$ such that every $X' \in U(X)$ is contained in $X_{\ve}$. Thus, $\{ U(X)  \}_{X\in \fC_0}$ gives an open covering of $\fC_0$. Since $\fC_0$ is second countable, there exists a countable subcovering $\{ U(X_n) \}_{n\in \N}$. Setting $Y_n \= (X_n)_{\ve}$, $n\in\N$, completes this corollary.
\end{proof}

\begin{lemma}\label{lem:pre.restriction}
	Let $(X,d)$ be a compact metric space and $Y\subseteq X$ be a compact subspace. Let $A$ be a shy subset of $\Lip(Y, d)$. Then the set $B \= \{ f \in \Lip(X, d) : f|_{Y} \in A \}$ is a shy subset in $\Lip(X,d)$. 
\end{lemma}

\begin{proof}
	Let $R \: \Lip(X, d) \to \Lip(Y, d)$ be defined as $R(f) \= f|_Y$. Clearly, $R$ is a continuous linear map and $B=R^{-1}(A)$. By the McShane extension theorem, $R$ is surjective. Since $\Lip(X, d)$ and $\Lip(Y, d)$ are Banach spaces, by the Bartle--Graves theorem, there exists a continuous linear map $E \: \Lip(Y,d) \to \Lip(X,d)$ such that $R \circ E = \id_{\Lip(Y, d)}$.
	
	Since $A$ is shy in $\Lip(Y,d)$, there exists a compactly supported measure $\nu$ such that $\nu (f+ A)=0$ for each $f\in \Lip(Y,d)$. Set $\nu_E \= E_*(\nu)$, where $E_*$ is the pushforward map of $E$. As $E$ is continuous, $\nu_E$ is also compactly supported. For each $g\in \Lip(X,d)$, we have
	\begin{equation*}
		\nu_E(g + B) = \nu \bigl(E^{-1}(g+B)\bigr) = \nu(R(g) + A) =0.
	\end{equation*}
	Therefore, $B$ is shy in $\Lip(X,d)$.
\end{proof}

\begin{lemma}\label{lem:prelocking}
    Under the hypotheses of Theorem~\ref{thm:main},  $(\sP \sU(\Sigma_{\bF}, \sigma) \cap \Lip(\Sigma_{\bF}, d_{\ba}) ) \smallsetminus \Lock_{d_{\ba}}(\Sigma_{\bF}, \sigma)$ is shy in $\Lip(\Sigma_{\bF}, d_{\ba})$.
\end{lemma}

\begin{proof}
    Since there are only countably many periodic orbits, it is sufficient to prove $\{ f\in \Lip(\Sigma_{\bF}, d_{\ba}) : \Mmax(f) = \{\mu_{\boo}\} \} \smallsetminus \Lock_{d_{\ba}}(\Sigma_{\bF}, \sigma)$ is shy in $\Lip(\Sigma_{\bF}, d_{\ba})$ for every $\boo \in \cO_*$. Now fix $\boo \in \cO_*$. It can be seen from the proofs in the note \cite{BZ15} that if $\Mmax(\sigma,f) = \{\mu_{\boo}\}$, then $f -td(\,\cdot\,, \boo) \in \Lock_{d_{\ba}}(\Sigma_{\bF},\sigma)$ for each $t>0$. We consider the one-dimensional compact space $\cH_{\boo} \= \{ td(\,\cdot\,, \boo) : t\in [0,1] \}$ equipped with the Lebesgue measure. For each $f\in \Lip(\Sigma_{\bF}, d_{\ba})$, there is at most one function $g$ in $f + \cH_{\boo}$ such that $\Mmax(g) = \{\mu_{\boo}\}$ but $g\notin \Lock_{d_{\ba}}(\Sigma_{\bF}, \sigma)$, and this point is of zero measure. Thus, this lemma is true.
\end{proof}

Finally, we prove Theorem~\ref{thm:main}.

\begin{proof}[\textbf{Proof of Theorem~\ref{thm:main} assuming Theorem~\ref{thm:maim.redu}}]
	Set 
    \begin{equation}\label{eq:htop.assume}
        \ve \= -M^{-1} \log \theta.
    \end{equation}
    By Corollary~\ref{cor:sep.lowentropy}, we obtain a sequence of subshifts of finite type $\{Y_n\}_{n\in \N}$ with  topological entropy smaller than $\ve$ such that every $X \in \fC_0$ is contained in one of them. 
	
	For each $f\in \Lip(\Sigma_{\bF}, d_{\ba})$, let $K(f)$ be the closure of the union of the supports of all maximizing measures of $f$. Define
	\begin{align}
		L_* &\= \{ f\in \Lip(\Sigma_{\bF}, d_{\ba}) : h_{\operatorname{top}}(\sigma|_{K(f)}) >0 \} \quad\text{ and} \\
		L_n &\= \{ f\in \Lip(\Sigma_{\bF}, d_{\ba}) : K(f) \subseteq Y_n\}
	\end{align}
	for each $n\in \N$. Thus, the choice of $\{Y_n\}_{n\in \N}$ gives us a division 
	\begin{equation}\label{eq:Lipdiv}
		\Lip(\Sigma_{\bF}, d_{\ba}) = L_* \cup   \bigcup_{n=1}^{+\infty} L_n .
	\end{equation}
	By Proposition~\ref{prop:pre.zero.entropy}, $L_*$ is shy. 
    Our choice (\ref{eq:htop.assume}) implies that every $(Y_n, \sigma)$ satisfies the hypotheses (\ref{eq:main.redu}) in Theorem~\ref{thm:maim.redu}.
    Thus, applying Theorem~\ref{thm:maim.redu} and Lemma~\ref{lem:pre.restriction}, the set 
	\begin{equation*}
		\wt{L}_n \= \{ f\in \Lip(\Sigma_{\bF}, d_{\ba}) : f|_{Y_n} \notin  \sP\sU (Y_n,\sigma) \}
	\end{equation*}
	is shy in $\Lip(X, d_{\ba})$. Note that $L_n \smallsetminus \sP\sU(X, \sigma) \subseteq \wt{L}_n$, so $L_n \smallsetminus \sP\sU(X, \sigma)$ is shy in $\Lip(\Sigma_{\bF}, d_{\ba})$. Thus, by (\ref{eq:Lipdiv}), we have
	\begin{equation*}
		\Lip(\Sigma_{\bF}, d_{\ba}) \smallsetminus \sP\sU(X, \sigma)= (L_* \smallsetminus \sP\sU(X, \sigma)) \cup   \bigcup_{n=1}^{+\infty} (L_n \smallsetminus \sP\sU(X, \sigma))   ,
	\end{equation*}
	where every term on the right-hand side above is shy in $\Lip(X, d_{\ba})$, so $\Lip(\Sigma_{\bF}, d_{\ba}) \smallsetminus \sP\sU(X, \sigma)$ is shy in $\Lip(X, d_{\ba})$. Finally, applying Lemma~\ref{lem:prelocking} finishes the proof of Theorem~\ref{thm:main}.
\end{proof}

\section{Cycle and path advantages}
\label{sec:advantage}

Now, we state a sufficient condition for locking which will concern only the periodic orbits. This way, we reduce our problem to studying periodic orbits, which are easier to analyze compared to general invariant measures.

We generally follow the approach of \cite{BZ16} and \cite{DLZ24}. The rough idea is that if some orbit $\boo \in \PO_*$ has much larger $\mpe(\boo, f)$ compared to all other $\mpe(\boo', f)$, then $\mu_{\boo}$ would be the unique maximizing measure. Not every $\boo'$ has to be considered: Lemma~4.1 of \cite{BZ16} shows that it is enough to estimate $\mpe(\boo, f) - \mpe(\boo', f
)$ from below for $\boo'$ with a large enough gap $\Delta(\boo')$ (i.e., orbits with corresponding in the same level de Bruijn--Good graph as $\boo$). The gap appears here due to Lemma~\ref{lemma.gapcompose}: smaller gap orbits' ergodic suprema can be approximated by ergodic suprema of larger gap orbits, so they must be small as well.

We modify this idea slightly. In theory, $\mpe(\boo, f) - \mpe(\boo', f)$ may be small just because $\boo'$ is very similar to $\boo$. To avoid this issue, we introduce an $m$-advantage of an orbit $\boo$ with respect to an orbit $\boo'$. Precisely, for $m\in\N$, we define
\begin{equation}\label{eq:betadef}
    \mpe_m(\boo, \boo', f) \= \mpe(\boo, A_mf) - \mpe(\boo', A_mf),
\end{equation}
and if $\mu_{\boo'}([\Sigma_m^c(\boo)]) \neq 0$ or $\abs{\boo'} < \abs{\boo}$, we define:
\begin{equation}    \label{eq:betastardef}
    \mpe^{\sharp}_m(\boo, \boo', f) \= \begin{cases} \frac{\mpe_m(\boo, \boo', f)}{\mu_{\boo'}([\Sigma_m^c(\boo))]}&\text{ if }\abs{\boo'} \ge \abs{\boo},\\
    \mpe_m(\boo, \boo', f)&\text{ if } \abs{\boo'} < \abs{\boo},
    \end{cases}
\end{equation}
where $A_mf$ was introduced in Subsection~\ref{subsec:AF}. Note that if $\Delta(\boo) > a_m$, Lemma~\ref{lem:sameconc} implies that $\mu_{\boo'} ([\Sigma_{m}^c(\boo)] )\neq 0$ holds, so $\mpe_m^\sharp (\boo,\boo', f)$ is well defined. 
\begin{lemma}\label{lem:betatransitivity}
    Suppose $m\in\N$ satisfies $m +1 \ge \sup_{\bw \in \bF}\abs{\bw}$. Let $\boo,\, \boo',\, \boo'' \in \PO_*$ have $\Delta(\boo),\,\Delta(\boo'),\, \Delta(\boo'') > a_m$, assume $\boo \ne \boo''$, and suppose $\mpe_m^\sharp(\boo, \boo', f), \, \mpe_m^\sharp(\boo', \boo'', f) \ge 0$. Then 
    \begin{equation}    
        \mpe_m^\sharp(\boo, \boo'', f) \ge \min\bigl\{\mpe_m^\sharp(\boo, \boo', f),\,\mpe_m^\sharp(\boo', \boo'', f)\bigr\}.
    \end{equation}
    \vspace{-\baselineskip}
\end{lemma}
\begin{proof}
    Since $\Delta(\boo), \, \Delta(\boo') > a_m$, all three $\mpe_m^\sharp$-s in this lemma are well defined. 
    We only consider the case $\abs{\boo} \le \abs{\boo'} \le \abs{\boo''}$. Indeed, we have $0 \le \mpe_m(\boo, \boo' ,f), \, \mpe_m(\boo', \boo'',f) \le \mpe_m(\boo, \boo'',f)$. If $\abs{\boo} > \abs{\boo'}$, then $\mpe_m^\sharp (\boo, \boo'' ,f) \ge \mpe_m (\boo, \boo'' ,f) \ge \mpe_m (\boo, \boo' ,f) \ge \mpe_m^\sharp (\boo, \boo' ,f)$. If $\abs{\boo} > \abs{\boo'}$, similarly, we have $\mpe_m^\sharp (\boo, \boo'' ,f) \ge \mpe_m^\sharp (\boo', \boo'' ,f)$. Both cases make this lemma trivial.
    
    Next, note that for all sets $A$, $B$, and $C$ we have 
    $\abs{C \smallsetminus A} \le \abs{C \smallsetminus B} + \abs{B \smallsetminus A}$.
    Thus, if $0<\abs{B}\le \abs{C}$, then
    $\frac{\abs{C \smallsetminus A}}{\abs{C}} \le \frac{\abs{C \smallsetminus B}}{\abs{C}} + \frac{\abs{B\smallsetminus A}}{\abs{B}}$.
    Applying this for $A \= \Sigma_m(\boo)$, $B \= \Sigma_m(\boo')$, and $C \= \Sigma_m(\boo'')$ gives
    \begin{equation}\label{eq:muoineq}
        \mu_{\boo''}([\Sigma_m^c(\boo)]) \le \mu_{\boo''}([\Sigma_m^c(\boo')])+\mu_{\boo'}([\Sigma_m^c(\boo)]).
    \end{equation}
    Finally, applying (\ref{eq:muoineq}) the inequality 
    $(a+c)/(b+d) \ge \min \{a/b, \, c/d \}$
    for nonnegative $a \= \mpe_m(\boo, \boo',f),\,b \= \mu_{\boo'}([\Sigma_m^c(\boo)]),\,c \= \mpe_m(\boo', \boo'', f),$ and $d \= \mu_{\boo''} ([\Sigma_m^c(\boo')])$ finishes the proof of the lemma.
\end{proof}

More generally, $\Delta(\boo')$ in Lemma~\ref{lem:betatransitivity} need not be large enough. For the case where $\Delta(\boo')$ is small, we have the following analogous lemma.

\begin{lemma}
    \label{lem:betatransitivity2}
    Suppose $m\in\N$ satisfies $m +1 \ge \sup_{\bw \in \bF}\abs{\bw}$. Let $\boo,\, \boo',\, \boo'' \in \PO_*$ with $\Delta(\boo),\, \Delta(\boo'') > a_m$ and $\boo\ne \boo''$. Let $\{\boo_i'\}_{i = 1}^N$ be a decomposition of $\boo'$ from Lemma~\ref{lemma.gapcompose} applied to $a_m$. Assume that $\mpe_m^{\sharp}(\boo, \boo_{i}', f) \ge 0$ for all $1\le i\le N$ and $\mpe_m^{\sharp}(\boo', \boo'', f) \ge 0$. Then 
    \begin{equation*}
        \mpe_m^{\sharp}(\boo, \boo'', f) 
        \ge \min\Bigl\{\min_{1\le i \le N} \mpe_m^{\sharp}(\boo, \boo_i', f), \, \mpe_m^{\sharp}(\boo', \boo'', f)\Bigr\}.
    \end{equation*}
    \vspace{-\baselineskip}
\end{lemma}
\begin{proof}
   Note that $\Delta(\boo)>a_m$, so $\mpe_m^\sharp(\boo, \boo'', f)$ is well defined.
   By (\ref{eq:shadowing2}), we have
   \begin{equation}\label{eq:shadowing3}
       \mpe_m(\boo, \boo' , f) \= \sum_{i=1}^N \frac{\abs{\boo_i'}}{ \abs{\boo'}} \cdot \mpe_m (\boo, \boo_i' ,f).
   \end{equation}
   So there exists an integer $1\le i \le N$ such that $\mpe_m (\boo, \boo_i' ,f) \le \mpe_m(\boo, \boo' , f)$. Note that if $\abs{\boo} > \abs{\boo'}$ implies $\abs{\boo} > \abs{\boo_i'}$ for every integer $1\le i \le N$. Using an argument analogous to Lemma~\ref{lem:betatransitivity}, 
   we only need to consider the case $\abs{\boo} \le \abs{\boo'} \le \abs{\boo''}$. Then, the idea is similar to Lemma~\ref{lem:betatransitivity}. First, using (\ref{eq:shadowing4}), we have
    \begin{equation*}
        \mu_{\boo''}([\Sigma_m^c(\boo)]) \le \mu_{\boo''}([\Sigma_m^c(\boo')]) + \sum_{i= 1}^N\frac{\abs{\boo_i'}}{\abs{\boo'}}\cdot \mu_{\boo_i'}([\Sigma_m^c(\boo)]). 
    \end{equation*}
    If some $\boo_i'$ are shorter than $\boo$, we can replace their concentrations by $1$, preserving the inequality. At the same time,   for $\mpe(\boo, \boo'', f) \= \mpe(\boo, A_mf) - \mpe(\boo'', A_mf)$, (\ref{eq:shadowing2}) gives
    \begin{equation*}
        \mpe(\boo, \boo'', f) = \mpe(\boo', \boo'', f) + \sum_{i = 1}^N\frac{\abs{\boo_i'}}{\abs{\boo'}}\cdot \mpe( \boo, \boo_i', f).
    \end{equation*}
    At last, we to conclude this lemma we only need to apply the following inequality for nonnegative $\{a_i\}_{i=0}^N$ and $\{b_i\}_{i=0}^N$:
    \begin{equation*}
        \frac{\sum_{i = 0}^Na_i}{\sum_{i = 0}^N b_i} \ge \min \Bigl\{ \frac{a_i}{b_i} : i\in \{0,\,1,\, \ldots,\, N\} \Bigr\}.\qedhere
    \end{equation*}
\end{proof}

Let $\widehat{\boo}$ be a periodic orbit that maximizes $\mpe(\boo, A_mf)$. Note that we can always choose $\widehat{\boo}$ with $\Delta > a_m$ due to Lemma~\ref{lemma.gapcompose}. Then we define the \emph{separation}
\begin{equation*}
    \ouv_m^{\sharp}(f) \= \min\bigl\{ \mpe_m^{\sharp}(\widehat{\boo}, \boo', f) : \boo' \in \PO_*, \, \Delta(\boo') > a_m \bigr\}.
\end{equation*}
This object is similar to $\text{gap}_n$ in \cite{BZ16, DLZ24}, we just change the notation to avoid confusion with $\Delta$. Note that $\ouv_m^{\sharp}$ is independent of the choice of $\widehat{\boo}$ because if there are multiple maximizing orbits then $\ouv_m^{\sharp} = 0$.

\subsection{Sufficient locking criterion}
Next, we give the analogue of Lemma~4.1 and Corollary~4.2 of \cite{BZ16}.
\begin{lemma}[$\sharp$-criterion]
    \label{lemma.betastar}
    Suppose $m > n_{\theta}$ is an integer and $f\in \Lip(\Sigma_{\bF}, d_{\ba})$ satisfies
    \begin{equation}        \label{eq:ouvcrit}
        \ouv_m^{\sharp}(f) > 4\LIP(f) \cdot a_{m} \cdot (1-\theta)^{-2}.
    \end{equation}
    Then $f \in \Lock$ and $\mu_{\widehat{\boo}}$ is the unique maximizing measure.
\end{lemma}
    The proof of Lemma~\ref{lemma.betastar} is similar to the one in \cite{BZ16}, so we include it in Appendix~\ref{sec:lockingproof}. The only main difference that allows us to get a bound using $\ouv^{\sharp}$ is the modified estimate on term $\circledd{1}$ in equation (4.8) of \cite{BZ16}, which controls
    \begin{equation}        \label{q.modifiedestimate}
        \sum_{\boo', \Delta(\boo') > a_m} p_m(\boo') \langle A_m f, \mu_{\widehat{\boo}} - \mu_{\boo'}\rangle.
    \end{equation}
    Instead of a trivial bound, we have
    \begin{equation*}
    \begin{aligned}
        &\ouv_m^{\sharp}(f) \Bigl( \sum_{\boo', \Delta(\boo')>a_m, \abs{\boo'} < \abs{\widehat{\boo}}}p_m(\boo')+\sum_{\boo', \Delta(\boo') > a_m, \abs{\boo'}\ge\abs{\widehat{\boo}}} p_m(\boo') \mu_{\boo'}([\Sigma_m^c(\widehat{\boo})]) \Bigr)\\ 
        &\qquad= \ouv_m^{\sharp}(f) \Bigl( 1 - \sum_{\boo', \Delta(\boo') > a_m, \abs{\boo'} \ge \abs{\widehat{\boo}}} p_m(\boo') \mu_{\boo'}([\Sigma_m(\widehat{\boo})])\Bigr) 
        =\ouv_m^{\sharp}(f)(1 - u_m).
    \end{aligned}
    \end{equation*}
The goal of the proof of Theorem~\ref{thm:main} would be to show that \eqref{eq:ouvcrit} holds. For this, we would need to bound $\mpe_m^{\sharp}(\widehat{\boo}, \boo'', f)$ for all $\boo''$ with $\Delta(\boo'') > a_m$. Note that we do not have to compare $\boo''$ to $\widehat{\boo}$: by Lemma~\ref{lem:betatransitivity} it is enough to compare $\boo''$ to some $\boo' \in \PO_*$, provided we later somehow bound $\mpe_m^{\sharp}(\widehat{\boo}, \boo', f)$.
\begin{definition}
    Let $m > n_{\theta}$ be an integer. If for some $\boo''$ with $\Delta(\boo'') > a_m$ there exists $\boo' \in \PO_*$ -- an orbit with $\mpe_m^{\sharp}(\boo', \boo'', f) > 4\LIP(f) \cdot a_{m} \cdot (1 - \theta)^{-2}$, we call $\boo''$ \emph{outadvantaged} for $m$. Let $\PO_{[m]}^{\sharp}(f)$ be the set of orbits $\boo''$ that are not outadvantaged for $m$ and that satisfy $\Delta (\boo'') > a_m$. 
    \label{def:advantagedorbit}
\end{definition}
Note that we do not require $\boo'$ to have a large $\Delta(\boo')$: if $\boo'$ has a small gap, we will use the analogue of Lemma~\ref{lem:betatransitivity} -- Lemma~\ref{lem:betatransitivity2}. 
\subsection{Path advantages}
One way of constructing such an orbit $\boo'$ for $\boo''$ is to replace a path which is a part of $\boo''$ by another path in the same $\Path_{n, \bw', \bw''}$. For two paths $\bp,\, \bp' \in \Path_{n, \bw', \bw''}$ with $\abs{\bw''} \ge n_{\theta}$ we define the advantage $\mpe_m (\bp,\bp' ,f)$ with ingredients defined (\ref{eq:betapaths}) in an analogous way to (\ref{eq:betadef}). Then, we simply let $\mpe_m^{\sharp}(\bp,\bp',f) \= \mpe_m (\bp,\bp' ,f)$.

\begin{figure}
    \begin{tikzpicture}
    \numbercircle{B1}{(0.0,0.0)}{black}{15}{1}
    \numbercircle{B2}{(0.9,0.0)}{black}{0}{0}
    \numbercircle{B3}{(1.8,0.0)}{black}{15}{1}
    \numbercircle{B4}{(2.7,0.0)}{black}{0}{0}
    \numbercircle{B5}{(3.6,0.0)}{black}{0}{0}
    \numbercircle{B6}{(4.5,0.5)}{black}{15}{1}
    \numbercircle{B7}{(5.4,0.5)}{black}{15}{1}
    \numbercircle{B8}{(6.3,0.5)}{black}{0}{0}
    \numbercircle{B9}{(7.2,0.5)}{black}{15}{1}
    \numbercircle{B6a}{(4.5,-0.5)}{black}{0}{0}
    \numbercircle{B7a}{(5.4,-0.5)}{black}{15}{1}
    \numbercircle{B8a}{(6.3,-0.5)}{black}{15}{1}
    \numbercircle{B9a}{(7.2,-0.5)}{black}{0}{0}
    \numbercircle{B10}{(8.1,0.0)}{black}{0}{0}
    \numbercircle{B11}{(9.0,0.0)}{black}{15}{1}
    \numbercircle{B12}{(9.9,0.0)}{black}{15}{1}
    \numbercircle{B13}{(10.8,0.0)}{black}{15}{1}
    \numbercircle{B14}{(11.7,0.0)}{black}{0}{0}

    \numbercircle{B12b}{(9.9,0.6)}{gray}{0}{0}
    \numbercircle{B13b}{(10.8,0.8)}{gray}{0}{0}
    \numbercircle{B14b}{(11.7,1.0)}{gray}{15}{1}
    
    \draw[->, thick] (B1) -- (B2);
    \draw[->, thick] (B2) -- (B3);
    \draw[->, thick] (B3) -- (B4);
    \draw[->, thick] (B4) -- (B5);
    \draw[->, thick] (B5) -- (B6);
    \draw[->, thick] (B6) -- (B7);
    \draw[->, thick] (B7) -- (B8);
    \draw[->, thick] (B8) -- (B9);
    \draw[->, thick] (B9) -- (B10);
    \draw[->, thick] (B10) -- (B11);
    \draw[->, thick] (B11) -- (B12);
    \draw[->, thick] (B12) -- (B13);
    \draw[->, thick] (B13) -- (B14);
    \draw[->, thick, gray] (B11) -- (B12b);
    \draw[->, thick, gray] (B12b) -- (B13b);
    \draw[->, thick, gray] (B13b) -- (B14b);
    \draw[->, thick, gray] (B14b) -- (12.3, 1.0);
    \draw[->, thick] (B5) -- (B6a);
    \draw[->, thick] (B6a) -- (B7a);
    \draw[->, thick] (B7a) -- (B8a);
    \draw[->, thick] (B8a) -- (B9a);
    \draw[->, thick] (B9a) -- (B10);
    \draw[->, thick] (B14) -- (12.3, 0.0);
    \draw[->, thick] (-0.6, 0) -- (B1);

    \draw[dashed, red, thick] (9.45, -1.0) -- (9.45, 1.0);
    \draw[dashed, red, thick] (2.25, -1.0) -- (2.25, 1.0);

     \node[] at (9.25,0.8) {$\bp$};
     \node[red] at (2.25,-1.2) {start};
     \node[] at (9.25,-0.8) {$\bp'$};
     \node[] at (2.45,0.8) {$\bp$};
     \node[] at (2.45,-0.8) {$\bp'$};
     \node[] at (0,0.8) {$\bz$};
     \node[] at (0,-0.8) {$\boo'$};
     \node[] at (3.15,0.3) {$\bw$};
     \node[] at (8.55,0.3) {$\bw'$};
     \node[gray] at (10.8,1.4) {$\bp_{\ext}(\bw')$};
\end{tikzpicture}
\caption{An illustration to Lemma~\ref{lem:pathadv}. Here, $n = 6$, $k=2$, $\bw = (00)$, and $\bw' = (01)$.}
\label{fig:pathadv}
\end{figure}
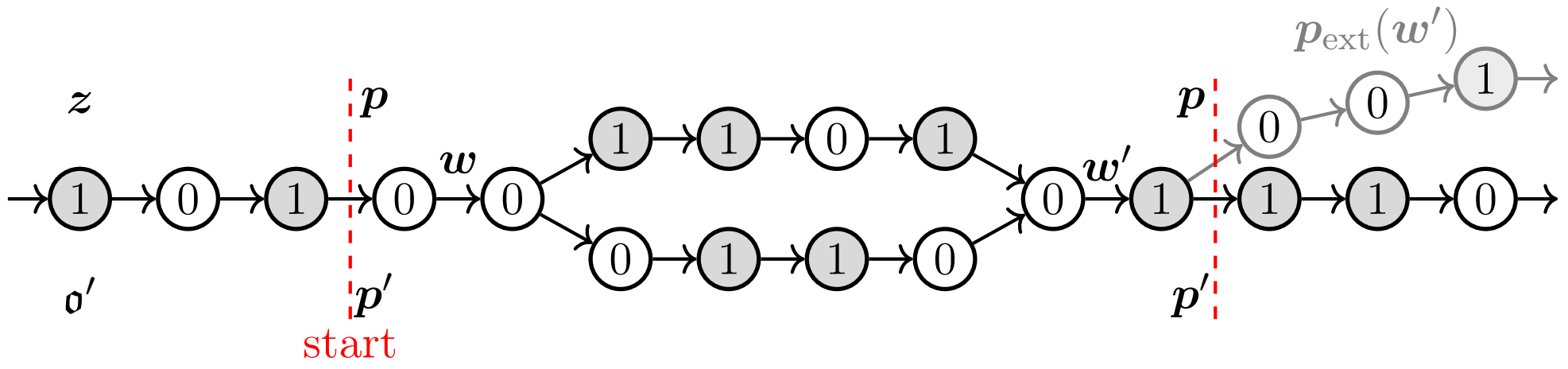
\vspace{-\baselineskip}

\medskip
\medskip

\begin{lemma}
    Suppose $n>k> n_{\theta}$ and $m> k$ are integers. Let $\bw,\, \bw' \in \Sigma_k(\bF)$, and let $\bp,\, \bp' \in \Path_{n, \bw, \bw'}$. Suppose that
    \begin{equation} \label{eq:nkm}
        n+m\le \theta^{-m+k}.
    \end{equation}
    Now assume that $\boo'$ is some periodic orbit with $\Delta(\boo') > a_m$ that contains the word $\bp' \in \Sigma_{n+k}(\bF)$. Let $\abs{\boo'} > m+n+k$. If
    \begin{equation}        \label{eq:betasharppath}
        \mpe^{\sharp}(\bp, \bp', f) > 8(\LIP(f)+1)\cdot a_k \cdot (1-\theta)^{-2} \big/ n,
    \end{equation}
    then $\boo' \notin \PO_{[m]}^{\sharp}(f)$.
    \label{lem:pathadv}
\end{lemma}
\begin{proof}
    Assume that $\bz' \in \boo'$ and that $\bz'_{[1,n+k]} = \bp'$. We consider a word $\bz \in \Sigma_{\bF}$ obtained by replacing this one visit to $\bp'$ modulo the period of $\boo'$ by the visit to $\bp$. In particular, $\bz_{[1,n+k]} = \bp$ and $\sigma^{\abs{\boo'}}(\bz) = \bz$. Moreover, within a period $\abs{\boo'}$, at most $n-k$ symbols of $\bz$ and $\bz'$ between positions $k+1$ and $n$ may differ, see Figure~\ref{fig:pathadv}. Since $k > n_{\theta}$ we did not add any words in $\bF$ when constructing $\bz$, so $\bz\in \Sigma_{\bF}$ is a valid word. Note, however, that $\abs{\boo'}$ does not have to be the minimal period of $\bz$, but let $\boo \in \PO_{\le \abs{\boo'}}$ be the periodic orbit of $\bz$. First, we bound $\mpe(\boo, A_mf) - \mpe(\boo', A_mf)$. Then,
    \begin{equation*}
        \mpe(\boo, A_mf) - \mpe(\boo', A_mf) = \frac{1}{\abs{\boo'}}\sum_{i \in \cI } \bigl(A_mf\bigl(\sigma^i(\bz)\bigr) - A_mf\bigl(\sigma^i(\bz')\bigr)\bigr),
    \end{equation*}
    where $\cI = \{0,\, 1,\, \ldots,\, n-1\}\cup\{\abs{\boo'} - m+k+1,\, \abs{\boo'} - m+k+2,\, \ldots, \, \abs{\boo'}-1\}$, as the middle iterates will be identical. After applying \eqref{eq:fAmfclose}, it is at least
    \begin{equation*}
        \frac{1}{\abs{\boo'}}\sum_{i \in \cI } \bigl(f\bigl(\sigma^i(\bz)\bigr) - f\bigl(\sigma^i(\bz')\bigr)\bigr) - \frac{2(n-k+m-1)\LIP(f)\cdot a_{m+1}}{\abs{\boo'}}.
    \end{equation*}
    We aim to reduce the first sum to the advantage of $\bp$ to $\bp'$. For an integer $i\in [\abs{\boo'}- m+k+1,  \abs{\boo'}-1]$ we use the $\LIP$ estimate. This gives us
    \begin{equation*}
        \frac{1}{\abs{\boo'}}\sum_{i = 0}^{n-1} \bigl(f\bigl(\sigma^i(\bz)\bigr) - f\bigl(\sigma^i(\bz')\bigr)\bigr) - \frac{(2(n-k+m-1)a_{m+1} + a_k)\LIP(f)}{\abs{\boo'}(1-\theta)}.
    \end{equation*}
    Next, we do another $\LIP$ estimate to replace $\bz$ and $\bz'$ by $\overline\bp$ and $\overline\bp'$:
    \begin{equation*}
        \frac{1}{\abs{\boo'}}\sum_{i = 0}^{n-1} \bigl(f\bigl(\sigma^i(\overline{\bp})\bigr) - f\bigl(\sigma^i(\overline{\bp}')\bigl)\bigr) - \frac{(2(n-k+m-1)a_{m+1} + 2a_k)\LIP(f)}{\abs{\boo'}(1-\theta)}.
    \end{equation*}
    Now use \eqref{eq:nkm} and the definition of $\mpe^{\sharp}$:
    \begin{equation*}
        \frac{n}{\abs{\boo'}}\mpe^{\sharp}(\bp, \bp', f) - \frac{4\LIP(f)\cdot a_k}{\abs{\boo'}(1-\theta)}.
    \end{equation*}
    By \eqref{eq:betasharppath} this is at least $4\LIP(f)\cdot a_{k} \cdot (1-\theta)^{-2} \big/\abs{\boo'}$. Since this is nonzero, $\boo \ne \boo'$, so $\mpe_m^{\sharp}(\boo, \boo', f)$ is defined. At the same time,
    \begin{equation*}
        \mu_{\boo'}([\Sigma_m^c(\boo)]) 
        \le  (n-k+m-1)/\abs{\boo'}
        < (n+m)/\abs{\boo'}.
    \end{equation*}
    Then, $\mpe_m^{\sharp}(\boo, \boo', f) > 4\LIP(f)\cdot a_{k} \cdot (1-\theta)^{-2} \big/(n+m)$. Applying \eqref{eq:nkm} gives us the expression in Definition~\ref{def:advantagedorbit}.
\end{proof}
\begin{definition}
    Assume that $\breve{\bp}' \in \Path_{\breve{n}, \breve{\bw}, \breve{\bw}'}$ is some path that contains some subpath $\bp' \in \Path_{n, \bw, \bw'}$ with $\bw,\, \bw' \in \Sigma_k(\bF)$ and assume that $(m, n, k)$ satisfy hypotheses of Lemma~\ref{lem:pathadv}. Then, if there exists $\bp \in \Path_{n, \bw, \bw'}$ with \eqref{eq:betasharppath}, we call $\breve{\bp}'$ \emph{outadvantaged} for $m$. We denote by $\Path_{\breve{n}, \breve{\bw}, \breve{\bw}', [m]}^{\sharp}$ the set of paths in $\Path_{\breve{n}, \breve{\bw}, \breve{\bw}'}$ that are not outadvantaged for $m$. Define
    \begin{equation}        \label{eq:sharpdef}
        \Path_{\breve{n}, \breve{\bw}, \breve{\bw}'}^\sharp \= \bigcup_{m = \breve{n}}^{+\infty} \Path_{\breve{n}, \breve{\bw}, \breve{\bw}', [m]}^\sharp,
    \end{equation}
    and let $\Path_{\breve{n}, \breve{k}}^\sharp$ be the union of that over all $\breve{\bw},\, \breve{\bw}'\in \Sigma_{\breve{k}}(\bF)$. 
    \label{def:outadvpath}
\end{definition}
\begin{lemma}
    A path $\breve{\bp} \in \Path_{\breve{n}, \breve{\bw}, \breve{\bw}'}^{\sharp}$ cannot contain outadvantaged subpaths of shorter length. In particular, if $\bp \in \Path_{n, \bw, \bw'}$ with $n < \breve{n}$ is a subpath of $\breve{\bp}$, then $\bp \in \Path_{n, \bw, \bw'}^{\sharp}$.
    \label{lem:sharpsubpath}
\end{lemma}
\begin{proof}
    If $\bp$ has a low-score subpath with \eqref{eq:betasharppath} then it is also a subpath of $\breve{\bp}$, so $\breve{\bp}$ also has a low-score subpath. We need to be careful to confirm that $m$ satisfies conditions of Lemma~\ref{lem:pathadv}. But since $\breve{n} > n$, in \eqref{eq:sharpdef} for $\Path_{\breve{n}, \breve{\bw}, \breve{\bw}'}^{\sharp}$ we are taking a union over a smaller number of $m$ than for $\Path_{n, \bw, \bw'}^{\sharp}$.
\end{proof}

\section{Splitting the bricks}
\label{sec:split}
Let $f$ be a function with a bounded $\Lip$ norm. From now on we focus on proving Theorem~\ref{thm:main} under the assumption of $M=3$ in \eqref{eq:main.redu}. In particular, we assume that for $n \ge n_\theta$ we have
\begin{equation}    \label{eq:seassumption}
    \theta^{n/3} < Z_n^{-1}. 
\end{equation}
Then, Theorem~\ref{thm:maim.redu} and thus Theorem~\ref{thm:main} for all $\theta<1$ will follow from this proof by the results of Section~\ref{sec:thetareduction}.

We make the problem quantitative: we prove that there is an $\ba$-admissible measure $\nu_j$ such that for $g_j \sim \nu_j$ the function $f+g_j$ has the locking property with arbitrary large probability as $j \rightarrow +\infty$. Then, we convolve all $\nu_j$ into a single measure $\nu$ for which the probability will be $1$.

\begin{definition}
    For $j\in\N$, we say that a Borel probability measure $\nu$ on $\Lip(\Sigma_{\bF}, d_{\ba})$ is \emph{$j$-well chosen} if 
    \begin{itemize}
        \smallskip \item $\nu$ is supported on a Hilbert brick $\cH[\bb_j]$, where $\bb_j$ is $2^{-j}$-$\log$-$\ba$-admissible, and
        \smallskip \item for $g \sim \nu$ and each $f \in \Lip(\Sigma_{\bF}, d_{\ba})$ with $\LIP(f) < j$ we have that $f+g $ has the locking property with $\nu$-probability at least $1 - 2^{-j}$. 
    \end{itemize}
\end{definition}
\begin{theorem}
    Let $\{\nu_j\}_{j \in\N}$ be a sequence of Borel probability measures on $\Lip(\Sigma_{\bF}, d_{\ba})$, where $\nu_j$ is $j$-well chosen for all $j \in \N$. Then the convolved measure $\nu \= \nu_1 \ast \nu_2 \ast \cdots$ is $\ba$-admissible and, for each fixed $f \in \Lip(\Sigma_{\bF}, d_{\ba})$, the function $f+g$ has the locking property with probability $1$ for $g\sim\nu$.
\end{theorem}
\begin{proof}
    Fix $f\in \Lip(\Sigma_{\bF}, d_{\ba})$. 
    Since $\sum_{j=1}^{+\infty}2^{-j}=1<+\infty$, the infinite convolution $\nu$ is $\ba$-admissible by Proposition~\ref{prop:infinite-convolution}.
    Let $s \in \N$ be sufficiently large. We claim that $f+g$ has the locking property with probability at least $1-2^{-s}$. Define
    \begin{equation*}
        \widetilde{\nu}_s \= \nu_1 \ast \cdots \ast \nu_{j-1} \ast \nu_{j+1} \ast \cdots.
    \end{equation*}
    Then, by Proposition~\ref{prop:infinite-convolution}, $\widetilde{\nu}_s$ is $1-2^{-s}$-$\log$-$\ba$-admissible. Hence for $\widetilde{g}_s \sim \widetilde{\nu}_s$ we have $\LIP(\widetilde{g}_s) < 1/4$, so $\LIP(f+\widetilde{g}_s) < s$ as $s$ is sufficiently large.
    So, by definition of $s$-well chosen measures, for every fixed choice of $\widetilde{g}_s$ the function $(f+\widetilde{g}_s)+g_s$ has the locking property with probability of at least $1-2^{-s}$ for $g_s \sim \nu_s$. As $\nu = \nu_s \ast\widetilde{\nu}_s$, the claim follows from Fubini's theorem. 
\end{proof}

Thus the proof of Theorem~\ref{thm:main} reduces to finding $j$-well chosen measures for all $j \in\N$. So, from now on, we focus on the following statement:
\begin{theorem}
    \label{thm:wellchosen}
    Let $j \in \N$. Then there exists a $j$-well chosen measure $\nu_j$.
\end{theorem}

The number $j$ will be fixed from now on. Thus we will drop all the subindices containing $j$ and denote the desired measure by $\nu$ and $g\sim \nu$. Let $N = N(j, \theta, \Sigma_{\bF})$ be some sufficiently large fixed number.  To guarantee that $f+g$ has the locking property, we are going to verify that two auxiliary properties of $f+g$ are likely. That is:
\begin{itemize}
	\smallskip
    \item Property~\textbf{A}: Let $\widehat{\boo} \in \PO_{*}$ be a periodic orbit that maximizes $\mpe(\boo,A_{N}(f+g))$ among the orbits $\boo$ in $\cO_{\Delta,N}$. Then
    $       \abs{\widehat\boo} \le \exp\bigl(N^{1/10}\bigr)$.
    
    \smallskip
    \item Property~\textbf{B}: Let $\widehat{\boo} \in \cO_*$ be a periodic orbit that maximizes $\mpe(\boo, A_{N} (f+g))$ among the orbits $\boo$ in $\cO_{\Delta,N}$. Then 
    $    	\Delta(\widehat{\boo}) > a_{N - \lceil N^{1/5} \rceil}$.
\end{itemize}
The main process of our discussion is to reduce a general case to the case with Property~\textbf{A} with a small loss of probability, then to the case with both Properties~\textbf{A} and~\textbf{B} with a small loss of probability, and at last to the case with the locking property with a small loss of probability.

To verify these properties we will construct three $\ba$-admissible measures: $\nu_{A}$, $\nu_B$, and~$\nu_{C}$.
\begin{theorem}
    There exist  $2^{-j-2}$-$\log$-$\ba$-admissible measures $\nu_A$, $\nu_B$, and $\nu_C$ such that all $f \in \Lip(\Sigma_{\bF},d_{\ba})$ with $\LIP(f) < j + 2$, we have
    \begin{enumerate}
    	\smallskip
    	\item[(i)] For $g_A \sim \nu_A$, a function $f+g_A$ fails Property~\textbf{A} with probability at most $2^{-j-2}$.
    	
    	\smallskip
    	\item[(ii)] For $g_B \sim \nu_B$, a function $f+g_B$ satisfies Property~\textbf{A} but fails Property~\textbf{B} with probability at most $2^{-j-2}$.
    	
    	\smallskip
    	\item[(iii)] For $g_C \sim \nu_C$, a function $f+g_C$ satisfies both Properties~\textbf{A} and~\textbf{B} but fails the locking property with probability at most $ 2^{-j-2}$.
    \end{enumerate}
    \label{thm:prop}
\end{theorem}

The measure $\nu_A$ will be constructed in Section~\ref{sec:A}, $\nu_B$ in Sections~\ref{sec:B}, and $\nu_C$ in Section~\ref{sec:C}. 

\begin{proof}[Proof of Theorem~\ref{thm:wellchosen} from Theorems~\ref{thm:prop}]
    Take $\nu \= \nu_A \ast \nu_B \ast \nu_C$. Then by Proposition~\ref{prop:finite-convolution}, the measure $\nu$ satisfies the admissibility requirement. Fix $f\in \Lip(\Sigma_{\bF}, d_{\ba})$ with $\LIP(f) < j$. Then, with full probability, we have $\LIP(f+g_B + g_C) < j+2$, so $(f+g_B + g_C)+g_A$ fails Property~\textbf{A} with probability at most $2^{-j-2}$ for fixed $g_B$ and $g_C$. Thus, by Fubini's theorem, $f+g$ fails Property~\textbf{A} with probability at most $2^{-j-2}$. Similarly, the probabilities of the two events (a) $f+ g$ satisfies Property~\textbf{A} but not \textbf{B}, (b) $f+g$ has Properties~\textbf{A} and~\textbf{B} but not the locking property both not exceed $2^{-j-2}$. Thus, the probability for $f+g$ to satisfy the locking property is at least $1-3\cdot 2^{-j-2} > 1-2^{-j}$. This concludes this theorem.
\end{proof}

In the last three sections, we construct admissible measures $\nu_A$, $\nu_B$, and $\nu_C$ respectively and prove the three statements in Theorem~\ref{thm:prop}.

\section{Verification of Theorem~\ref{thm:prop}~(i)}
\label{sec:A}
\subsection{Construction of $\nu_A$}
The goal of this section is to prove Theorem~\ref{thm:prop} for Property~\textbf{A}.

To construct $\nu_A$ we first define a $2^{-j-2}$-$\log$-$\ba$-admissible gauge $\{b_{A, \bw}\}_{\bw\in\Sigma_*(\bF)}$ and a sequence $\{b_{A,n}\}_{n\in \N}$ by
\begin{equation}    \label{eq:bbdef}
    b_{A, \bw} = b_{A, \abs{\bw}}\= \begin{cases}
        \frac{2^{-j-2}a_{\abs{\bw}}}{\log \abs{\bw}}&\text{ if }\abs{\bw} >n_\theta\\
        0& \text{ otherwise}.
    \end{cases}
\end{equation}
We define  $\nu_{A, \bw}$ to be a uniform measure on $[-b_{A, \bw}, b_{A, \bw}]$ for $b_{A, \bw} \ne 0$. For $b_{A, \bw} = 0$ the measure $\nu_{A, \bw}$ is defined trivially. Then, the product measure is denoted by $\overline{\nu}_{A} \= \prod_{\bw \in \Sigma_*(\bF)} \nu_{A,\bw}$, and its pushforward by $\nu_A \= \Pi_*(\overline{\nu}_A)$. Note that the choice \eqref{eq:bbdef} satisfies Definition~\ref{def:logadmis} with $C \= 2^{-j-2}$, so it is $2^{-j-2}$-$\log$-$\ba$-admissible, as claimed by Theorem~\ref{thm:prop}.

Next, we prepare a $\delta$-argument. Let $\{\marg_i\}_{i\in \N}$ be a sequence of numbers in $[0,b_{A,i}]$. For all $\bw \in \Sigma_*(\bF)$ with $\abs{\bw} > n_\theta$, we define a cut version of $b_{\bw}$:
\begin{equation*}
     b_{A, \bw}^{(i)} \=
    \begin{cases} b_{A, \bw} - \marg_i &\text{ if } \abs{\bw} = i,\\
    b_{A, \bw}& \text{ otherwise}.
    \end{cases}
\end{equation*}
Then, let $\nu_{A, \bw}^{(i)}$ be the restriction of $\nu_{A, \bw}$ on $\bigl[-b_{A, \bw}^{(i)}, b_{A, \bw}^{(i)}\bigr]$.
\begin{lemma}
    \label{lem:Bdeltaarg}
    For $i > n_\theta$, the cut measure $\nu_{A}^{(i)} \= \Pi_* \bigl( \prod_{\bw \in \Sigma_{i}(\bF)}\nu_{A, \bw}^{(i)} \bigr)$ satisfies 
    \begin{equation}    \label{eq:givenprobB}
    \int\!  \, \mathrm{d}\nu_A^{(i)} 
    \ge 1 - Z_i b_{A,i}^{-1}\marg_i.
    \end{equation}
    \vspace{-\baselineskip}
\end{lemma}
\begin{proof}
    We have
    \begin{align*}
        \int\! \, \mathrm{d}\nu_A^{(i)} 
        = \prod_{\bw \in \Sigma_{i}(\bF)} \int\! \, \mathrm{d}\nu_{A, \bw}^{(i)} 
        = \prod_{\bw \in \Sigma_{i}(\bF)}\biggl(1 - \frac{\marg_i}{b_{A,i}} \biggr) 
        \ge \biggl(1 - \frac{\marg_i}{b_{A,i}} \biggr)^{Z_{i}}
         \ge 1 - Z_i b_{A,i}^{-1} \marg_i.  
    \end{align*}
    This concludes the lemma.
\end{proof}

\subsection{Same path concentration}\label{subs:concA}
\begin{prop}
    Let $n > k > l > n_\theta$ be integers, $\bw, \, \bw' \in \Sigma_{k}(\bF)$, $f\in \Lip(\Sigma_{\bF}, d_{\ba})$ with $\LIP(f) < j+1$, and $g_{A} \sim \nu_A$. For every $\ve>0$, let $E = E_{n,l,\bw,\bw'}$ be the event of $g_A \in \supp \nu_A$ that there if $\mpe(\cdot, f+g_A)$ of two paths $\bp, \, \bp' \in \cP_{n,\bw, \bw'}$ are both $2^{-1}\marg_l n^{-2}$-close to the maxima, then they have the same $l$-concentration. Then we have
    \begin{equation*}
        \nu_A(E_{\ve}) \ge 1-  Z_l b_{A,l}^{-1} \marg_l
    \end{equation*}
    \label{prop:eventbstar'}
\end{prop}
    
    To estimate this probability, we will employ the $\delta$-argument with $\nu_{A}^{(l)}$ from Lemma~\ref{lem:Bdeltaarg}, as \eqref{eq:requiredprobc} coincides with \eqref{eq:givenprobB}. So, our next goal is to define $\cF(g_A)$.

    Let $\widehat{\bp} = \widehat{\bp}(g_A) \in \cP_{n, \bw, \bw'}$ be the path that maximizes the following expression:
    \begin{equation*}
        \widehat{\mpe}(\bp) \= \mpe(\bp, f+g_A) + 2^{-1}  \marg_l \norm{\conc_{l , \bp}}_{2}^2.
    \end{equation*}

    Let $\widehat{\bc} = \widehat{\bc}(g_A) \= \conc_{l, \widehat{\bp}}$. We define $\cF^{(l)}$ by
    \begin{equation}        \label{eq:mathcalFBdef}
        \cF_{\bw}^{(l)}(g_{A, \bw}) \= \begin{cases}
            g_{A, \bw} +  \marg_l \widehat{\bc}[\bw]&\text{ if } \abs{\bw} = l, \\
            g_{A, \bw}& \text{ otherwise}.
        \end{cases}
    \end{equation}
    The following lemma will show that for $\nu_A$-almost surely, $\cF^{(l)}$  is well defined.
    
    \begin{lemma}
        \label{lem:mathcalFBwell}
        Transformation $\cF^{(l)}$ is well defined $\nu_A$-almost everywhere and $\widehat{\bc}(g_A)$ is unique  $\nu_A$-almost everywhere. 
    \end{lemma}
    \begin{proof}
    Note that for $\cF^{(l)}$ to be defined almost everywhere, we just need that $\widehat{\bc}(g_A)$ is unique almost everywhere. If this does not happen, then for two paths $\bp,\, \bp'$ with different $\conc_{l}$, with positive probability
    \begin{equation*}
        \mpe(\bp, f+g_A) - \mpe(\bp', f+g_A) = \marg_l \bigl[\norm{\conc_{l, \bp'}}_2^2 - \norm{\conc_{l, \bp}}_2^2 \bigr].
    \end{equation*}
    The right-hand side is independent of $g_A$ and is thus a constant random variable. Meanwhile, the left-hand side is a linear combination of $g_{A, \bw}$, which is nontrivial as the concentrations are different. Thus the probability is $0$ and we get a contradiction.

    We remark that several combinations of $g_{A, \bw}$ can give the same $g_A$, so we should not associate the Hilbert brick $\cH_A$ of $g_A$ with a direct product $\cI_A$ of segments of $g_{A, \bw}$. Still, in this case \eqref{eq:mathcalFBdef} extends well to $\cF^{(l)}(g_A)$ as the shift amount and $\widehat{\bc}$ do not depend on the choice of $\{g_{A, \bw}\}$ for a given $g_A$.
    \end{proof}

    Note that there are finitely many choices of $\widehat{\bp}(g_A)$, and that $\cF^{(l)} \bigl( \supp \nu_A^{(l)} \bigr) \subseteq \supp \nu_A$, because $\widehat{\bc}[\bw] \in [0,1]$ for every $\bw \in \Sigma_{l}(\bF)$.
    So, using Lemma~\ref{lem:Bdeltaarg}, to prove Proposition~\ref{prop:eventbstar} we only require two things:
    \begin{itemize}
        \item $\cF^{(l)}$ is a bijection on its image $\nu_{A}^{(l)}$-almost everywhere.
        \item For $\nu_{A}^{(l)}$-almost every $g_A$ the event $E$ happens for $f+\cF^{(l)}(g_A)$.
    \end{itemize}
    \begin{lemma}
        For $\nu_{A}^{(l)}$-almost every $g_A$, the event $E$ happens for $f+\cF^{(l)}(g_A)$.
        \label{lem:modifiedbbstar}
    \end{lemma}
    \begin{proof}
        Let $\bp' \in \Path_{n, \bw, \bw'}$ have a different concentration, compared to $\widehat{\bp}$. Call it $\bc' \= \conc_{l, \bp'}$. Then 
        \begin{equation}\label{eq:conc}
            \begin{aligned}
                \mpe\bigl(\bp', f+\cF^{(l)}(g_A)\bigr) 
                &= \mpe(\bp', f+g_A) + \marg_l\langle \bc', \widehat{\bc}\rangle  \\ 
                &= \widehat{\mpe}(\bp', f+g_A) - 2^{-1}\marg_l \norm{\bc'}_2^2 + \marg_l\langle \bc', \widehat{\bc}\rangle \\
                &\le  \widehat{\mpe}(\widehat{\bp}, f+g_A) - 2^{-1}\marg_l \norm{\bc'}_2^2 + \marg_l\langle \bc', \widehat{\bc}\rangle \\
                &= \mpe\bigl(\widehat{\bp},f+\cF^{(l)}(g_A)\bigr) - 2^{-1}\marg_l\norm{\bc' - \widehat{\bc}}_2^2,
            \end{aligned}
        \end{equation}
        where we have used \eqref{eq:mathcalFBdef} and the fact that $\widehat{\bp}$ is maximizing. Now note that $\widehat{\bc}$ and $\bc'$ are bounded away from one another as they are both actual concentrations of paths of length $n$. Thus, $\bc' - \widehat{\bc}$ has at least one element of type $k/n$ with $k \ne 0$ an integer. So,
        \begin{equation*}
            \norm{\bc' -\widehat{\bc}}_2^2 \ge n^{-2}.
        \end{equation*}
        Thus, the event $E$ happens for $f+ \cF^{(l)}(g)$.
    \end{proof}
    \begin{lemma}
        $\cF^{(l)}$ is a bijection on its image $\nu_{A}^{(l)}$-almost everywhere.
        \label{lem:bbijection}
    \end{lemma}
    \begin{proof}
        Assume that $\cF^{(l)}$ is not injective at some point: $g_{A}' \ne g_{A}''$, but $\cF^{(l)}(g_A') = \cF^{(l)}(g_A'')$. Here we assume that $\widehat{\bc}(g_A')$ and $\widehat{\bc}(g_A'')$ are both uniquely chosen. By (\ref{eq:mathcalFBdef}), we have $\widehat{\bc}(g_A') \neq \widehat{\bc}(g_A'')$. By (\ref{eq:conc}), every path maximizing $\mpe(\bp, f+ \cF^{(l)}(g_A'))$ have the same concentration, and both $\widehat{\bp}(g_A')$ and $\widehat{\bp}(g_A'')$ maximize $\mpe(\bp, f+ \cF^{(l)}(g_A'))$. Thus, $\widehat{\bc}(g_A') = \widehat{\bc}(g_A'')$, leading to a contradiction. At last, the $\nu_A^{(l)}$-almost everywhere statement follows from Lemma~\ref{lem:mathcalFBwell}.
    \end{proof}
    Thus we have proven Proposition~\ref{prop:eventbstar'}. We remark that there are other ways of choosing $\widehat{\bp}$ that use different $\widehat{\beta}$ and $\cF$.

    Let $r \ge 1$ be some integer, quantified later, and let $\{n_i\}_{i = 0}^r$ be a sequence of increasing positive integers. We require that $n_0 > n_\theta$ and that $n_i \= K^in_0$ for some integer $K$ and all $i\in\N$; also define $k_i \= Ln_i$ for some integer $L$. We consider
\begin{equation}   \label{eq:LKineq}
   16 < L < K < n_0.
\end{equation}
The exact values of $n_0$, $K$, and $L$ will also be chosen later. We will set $N \= n_r$.
Now, we use Proposition~\ref{prop:eventbstar'} to show that under $g_{A} \sim \nu_A$ deformation it is likely that all the paths of length $n_{i+1}$ that are not outadvantaged will have the same concentration $\conc_{n_i}$. We consider the event $\textbf{A}^*$. Namely, for all $i = 0,\, 1,\, \ldots,\, r-1$ and all $\bw,\, \bw' \in \Sigma_{k_i}(\bF)$, all the paths $\bp'$ in $\Path_{n_{i+1}, \bw, \bw'}^{\sharp}(f+g_A)$ share the same $n_i$-concentration $\conc_{n_i, \bp'}$.
    \begin{lemma}
    Assume that $n_0$ is large enough and that
    \begin{equation}        \label{eq:eventbsterreq}
        rn_i^5 \cdot \theta^{Ln_i/4} < 1
    \end{equation}
    holds for every $i \in \{0, \, 1,\, \dots, \, r-1\}$. Let $f \in \Lip(\Sigma_{\bF}, d_{\ba})$ satisfy $\LIP(f) < j+1$ and let $g_A \sim \nu_{A}$. Then, with probability $\ge1 - 2^{-j-2}$ Event $\textbf{A}^*$ happens. 
    \label{prop:eventbstar}
\end{lemma}

We fix $i$, $\bw$, and $\bw'$. We will estimate the probability for this triple and then combine the probabilities over all triples. For this triple, call the respective event $\textbf{A}_{i, \bw, \bw'}^*$. Then, it is enough to prove that
    \begin{equation*}        \label{eq:requiredprobB}
        \nu_A(\textbf{A}_{i, \bw, \bw'}^*) \ge 1 - 2^{-j-2}\cdot Z_{k_i}^{-2}r^{-1}
    \end{equation*}
    as there are at most $Z_{k_i}$ choices for both $\bw$ and $\bw'$, and $r$ choices for $i$. Since $\theta^{k_i/3}<Z_{k_i}^{-1}$ by \eqref{eq:seassumption}, we will just prove that
    \begin{equation}        \label{eq:requiredprobc}
        \nu_A(\textbf{A}_{i, \bw, \bw'}^*) \ge 1 - 2^{-j-2}\theta^{2Ln_i/3}r^{-1}.
    \end{equation}

    \begin{proof}[Proof of Lemma~\ref{prop:eventbstar}]
        For each $n\in \N$, we set
        \begin{equation}
            \marg_n \= 2^{-j-2}  Z_n^{-1} b_{A,n} \theta^{2nL/3} r^{-1}. 
        \end{equation}
        Applying Proposition~\ref{prop:eventbstar'} for $n \= n_{i+1}$, $k \= k_i$, and $l \= n_i$, we have $\nu_A(E) \ge 1- 2^{-j-2} \theta^{2n_iL/3} r^{-1}$. So it is sufficient to show $E \subseteq \textbf{A}^*_{i,\bw, \bw'}$. 

        Let $g_A \in E$ and $F \= f+g_A$. Let $\widehat{\bp} \in \cP_{n_{i+1}, \bw,\bw'}$ be (one of) the path maximizing $\mpe(\cdot, F)$. Below we check the conditions in Lemma~\ref{lem:pathadv} for $n \= n_{i+1}$ and $k \= k_i$. Let $\bp' \in \cP_{n_{i+1}, \bw,\bw'}$ be another path with different $n_i$-concentration from $\widehat{\bp}$. Using the definitions of $b_{A,\bw}$ and $\marg_{n_i}$, we get
        \begin{equation}\label{eq:mpepp'}
            \mpe\bigl(\widehat{\bp}, \bp', F\bigr) 
            \ge \frac{2^{-2j-5}  a_{n_i} \theta^{2Ln_i/3}}{rn_{i+1}^2  Z_{n_i} \log n_i}. 
        \end{equation}
        Firstly, we desire the inequality: $\mpe(\widehat{\bp} , \bp', F) > 8(j+3)a_{k_i} n_{i+1}^{-1} (1-\theta)^{-2}$, coming from \eqref{eq:betasharppath}. Substituting $K$, $L$ and $a_n = \theta^n$ and using (\ref{eq:mpepp'}) reduces the desired inequality into
        \begin{equation*}
              \frac{8(j + 3)r\cdot Kn_i \log n_i}{(1 - \theta)^2} \cdot \theta^{Ln_i} < 2^{-2j-5}Z_{n_i}^{-1}\theta^{(1+2L/3)n_i}.
        \end{equation*}
        The main terms will be the ones containing $\theta^{Ln_i}$. Canceling them and using the upper bound \eqref{eq:LKineq} on $L$ along with \eqref{eq:seassumption} further reduces the desired inequality into
        \begin{equation}            \label{eq:condBtoverify}
            \frac{8(j + 3)r\cdot Kn_i \log n_i}{(1 - \theta)^2} \cdot \theta^{Ln_i/4} < 2^{-2j-5}.
        \end{equation}
        Doing some final bounds that only use that $n_i$ is large enough, we get \eqref{eq:eventbsterreq}, verifying \eqref{eq:condBtoverify}. Secondly, by Definition~\ref{def:outadvpath}, we need to check \eqref{eq:nkm} for all $m \ge n_{i+1}$, which requires the following:
        \begin{equation}            \label{eq:niplusone}
            n_{i+1}+m \le \theta^{-m + k_i}.
        \end{equation}
        Since $n_0$ is large, it is enough to verify \eqref{eq:niplusone} only for $m = n_{i+1}$, giving us
        \begin{equation*}
            2Kn_i \le \theta^{(L-K)n_i} .
        \end{equation*}
        Finally, \eqref{eq:LKineq} reduces this to $2n_i^2 \le \theta^{-n_i}$, which is true as $n_0$ is large.

        Therefore, $g_A \in \textbf{A}^*_{n_{i+1}, \bw, \bw'}$ and this lemma is proved.
    \end{proof}

\subsection{Counting $P^{\sharp}$ paths}\label{subs:counting}
The rest of the section will be devoted to the following proposition.
    \begin{prop}
        In the setting of Lemma~\ref{prop:eventbstar} assume that Event \textbf{A}* happens. Then Property~\textbf{A} holds.
        \label{prop:bstarb}
    \end{prop}
    To start the proof, we introduce
    \begin{equation*}
        \pi_i \= \Absbig{\Path_{n_{i}, k_{i-1}}^{\sharp}(f+g_A)}
    \end{equation*}
    for $i = 1,\, 2,\, \ldots,\, r$, and set $\pi_0 \= q^{n_0}$. We use Proposition~\ref{prop:mbp} to prove the following.
    \begin{lemma}        \label{lem:recursive}
        In the setting of Lemma~\ref{prop:eventbstar} assume that Event \textbf{A}* happens. Then, for all $i = 1,\, 2,\, \ldots,\, r$,
        \begin{equation}            \label{eq:recursive}
            \pi_i \le n_{i}^{K+2}\cdot \pi_{i-1}^{2L}.
        \end{equation}
        \vspace{-\baselineskip}
    \end{lemma}
    \begin{proof}
        Let $\bw, \, \bw' \in \Sigma_{k_{i-1}}(\bF)$. We bound the size of $\Path_{n_{i}, \bw, \bw'}^{\sharp}(f+g_A)$. Let $\overline{\bw} \in \Sigma_{n_{i-1}-1}(\bF)$ be the $(n_{i-1}-1)$-suffix of $\bw$ and let $\overline{\bw}'$ be the $(n_{i-1}-1)$-prefix of $\bw'$. Then, there is an injective projection
        \begin{equation*}
            \Pi_{\Path}\:\Path_{n_{i}, \bw, \bw'}^{\sharp}(f+g_A)\rightarrow \Path_{n_{i}, \overline{\bw}, \overline{\bw}'}
        \end{equation*}
        that removes the prefix of $\bw$ and the suffix of $\bw'$ from paths. The map $\Pi_{\Path}$ changes $\conc_{n_i}$, but if $\bp$ and $\bp'$ shared the same $\conc_{n_i}$, then $\Pi_{\Path}(\bp)$ and $\Pi_{\Path}(\bp')$ would as well. Since Event \textbf{A}* happens, all the elements of the image of $\Pi_{\Path}$ have the same concentration. Thus, by Proposition~\ref{prop:mbp} with $k \= n_{i-1}-1$ and $n \= n_i - k_{i-1} + n_{i-1}-1$,
        \begin{equation*}
            \Absbig{\Pi_{\Path} \bigl(\Path_{n_{i}, \bw, \bw'}^{\sharp}(f+g_A)\bigr)} \le (n_i - k_{i-1} + 2n_{i-1}-2)^{(n_i - k_{i-1}) (n_{i-1}-1)^{-1}+2}.
        \end{equation*}
        We use \eqref{eq:LKineq} to simplify and remind that $\Pi_{\Path}$ was injective. Hence,
        \begin{equation}            \label{eq:pathsharpsmall1}
            \Absbig{\Path_{n_{i}, \bw, \bw'}^{\sharp}(f+g_A)} \le n_i^{K+2}.
        \end{equation}
        We want to sum \eqref{eq:pathsharpsmall1} over all $\bw$ and $\bw'$. When $i = 1$, there are at most $q^{k_0} = q^{Ln_0}$ choices for both $\bw$ and $\bw'$, giving us \eqref{eq:recursive}. When $i > 1$, some choices of $\bw$ and $\bw'$ will cause $\Path_{n_{i}, \bw, \bw'}^{\sharp}$ to be empty. For example, the first $n_{i-1}+k_{i-2}$ symbols of $\bw$ constitute a subpath of every element of $\Path_{n_{i}, \bw, \bw'}^{\sharp}(f+g_A)$. By Lemma~\ref{lem:sharpsubpath} this subpath must lie in $\Path_{n_{i-1}, k_{i-2}}^{\sharp}(f+g_A)$ for \eqref{eq:pathsharpsmall1} to be nonempty. So, there are at most $\pi_{i-1}$ choices for these symbols of $\bw$. Since $k_{i-1} < L(n_{i-1}+k_{i-2})$, there are at most $\pi_{i-1}^L$ choices for $\bw$.
    \end{proof}
    \begin{proof}[Proof of Proposition~\ref{prop:bstarb}]
        Our first goal is to prove
        \begin{equation}            \label{eq:condbgoal}
            \pi_r \le \exp\bigl( n_r^{1/10}\bigr).
        \end{equation}
        We apply Lemma~\ref{lem:recursive} and rewrite \eqref{eq:recursive} in $\log$ form:
        \begin{equation*}
            \log \pi_i \le (K+2)\log n_i + 2L\log \pi_{i-1}.
        \end{equation*}
        We apply this $r$ times to bound $\pi_r$:
        \begin{equation*}
            \log \pi_r \le (K+2)\sum_{i = 1}^r (2L)^{r-i} \log n_i + (2L)^r \log \pi_0.
        \end{equation*}
        We substitute $\pi_0$ and $n_i \= K^i n_0$:
        \begin{equation*}
            \log \pi_r \le (K+2)\sum_{i = 1}^r (2L)^{r-i} (i\log K + \log n_0) + (2L)^r n_0\log q.
        \end{equation*}
        Evaluating and bounding the sum in $i$ gives
        \begin{equation*}
            \log \pi_r \le (K+2)(2L)^r(\log n_0+\log K) + (2L)^r n_0\log q.
        \end{equation*}
        We also use $\log K \le \log n_0$. Meanwhile,
        \begin{equation*}
            n_r^{1/10} = K^{r/10} n_0^{1/10}.
        \end{equation*}
        Thus, to prove \eqref{eq:condbgoal}, we need to satisfy two inequality simultaneously:
        \begin{equation*}
            \begin{cases}
                2(K+2)(2L)^r\log n_0 + (2L)^r n_0\log q \le K^{r/10} n_0^{1/10},\\
                rn_0^5 \cdot \theta^{Ln_0/4} < 1 \qquad \text{ (the equation required in Proposition~\ref{prop:eventbstar}).}
            \end{cases}
        \end{equation*}
        Setting $r \= n_0$ to be large and setting $L \ll K^{1/20}\ll n_0$ will satisfy both equations. 

        We have proven \eqref{eq:condbgoal}. Now consider the orbit $\widehat{\boo}$ from Property~\textbf{A}. It maximizes $\mpe(\boo, A_{n_r}(f+g_A))$, so clearly $\widehat{\boo} \in \PO_{[n_r]}^{\sharp}(f+g_A)$. Then by Lemma~\ref{lem:pathadv} for $m = n_r$ and $n \le n_r$ we either have $\abs{\widehat{\boo}} < m+n+k <3m=3n_r$ or $\widehat{\boo}$ contains no $(m, n, k)$-outadvantaged paths. The first alternative immediately gives Property~\textbf{A} whereas the second implies that $\widehat{\boo}$ cannot visit outadvantaged paths in $\Path_{n_r, k_{r-1}}$ by Definition~\ref{def:outadvpath}. Since $\Delta(\widehat{\boo}) > a_{n_r}$ the paths of length $n_r$ that $\widehat{\boo}$ visits cannot repeat. So, we get
        \begin{equation*}
            \abs{\widehat{\boo}} 
            \le \Absbig{\Path_{n_{r}, k_{r-1}}^{\sharp}(f+g_A)}
            = \pi_r \le \exp\bigl( n_r^{1/10}\bigr).
        \end{equation*}
        Thus we have proven Property~\textbf{A} for this $g_A$.
    \end{proof}
    Combining Propositions~\ref{prop:eventbstar} and~\ref{prop:bstarb} gives Theorem~\ref{thm:prop} for Property~\textbf{A}.
   
 \section{Verification of Theorem~\ref{thm:prop}~(ii)} \label{sec:B}
 
 In this section, we mainly prove Theorem~\ref{thm:prop}~(ii), reducing the case with Property~\textbf{A} to that with both Properties~\textbf{A} and~\textbf{B} with a small loss of probability.
 
 Following Section~\ref{sec:A}, we let $N$ be the integer defined in Section~\ref{sec:A}, and $N$ can be chosen large enough. For convenience, we define values 
 \begin{equation}\label{eq:Bvalues}
 	P \= \exp\bigl(N^{1/10}\bigr), \quad 
    s \= \bigl\lceil N^{1/5} \bigr\rceil, \quad 
    R \= \lfloor \log_2 P \rfloor +1, \quad \text{ and } 
    k_i \= N-is
 \end{equation}
 for each integer $1\le i \le R$.
 Clearly, if $N$ is large enough, $Rs<N$, so $\{k_i\}_{i=1}^R$ is well defined, and we can also assume that $k_R = N-sR$ is large enough. Recall that a Lipschitz function $f$ satisfies Property~\textbf{A} if for every orbit $\widehat{\boo} \in \cO_*$ that is maximizing $\mpe(\boo, A_N(f))$ among the ones with $\Delta(\boo) > a_N$, it holds that $\abs{\widehat{\boo}} \le P$. 
 
 \subsection{Construction of $\nu_B$}
 
 Similarly to (\ref{eq:bbdef}), we define
 \begin{equation}\label{eq:bbdef'}
 	b_{k} \= \begin{cases}
 		\frac{2^{-j-2} }{\log k}\cdot a_k &\quad \text{ if } k \in \{k_i\}_{i=1}^R, \\
 		0 &\quad \text{ otherwise}
 	\end{cases}
 \end{equation}
 and $b_{\bw} \= b_{\abs{\bw}}$ for each $\bw\in \Sigma_*(\bF)$.
 As in (\ref{eq:deffk}), we also define $Z_k \= \abs{\Sigma_k(\bF)}$ for every $k\in \N$.
 
 Like in Section~\ref{sec:thetareduction}, we define $\cI_{B,k} \= \prod_{\bw \in \Sigma_k(\bF)} [-b_{\bw}, b_{\bw}]$ for every $k\in \{k_i\}_{i=1}^R$ and $\cI_B \= \prod_{i=1}^R \cI_{B,k_i}$. Define the projection $\Pi \: \cI_B \to \Lip(\Sigma_{\bF}, d_{\ba})$ by
 \begin{equation*}
 	\Pi (\{g_{\bw}\}_{\bw\in \Sigma_*(\bF)}) \= \sum_{\bw\in \Sigma_*(\bF)} g_{\bw} \mathbbold{1}_{[\bw]}
 \end{equation*}
 and $\cH_B \= \Pi(\cI_B)$. By Lemma~\ref{lem:hilbert-brick}, $\cH_B$ is compact. Moreover, for each $k \in \N$, define projections $\Pi_k \: \R^{ \Sigma_k (\bF) }  \to \Lip ( \Sigma_{ \bF },  d_{\ba} )$ by
 \begin{equation*}
 	\Pi_k ( \{g_{\bw}\}_{\bw\in \Sigma_k(\bF)} ) \= \sum_{\bw \in \Sigma_k(\bF)} g_{\bw} \mathbbold{1}_{[\bw]}.
 \end{equation*}
 
 However, in this section we do not consider the full Hilbert brick $\cH_B$, we consider a subset of it. Let $k\in \{k_i\}_{i=1}^R$, we define 
 \begin{equation}\label{eq:deflk}
 	l_k \= \frac{b_k Z_k}{2k(1+\log q )}.
 \end{equation}
 We consider the following polytope
 \begin{equation}\label{eq:defDlk}
 	D_{k, r} \= \Bigl\{ (g_{\bw})_{\bw\in \Sigma_k(\bF)} : \sum_{\bw\in \Sigma_k(\bF)} \abs{g_{\bw}} \le r \Bigr\}
 \end{equation}
 for every $r>0$ and
 \begin{equation}\label{eq:def.omegak}
 	\Omega_k \= \cI_{B,k} \cap D_{k,l_k}
 \end{equation}
 for every $k\in \{k_i\}_{i=1}^R$. For a Borel set $X$ in $\R^d$ for some $d\in \N$ and another integer $0\le c \le d$, we denote the $c$-dimensional Lebesgue measure of $X$ by $V_c(X)$. For every $k\in \{k_i\}_{i=1}^R$, let $\nu_{B,k}$ be the normalized Lebesgue measure on $\Omega_k$, more precisely, for every Borel set $E \subset \R^{Z_k}$, 
 \begin{equation}
 	\nu_{B,k} (E) \= \frac{ V_{Z_k}(E \cap \Omega_k)}{ V_{Z_k} (\Omega_k)}.
 \end{equation}
 Define $\Omega \= \prod_{i=1}^R \Omega_{k_i}$, $\cH_{B}' \= \Pi(\Omega) \subseteq \cH_{B}$, $\overline{\nu}_{B} \= \prod_{i=1}^R \nu_{B,k}$, and $\nu_{B} \= \Pi_{*}(\overline{\nu}_B)$. Clearly, $\supp \nu_{B} = \cH_{B}' \subseteq \cH_B$, so $\nu_{B}$ is compactly supported.
 
 For every $k\in \{k_i\}_{i=1}^R$, let $\nu_k$ be the normalized Lebesgue measure on $D_{k,l_k}$, we establish the following computational lemma to reduce the probability estimation of $\nu_{B,k}$ to that of $\nu_k$:
 
 \begin{lemma}\label{lem:probreduce}
 	Let $k\in \{k_i\}_{i=1}^R$, if $k_R$ is large enough, then for every Borel set $E \subseteq \R^{Z_k}$, we have $\nu_{B,k}(E) \le 2 \nu_{k}(E)$.
 \end{lemma}
 
 \begin{proof}
 	Note that for every $\bw_0 \in \Sigma_k(\bF)$, we have
 	\begin{align*}
 		\{ \bg \in D_{k,l_k} : g_{\bw_0} > b_k \} - b_k  \be_{\bw_0} &= \{ \bg \in D_{k,l_k-b_k} : g_{\bw_0} >0 \} \quad \text{and} \\
 		\{ \bg \in D_{k,l_k} : g_{\bw_0} < -b_k \} + b_k \be_{\bw_0} &= \{ \bg \in D_{k,l_k-b_k} : g_{\bw_0} <0 \},
 	\end{align*}
 	where $\be_{\bw_0} \= (0,  0,  \dots,  1,  \dots,  0)$ with the only ``$1$'' at the $\bw_0$-th coordinate. Thus,
 	\begin{equation}\label{eq:xw0large}
 		\nu_k (\{ \bg \in D_{k,l_k} : \abs{g_{\bw_0}} > b_k \}) = \frac{V_{Z_k}(D_{k,l_k-b_k})}{V_{Z_k}( D_{k,l_k})} = \biggl( 1- \frac{b_k}{l_k} \biggr)^{Z_k} \le  \exp\biggl( -\frac{Z_k b_k}{l_k} \biggr).
 	\end{equation}
 	Summing up (\ref{eq:xw0large}) over $\bw_0 \in \Sigma_{k}(\bF)$, using (\ref{eq:deflk}) and the fact that $Z_k \le q^k$, we obtain that
 	\begin{equation*}
 		\nu_k(\Omega_k) 
        \ge 1- Z_k \exp ( - Z_k b_k/l_k ) 
        \ge 1- Z_k\exp (-2(1+\log q )k) 
        \ge 1- \exp(-(2+\log q )k ).
 	\end{equation*}
 	Since $k_R$ can be chosen large enough and $k \ge k_R$, we have $\nu_k(\Omega_k) \ge 1- \exp(-(2+\log q )k ) > 1/2$. This completes this lemma.
 \end{proof}
 
 \subsection{Concentration based on Property~\textbf{A}} \label{subsec:concentrationbasedon}
 
 In this section, we mainly use the concentration idea developed in Subsection~\ref{subs:concA} based on Property~\textbf{A} to prove Theorem~\ref{thm:prop}~(ii).
 
 We denote the set of periodic orbits $\boo$ with $\Delta(\boo) > a_N$ by $\cO_{\Delta, N}$.
 Fix $k\in \{k_i\}_{i=1}^R$ and $f\in \Lip(\Sigma_{\bF}, d_{\ba})$. Analogously to Subsection~\ref{subs:concA}, for each $r>0$ and periodic orbit $\boo \in \cO_{\Delta,N}$, we define  
 \begin{equation}
 	\widehat{\mpe}_r (\boo) \= \mpe( \boo, A_N f ) + 2^{-1}r \norm{\conc_{k,\boo}}_2^2.
 \end{equation}
 Let $ \widehat{\boo} = \widehat{\boo}(f,r)$ be  the pair maximizing $\widehat{\beta}_r (\boo)$ among orbits in $\cO_{\Delta, N}$ and $\widehat{\bc} = \widehat{\bc} (f,r) \= \conc_{k,\widehat{\boo}}$. Moreover, for $\bg \in \R^{Z_k}$, define $\widehat{\boo}_{\bg} \= \widehat{\boo}(f+ \Pi_k(\bg) ,r)$ and $\widehat{\bc}_{\bg} \= \widehat{\bc} (f+\Pi_k(\bg), r)$.
 Then, we define the translation  $\cF_r \: \R^{Z_k} \to \R^{Z_k} $  by 
 \begin{equation}\label{eq:Btrans}
 	\cF_{r} (\bg) \= 
 		\bg + r \widehat{\bc}_{\bg}.
 \end{equation}

 Following an analogous argument to Lemmas~\ref{lem:modifiedbbstar} and~\ref{lem:bbijection}, we obtain the following two statements:
 \begin{enumerate}
 	\smallskip
 	\item[\textbf{B}(i)] $\cF_r$ is a bijection to its image $\nu_k$-almost everywhere. 
 	
 	\smallskip
 	\item[\textbf{B}(ii)] For every $\bg \in \R^{Z_k}$, and periodic orbit $\boo \in \cO_{\Delta, N} \smallsetminus \{\widehat{\boo}_{\bg}\}$, we have
 	\begin{equation*}
 		\beta(A_N f + \Pi_k(\cF_r(\bg)) , \conc_{k,\boo}) 
        \le \beta \bigl(A_N f + \Pi_k(\cF_r(\bg)) , \conc_{k,\widehat{\boo}_{\bg}} \bigr) - 2^{-1} r \norm{\conc_{k,\boo}-\conc_{k,\widehat{\boo}_{\bg}}}^2.
 	\end{equation*}
 \end{enumerate}
 
 Next, we use the above concentration results to prove the following lemma:
 \begin{lemma}\label{lem:concB}
 	Let $k\in \{k_i\}_{i=1}^R$ and $f\in \Lip(\Sigma_{\bF}, d_{\ba}).$ For every $\ve>0$, let $E_{\ve}$ be the event of $\bg \in D_{k,l_k}$ that there exist $\boo_1, \, \boo_2 \in \cO_{\Delta, N}$ with $\abs{\boo_1}, \, \abs{\boo_2} \le P$, $\conc_{k,\boo_1} \neq \conc_{k, \boo_2} $, and $\mpe(A_N f + \Pi(\bg) , \boo_1)$ and $\mpe(A_N f + \Pi(\bg), \boo_2)$ are both $\ve$-close to the maximum of $\mpe(A_N f + \Pi(\bg), \,\cdot\,)$ among $\boo \in \cO_{\Delta,N}$. Then 
 	\begin{equation*}
 		\nu_k ( E_{\ve} ) \le  8(1+\log q ) k \ve P^4/b_k .
 	\end{equation*}
 \end{lemma}
 
 \begin{proof}
 	Let $r \= 4\ve P^4 $. For every $\boo \in \cO_{\Delta,N}$, let $S(\boo)$ be the event of $\bg \in D_{k,l_k-r}$ such that $\boo$ is the unique orbit in $\cO_{\Delta,N}$ maximizing $\beta (A_N f + \Pi(\bg), \boo)$. An argument analogous to Lemma~\ref{lem:mathcalFBwell} gives $\nu_k (D_{k,l_k-r}) = \sum_{\boo \in \cO_{\Delta,N}} \nu_k(S(\boo))$. Since $\norm{\bc}_1$ for every concentration $\bc$, (\ref{eq:Btrans}) implies that for every $\boo \in \cO_{\Delta,N}$, we have $\cF_r (S(\boo)) \subseteq D_{k,l_k}$. By translation invariance, $\nu_k (S(\boo)) = \nu_k (\cF_r(\boo))$. Using those above and \textbf{B}(i), we have $\nu_k (D_{k,l_k-r}) = \sum_{\boo \in \cO_{\Delta,N}} \nu_k(\cF_r(S(\boo)))$.
 	
 	Let $\bg \in \cF_r(S(\boo) )$ for some $\boo \in \cO_{\Delta,N}$. We are going to prove that $\bg \notin E_{\ve}$ and argue by contradiction. Suppose $g\in E_{\ve}$. Let $\boo_1 $ and $\boo_2$ be the two orbits as in the statement. Conditions $\abs{\boo_1}, \, \abs{\boo_2} \le P$ and $\conc_{k,\boo_1} \neq \conc_{k,\boo_2}$ show that there exists $\bw \in \Sigma_k(\bF)$ such that $\abs{(\conc_{k,\boo_1})_{\bw} - (\conc_{k,\boo_2})_{\bw}} \ge P^{-2}$. This implies $\norm{\conc_{k,\boo_1} - \conc_{k,\boo_2}}_2^2 \ge P^{-4}$. However, since $\boo_1$ and $\boo_2$ are $\ve$-close to the maximum, \textbf{B}~(ii) and $r=4 \ve P_4$ show that $\norm{\conc_{k,\boo_1} - \conc_{k,\boo_2}}^2 \le 2^{-1}P^{-4}$, yielding a contradiction.
 	
 	Finally, we have 
 	\begin{equation*}
 		\begin{aligned}
 			\nu_k(E_{\ve}) 
            \le 1-\nu_k(D_{k,l_k-r})
            = 1- \bigl( 1- 4\ve P^4 l_k^{-1} \bigr)^{Z_k} 
            \le  4Z_k \ve P^4 l_k^{-1} 
            =  8(1+\log q ) k \ve P^4/b_k .
 		\end{aligned}
 	\end{equation*} 
 	where the last equality above follows from (\ref{eq:deflk}).
 	Therefore, this lemma is true.
 \end{proof}
 
 Using Lemma~\ref{lem:probreduce}, the Fubini theorem, and the push-forward $\nu_B \= \Pi_{*}(\nu_B')$, we conclude from Lemma~\ref{lem:concB} that
 
 \begin{cor}\label{cor:concB}
 	Let $k\in \{k_i\}_{i=1}^R$ and $f\in \Lip(\Sigma_{\bF}, \ba).$ For every $\ve>0$, let $F_{\ve,k}$ be the event of $g \in \cH_{B}'$ that there exist $\boo_1, \, \boo_2 \in \cO_{\Delta,N}$ with $\abs{\boo_1}, \, \abs{\boo_2} \le P$, $\conc_{k,\boo_1} \neq \conc_{k, \boo_2} $, and $\mpe(A_N f + g , \boo_1)$ and $\mpe(A_N f + g, \boo_2)$ are both $\ve$-close to the maximum of $\mpe(A_N f + g, \,\cdot\,)$ among $\boo \in \cO_{\Delta,N}$. Then 
 	\begin{equation*}
 		\nu_B ( F_{\ve,k} ) 
        \le  16(1+\log q ) k \ve P^4/b_k .
 	\end{equation*} 
 \end{cor}
 
 \begin{remark}
 	One might ask why we ``change the shape'', in other words, use $\Omega_k$ instead of $\cI_{B,k}$ to estimate probability. Comparing the result, one can try to use \textbf{B}~(i) and (ii) and an analogous argument to Section~\ref{sec:thetareduction} or~\ref{sec:A} and find the estimation in Corollary~\ref{cor:concB} will be $\operatorname{O}\bigl( \frac{\ve P^4 Z_k}{b_k} \bigr)$, which is more rough than the result in Corollary~\ref{cor:concB}, since $Z_k$ is far larger than $k$. Below we explain the main reason why the $\Omega_k$ do better. Indeed, essentially we estimate probability on $D_{k,l_k}$, a ball of $l^1$ norm, while $\cI_{B,k}$ is a ball of $l^\infty$ norm. The estimation using $\cI_{B,k}$ depends on the $l^\infty$ norm of the perturbation term, while the estimation using $D_{k,l_k}$ depends on the $l^1$ norm. The concentration arguments use the perturbation term $\widehat{\bc}$, with both $l^1$ and $l^\infty$ norms not exceed $1$. But for most of cases, the $l^\infty$ norm of $\widehat{\bc}$ is far smaller than $1$ and the $l^1$ norm always equals to $1$, so using $\cI_{B,k}$ causes a large probability loss. To give a comparison, we consider the estimation in Section~\ref{sec:thetareduction}. The perturbation term in this section the perturbation term is $\sum_{\bw \in \Sigma_m(\boo)}\mathbbold{1}_{[\bw]}$ (see the definition of $\cF$ below (\ref{eq:probeq})), whose $l^\infty$ norm always equals to $1$ and $l^1$ norm is upper bounded by $Z_m = \abs{\Sigma_m(\bF)}$. But for most orbits $\boo$, this $l^1$ norm is far smaller than the upper-bound, so using $l^\infty$ norm and a a space like $\cI_{B,m}$ is a more sensible choice. 
 \end{remark}
 
 \subsection{Verification of Property~\textbf{B}}
 
 At last, we use the concentration result Corollary~\ref{cor:concB} to show Property~\textbf{B} with a small loss of probability based on Property~\textbf{A}. 
 
 We first use the condition in Corollary~\ref{cor:concB} and the recursively half-cutting periodic orbits idea according to \cite[Proposition~2.1]{HLMXZ25} to give a sufficient condition of Property~\textbf{B} based on Property~\textbf{A}.

 Fix $f\in \Lip(\Sigma_{\bF}, d_{\ba})$ with $\LIP(f) \le j+2$. For integer $2 \le i \le R$, let $F_{\ve_i, k_i}$ be the event in Corollary~\ref{cor:concB}, define 
 \begin{equation}\label{eq:B2value}
 	\ve_i \= 2KP a_{k_{i-1} -1} \quad \text{ and } \quad S_i \= \cH_{B}' \smallsetminus F_{\ve_i, k_i}, \quad\text{ where }\quad K \=  2(j+3) \theta (1-\theta)^{-1}.
 \end{equation}
 
 \begin{prop}\label{prop:halfcut}
 	Suppose $k_R \ge n_\theta +1$, $s\ge 2$, $\theta^s \le 1/2$, and $g\in \cH_{B}'$ with $f + g \in \bigcap_{i=2}^R S_i$. Then if there exists $\boo_0 \in \cO_{\Delta,N}$ such that $\abs{\boo_0} \le P$ and $\boo_0$ maximize $\mpe(A_N f + g , \,\cdot\,)$ among orbits in $\cO_{\Delta,N}$, then $\boo_0 \in \cO_{\Delta,k_1}$.
 \end{prop}
 
 \begin{proof}
 	We first note that for every $g \in \cH_B'$, straightforward calculations give $\LIP(g) \le 1/4<1$. So we have $\Lip(f+g) \le j+3$. Then, we argue by contradiction and assume that $\boo_0 \notin \cO_{\Delta, k_1}$, so $\Delta(\boo_0) \le a_{k_1} < a_{k_1 -1}$. Applying Lemma~\ref{lemma.gapcompose} to $\boo_0$, we obtain several periodic orbits satisfying (\ref{eq:shadowing1}), and let $\boo_1$ be the shortest one of them. Thus, $\Delta(\boo_1) > a_{k_1}$, and (\ref{eq:shadowing1}) and that $\boo_0$ maximize $\mpe(A_N f + g, \,\cdot\,)$ among orbits in $\cO_{\Delta,N}$ give 
 	\begin{equation}\label{eq:cut1}
 		\abs{\boo_1} (\mpe(A_Nf +g , \boo_0) - \mpe(A_N f+g, \boo_1)) \le  2(j+3) \theta (1-\theta)^{-1}a_{k_1-1} .
 	\end{equation}
 	
 	Letting the construction of $\boo_1$ above as the base step, we recursively construct $\{\boo_i\}_{i=1}^{R'}$ for some integer $1\le R' \le R$. Assume that $\boo_{i-1}$ is constructed.
 	If $\boo_{i-1} \notin \cO_{\Delta,k_i}$, then $\Delta(\boo_{i-1}) \le a_{k_i} < a_{k_i-1}$. We continue the induction and apply Lemma~\ref{lemma.gapcompose} to $\boo_{i-1}$ and let $\boo_i$ be the shortest orbit obtained. Thus, $\abs{\boo_i} \le 2^{-i} P$. Letting $L_t \= \mpe(A_N f+g, \boo_0) - \mpe (A_N f+g, \boo_t)$, using (\ref{eq:shadowing1}), we obtain 
 	\begin{equation}\label{eq:cut2}
 		\abs{\boo_i} L_i 
        \le \abs{\boo_{i-1}} L_{i-1} +  2\abs{ \boo_{i-1} } (j+3) \theta (1-\theta)^{-1}a_{k_i-1} 
        = \abs{\boo_{i-1}} L_{i-1} + K\abs{\boo_{i-1}} a_{k_i-1}.
 	\end{equation}
 	Combining (\ref{eq:cut1}) and iteration of (\ref{eq:cut2}), we obtain
 	\begin{equation}\label{eq:cut3}
 		L_i \le \abs{\boo_i} L_i \le KP \sum_{t=1}^i a_{k_t-1} \le KP a_{k_i-1} \sum_{t=1}^{i} \theta^{(i-t) s} \le \frac{KP a_{k_i-1}}{1-\theta^s} \le  2KP a_{k_i-1} = \ve_{i+1},
 	\end{equation}
 	where the last but two inequality follows from $\theta^s \le  1/2$ and the last one follows from (\ref{eq:B2value}).
 	If $\boo_{i-1} \in \cO_{\Delta,k_i}$, our assumption $f+g \notin F_{\ve_i, k_i}$ and (\ref{eq:cut3}) for $\boo_{i-1}$ give $\conc_{k_i , \boo_{i-1}} = \conc_{k_i, \boo_0}$. Using Lemma~\ref{lem:sameconc}, we have $\boo_0 = \boo_{i-1} \in \cO_{\Delta,k_i} \subseteq \cO_{\Delta,k_1}$. This leads to a contradiction to $\boo_0 \notin \cO_{\Delta,k_1}$. So the induction will continue till $i=R$ and we obtain the orbit $\boo_R$. However, by (\ref{eq:Bvalues}), we have $1\le \abs{\boo_R} \le 2^{-R} P <1$, leading to a contradiction. Therefore, this proposition is true. 
 \end{proof}
 
 \begin{proof}[Proof of Theorem~\ref{thm:prop}~(ii)]
 	Let $f \in \Lip(\Sigma_{\bF} , d_{\ba})$ with $\LIP(f) \le j+2$ and $g \sim \nu_B$. Note that assumptions $k_R \ge n_\theta+1$, $s\ge 2$, and $\theta^s \le 1/2$ in Proposition~\ref{prop:halfcut} are satisfied if $N$ is large enough. Proposition~\ref{prop:halfcut} gives that if $f+g$ satisfies Property~\textbf{A} but not Property~\textbf{B}, then $f+g \in \bigcup_{i=2}^R F_{\ve_i, k_i}$. Corollary~\ref{cor:concB},~(\ref{eq:bbdef'}), and~(\ref{eq:B2value}) give
 	\begin{align*}
 		\sum_{i=2}^{R} \nu_B(F_{\ve_i, k_i}) 
        \le \sum_{i=2}^{R} \frac{16(1+\log q ) k_i \ve_i P^4}{b_{k_i}} 
 		&= 2^{j+2} \cdot32(1+\log q ) KP^5 \sum_{i=2}^{R} \biggl(  \frac{a_{k_{i-1} -1}}{a_{k_i}} k_i \log k_i \biggr) \\
 		& \le 2^{j+2} \cdot32(1+\log q ) KP^5 (R-1) \theta^{s-1} N \log N.
 	\end{align*}
 	By (\ref{eq:Bvalues}), the logarithm of the last term above is 
 	\begin{equation*}
 		\operatorname{O}(\log N) + 5 N^{1/10}  +   \bigl\lceil N^{1/5} -1 \bigr\rceil \log \theta,
 	\end{equation*}
 	which tends to $-\infty$ when $N$ tends to $+\infty$. Thus, $N$ can be chosen sufficiently large such that $\sum_{i=2}^{R} \nu_B(F_{\ve_i, k_i}) \le 2^{-j-2}$. Therefore, Theorem~\ref{thm:prop}~(ii) is verified.
 \end{proof}

 \section{Verification of Theorem~\ref{thm:prop}~(iii)} \label{sec:C}
 
 In this section, we mainly prove Theorem~\ref{thm:prop}~(iii), reducing the case with both Properties~\textbf{A} and~\textbf{B} to that with the locking property with a small loss of probability.
 
 We use the integer $N$ defined in Section~\ref{sec:A} which is large enough and the values defined in (\ref{eq:Bvalues}) and (\ref{eq:deflk}), that is, $P = \exp \bigl( N^{1/10} \bigr)$, $s = \bigl\lceil N^{1/5} \bigr\rceil$, $k \= k_1 = N-s$, and $l_k \= \frac{b_k Z_k}{2k(1+\log q )}$. In this section, we only use the polytope $\Omega_k$ defined in (\ref{eq:def.omegak}) to estimate probability. Following Section~\ref{sec:B}, let $\nu_{k}$ be the normalized Lebesgue measure on $D_{k,l_k}$, $\overline{\nu}_{C}$ be the normalized Lebesgue measure on $\Omega_k$, $\cH_{C} \= \Pi(\Omega_k)$, and $\nu_C \= \Pi_{*}(\nu_C')$.
 
 For a function $f\in \Lip(\Sigma_{\bF}, d_{\ba})$, let $\widehat{\boo} = \widehat{\boo}(f)$ be the orbit maximizing $\mpe(A_N f, \,\cdot\,)$ among orbits on $\cO_{\Delta,N}$. To make $f$ satisfying the locking property, we are going to use the sufficient condition in Lemma~\ref{lemma.betastar}, that is 
 \begin{equation}\label{eq:Csufcon}
 	\mpe^{\sharp}_N (f, \widehat{\boo}, \boo) \ge 4 \LIP(f) a_N (1-\theta)^{-2} 
 \end{equation}
 for every $\boo \in \cO_{\Delta,N} \smallsetminus \{\widehat{\boo}\}$. Note that $\widehat{\boo} \in \cO_{\Delta,N}$, by Lemma~\ref{lem:sameconc}, $\mpe_N^\sharp (f, \widehat{\boo}, \boo)$ is always well defined. For $g \sim \nu_B$, we mainly use the $\delta$-argument like in \cite[Proposition~5.1]{DLZ24} and in Section~\ref{sec:thetareduction}.
 
 Fix $f\in \Lip(\Sigma_{\bF}, d_{\ba})$. Let $h>0$, $\bg \in D_{k,l_k}$, and $\widehat{\boo}_{\bg}$ be the orbit maximizing $\mpe(A_N f + \Pi(\bg) , \,\cdot\,)$ among orbits in $\cO_{\Delta,N}$. Likewise the translation defined below (\ref{eq:probeq}), we define $\cF_h \: \R^{Z_k} \to \R^{Z_k}$ by
 \begin{equation}\label{eq:Ctrans}
 	(\cF_h)_{\bw} (\bg) \= \begin{cases}
 		g_{\bw} + \frac{h}{\abs{\widehat{\boo}_{\bg}}} \quad &\text{ if } \bw \in \Sigma_k(\widehat{\boo}_{\bg}), \\
 		g_{\bw} \quad & \text{ otherwise}.
 	\end{cases}
 \end{equation}
 Analogous to Lemma~\ref{lem:sN.fullprob}, $\widehat{\boo}_{\bg}$ is unique and $\cF_h$ is well defined $\nu_k$-almost surely. 
 
 \begin{lemma}\label{kem:Clem}
 	For the translation $\cF_h$, the following hold:
 	\begin{enumerate}
 		\smallskip
 		\item[(i)] $\cF_h (D_{k,l_k-h}) \subseteq D_{k,l_k}$.
 		
 		\smallskip
 		\item[(ii)] $\cF_h$ is a $\nu_k$-almost surely bijection to its image.
 		
 		\smallskip
 		\item[(iii)] $f+ \Pi(\bg)$ satisfies Properties~\textbf{A} and~\textbf{B} if and only if $f + \Pi(\cF_h(\bg))$ satisfies Properties~\textbf{A} and~\textbf{B}. Moreover, if $f+ \Pi(\bg)$ satisfies Properties~\textbf{A} and~\textbf{B}, then for each $\boo \in \cO_{\Delta,N}$,
 		\begin{equation*}
 			\mpe_N^{\sharp} (f+ \Pi(\cF_h(\bg)) , \widehat{\boo}_{\bg}, \boo) \ge h/(P^2(s-1)).
 		\end{equation*} 		
 	\end{enumerate}
 \end{lemma}
 
 \begin{proof}
 	\smallskip
 	(i) By (\ref{eq:Ctrans}), the $l^1$ norm of the perturbation vector is exactly $h$, then (i) follows.
 	
 	\smallskip
 	(ii) For $\nu_k$-almost every $\bg$, $\widehat{\boo}_{\bg}$ is uniquely chosen, and for every $\boo\in \cO_{\Delta,N} \smallsetminus \{\widehat{\boo}_{\bg}\}$, 
    \begin{equation*}
       \mpe(f+ \Pi(\cF_h(\bg)) , \widehat{\boo}_{\bg}) - \mpe(f+ \Pi(\cF_h(\bg)) , \boo) 
       \ge \mpe(f+ \Pi(\bg) , \widehat{\boo}_{\bg}) - \mpe(f+ \Pi(\bg) , \boo) 
       >0,
    \end{equation*}
    so $\widehat{\boo}_{\bg}$ also uniquely maximize $\mpe(f+ \Pi(\cF_h(\bg)) , \,\cdot\,)$ among orbits in $\cO_{\Delta,N}$. Then (ii) follows.
 	
 	\smallskip
 	(iii) Note that $\widehat{\boo}(\Pi(\bg)) = \widehat{\boo} (\Pi(\cF_h(\bg)))$. Both Properties~\textbf{A} and~\textbf{B} holds equivalently for $f+\Pi(\bg)$ and $f + \Pi(\cF_h(\bg))$.
    
    Below we prove the last statement. Let $\boo \in \cO_{\Delta,N} \smallsetminus \{ \widehat{\boo}_{\bg} \}$ and $p \= \abs{\widehat{\boo}_{\bg}}$.  If $\abs{\boo} \ge p$, then $\mpe(A_N f+ \Pi(\cF_h(\bg)), \widehat{\boo}_{\bg}) - \mpe(A_N f+ \Pi(\cF_h(\bg)), \boo) = (h/p) \mu_{\boo} (\Sigma_k^c (\widehat{\boo}_{\bg}) ) \ge \mu_{\boo} (\Sigma_N^c (\widehat{\boo}_{\bg}) ) h /((s-1)p)$. The last inequality can be shown as follows. Since $\widehat{\boo}_{\bg} \in \cO_{\Delta,k}$, by \cite[Lemma~4.3]{BZ16}, we have $\mu_{\boo} (\Sigma_N^c (\widehat{\boo}_{\bg}) ) \le (N-k+1) \mu_{\boo} (\Sigma_k^c (\widehat{\boo}_{\bg}) ) = (s-1)\mu_{\boo} (\Sigma_k^c (\widehat{\boo}_{\bg}) )$. If $\abs{\boo} < p$, since $\widehat{\boo}_{\bg} \in \cO_{\Delta,k}$, Lemma~\ref{lem:sameconc} gives $\mu_{\boo} (\Sigma_k^c(\widehat{\boo}_{\bg})) \ge 1/\abs{\boo} >1/p $. So $\mpe(A_N f+ \Pi(\cF_h(\bg)), \widehat{\boo}_{\bg}) - \mpe(A_N f+ \Pi(\cF_h(\bg)), \boo) = (h/p)  \mu_{\boo} (\Sigma_k^c (\widehat{\boo}_{\bg}) ) > h/p^2$. In both cases, 
 	\begin{equation*}
 		\mpe(A_N f+ \Pi(\cF_h(\bg)), \widehat{\boo}_{\bg}) - \mpe(A_N f+ \Pi(\cF_h(\bg)), \boo) \ge h\big/\bigl(p^2(s-1)\bigr).
 	\end{equation*}
 	Since $f+ \Pi(\bg)$ satisfies Property~\textbf{B}, we have $p\le P$. Therefore, (iii) is true.
 \end{proof}
 
 \begin{proof}[\bf Proof of Theorem~\ref{thm:prop}~(iii)]
 	Fix $f\in \Lip(\Sigma_{\bF}, d_{\ba})$ with $\LIP(f) \le j+2$ and $\bg \sim \nu_{k}$. Direct calculations give $\LIP(f+ \Pi(\bg)) \le j+3$. Let $\cC $ denote the set of Lipschitz functions satisfying Properties~\textbf{A} and \textbf{B} but not the locking property. Set
 	\begin{equation}\label{eq:Ch}
 		h \= 4 (j+3) a_N (1-\theta)^{-2} P^2 (s-1).
 	\end{equation}
 	Applying Lemma~\ref{kem:Clem}~(iii) for $h$ and using (\ref{eq:Csufcon}), if $\bg \in \cF_h(D_{k,l_k-h})$, then $f+ \Pi(\bg) \notin \cC$. Thus, by Lemma~\ref{kem:Clem}~(i) and~(ii), we have
 	\begin{equation*}
 		\nu_k (\bg :f+\Pi(\bg) \in \cC) 
        \ge 1- \nu_{k} (\cF_h(D_{k,l_k}) )
        = 1- \nu_k(D_{k,l_k -h}) 
 		 =1- ( 1-h/l_k )^{Z_k}  
         < Z_k h /l_k.
 	\end{equation*}
 	Using (\ref{eq:deflk}), (\ref{eq:bbdef'}), and (\ref{eq:Ch}), we have
 	\begin{equation*}
 		\log( Z_k h /l_k ) = \operatorname{O} (\log N) + 2 N^{1/10} +  \bigl\lceil N^{1/5} -1 \bigr\rceil \log \theta,
 	\end{equation*}
 	which tends to $-\infty$ when $N$ tends to $+\infty$. Thus, $N$ can be chosen large enough such that $\nu_k (\{\bg :f+\Pi(\bg) \in \cC\}) \le 2^{-j-3}$. Using Lemma~\ref{lem:probreduce} and the push-forward $\nu_C (\{ g: f+g \in \cC \}) \le 2^{-j-2}$. Therefore, Theorem~\ref{thm:prop}~(iii) is verified.
 \end{proof}

\appendix
\section{Gauge estimates}
\begin{lemma}\label{lem:sameconc}
 	Let $k \in \N$, for every $\boo \in \cO_*$, if $\boo_1 \in \cO_{\Delta,k}$ and $\boo_2 \in \cO_*$ satisfy $\Sigma_k(\boo_2) \subseteq \Sigma_k(\boo_1)$, then $\boo_1 = \boo_2$.
 \end{lemma}
 
 \begin{proof}
 	Let $p \= \abs{\boo_1}$. Since $\Delta(\boo_1) > a_k$, for every $[\bw] \in \Sigma_k(\boo_1)$, $[\bw] \cap \boo_1$ consists of exactly one point. Write $\boo_1 = \{ x_i \}_{i\in \N}$ with $\sigma(x_i) = x_{i+1}$ and $x_1 = x_{p+1}$. For every  $i\in \N$, let $[\bw_i] \in \Sigma_k(\boo_1)$ such that $x_i = [\bw_i] \cap \boo_1$. Clearly, $[\bw_i] = [\bw_{i+p}]$ for every $i \in \N$. For every $i\in \N$, $\sigma([w_i]) \cap [\Sigma_k(\boo_1)] \subseteq [w_{i+1}]$. Otherwise, there exists $x\in [w_i]$ but $\sigma(x) \in [w_{i'}] \neq [w_{i+1}]$, implying that $d_{\ba}(x_{i+1}, x_{i'}) \le a_k$, contradicting to $\Delta(\boo_1) > a_k$. Now, let $x\in \boo_2$ and $i\in \N$ such that $x \in [w_i]$. For every $t \in \N$, $\sigma^t(x) \in \sigma^{t-1}([w_{i+1}]) = \dots = [w_{i+t}]$. This implies $x = x_i$. Hence, $\boo_2 \subseteq \boo_1$, implying that $\boo_2 = \boo_1$.
 \end{proof}

\begin{lemma}\label{lem:hilbert-brick}
Let $\ba=\{a_n\}_{n\in \N}$ be a strictly decreasing sequence of positive reals, $\theta\in(0,1)$, and assume that there exists $n_{\theta}\in \N$ such that $a_{n+1}\le\theta a_n$ for each $n\ge n_\theta$.
Let $\bb=\{b_{\bw}\}_{\bw\in\Sigma_*(\bF)}$ be a family of nonnegative numbers.
For each $n\in\N_0$, set
\begin{equation*}
    b_n\=\max\{b_{\bw}:\bw\in\Sigma_n(\bF)\}.
\end{equation*}
Suppose $b_n/a_n\to0$ as $n\to+\infty$.
Then for every choice of real coefficients $(g_{\bw})_{\bw\in\Sigma_*(\bF)}$ with $\abs{g_{\bw}}\le b_{\bw}$, the series
\begin{equation*}
    \sum_{n=0}^{+\infty}\sum_{\bw\in\Sigma_n(\bF)}g_{\bw}\mathbbold1_{[\bw]}
\end{equation*}
converges in $\Lip(\Sigma_{\bF},d_{\ba})$.
Moreover, the map $\Pi \: \cI[\bb] \to \Lip(\Sigma_{\bF}, d_{\ba})$ is continuous and $\cH[\bb]$ is compact in $\Lip(\Sigma_{\bF},d_{\ba})$.
In particular, $\bb$ is an $\ba$-admissible gauge.
\end{lemma}

\begin{proof}
Equip the set $\cI[\bb]\=\prod_{\bw\in\Sigma_*(\bF)}[-b_{\bw},b_{\bw}]$ with the product topology.
The space $\cI[\bb]$ is compact by Tychonoff's theorem.
For every $N\ge n_\theta$, set $\ve_N\=\sup\{b_n/a_n : n \ge N \}$.
The assumption $b_n/a_n\to0$ gives $\ve_N\to0$.
For each $N\ge n_\theta$, the bound on $a_{n+1}/a_n$ gives $a_n \le \theta^{n-N} a_N$ for every integer $n\ge N$, implying $\sum_{n=N}^{+\infty} a_n \le a_N/(1-\theta)$.
In particular, $\sum_{n=0}^{\infty}b_n<+\infty$.
Fix $\bg=(g_{\bw})_{\bw\in\Sigma_*(\bF)}\in \cI[\bb]$.
For every $n\in\N$, define
\begin{equation*}
    g_n\=\sum_{\bw\in\Sigma_n(\bF)}g_{\bw}\mathbbold{1}_{[\bw]}.
\end{equation*}
The cylinders $[\bw]$, with $\bw\in\Sigma_n(\bF)$, form a partition of $\Sigma_{\bF}$.
Hence, $\norm{g_n}_\infty\le b_n$, implying that $\Pi(\bg)\=\sum_{n=0}^{+\infty}g_n$
converges uniformly on $\Sigma_{\bF}$.
For each integer $N\ge1$, define
\begin{equation*}
    \Pi_N(\bg)\=\sum_{n=0}^{N-1}g_n \quad \text{ and } \quad R_N(\bg)\= \Pi(\bg)-\Pi_N(\bg).
\end{equation*}
For each integer $N\ge n_\theta$, direct calculations give
\begin{equation}\label{eq:hilbert-brick-2}
    \norm{R_N(\bg)}_\infty \le\sum_{n=N}^{+\infty}b_n \le \ve_N \sum_{n=N}^{+\infty} a_n \le\frac{\ve_Na_N}{1-\theta}.
\end{equation}
Fix $N\ge n_\theta$ and distinct $\bw,\, \bw'\in\Sigma_{\bF}$.
Let $m\ge1$ be the first coordinate at which $\bw$ and $\bw'$ differ.
Then $d_{\ba}(\bw,\bw')=a_m$.
For each integer $1\leq n<m$, the words $\bw$ and $\bw'$ belong to the same cylinder of length $n$, so $g_n(\bw)=g_n(\bw')$.
Set $n'\=\max\{N,\,m\}$.
Since $n'\ge N\ge n_\theta$,
\begin{equation*}
\begin{aligned}
    \abs{R_N(\bg)(\bw)-R_N(\bg)(\bw')} 
    &\le\sum_{n=n'}^{+\infty}\abs{g_n(\bw)-g_n(\bw')} 
    \le2\sum_{n=n'}^{+\infty}b_n\\
    &\le2\ve_N\sum_{n=n'}^{+\infty}a_n 
    \le\frac{2\ve_Na_{n'}}{1-\theta} 
    \le\frac{2\ve_N}{1-\theta}d_{\ba}(\bw,\bw').
\end{aligned}
\end{equation*}
Together with \eqref{eq:hilbert-brick-2}, the last estimate implies
\begin{equation}\label{eq:hilbert-brick-3}
    \norm{R_N(\bg)}_{\Lip} \le \ve_N(a_N+2)/ (1-\theta).
\end{equation}
Clearly, each cylinder indicator function $\mathbbold{1}_{[\bw]}$ is Lipschitz. Consequently, as a finite linear combination of cylinder indicator functions, $\Pi_N(\bg)$ is Lipschitz.
Equation \eqref{eq:hilbert-brick-3} now shows that $\Pi(\bg)$ is Lipschitz and that $\Pi_N(\bg)\to \Pi(\bg)$ in the Lipschitz norm.

For each $N\in\N$, the map $\Pi_N\: \cI[\bb]\rightarrow\Lip(\Sigma_{\bF},d_{\ba})$ is continuous, because $\Pi_N$ is a finite linear combination of coordinate projections with fixed Lipschitz functions as coefficients.
The bound in \eqref{eq:hilbert-brick-3} is independent of $\bg\in \cI[\bb]$ and tends to zero as $N\to +\infty$.
Hence the maps $\Pi_N$ converge uniformly to $\Pi\: \cI[\bb]\rightarrow\Lip(\Sigma_{\bF},d_{\ba})$, and $\Pi$ is continuous.
Since $\Pi(\cI[\bb])=\cH[\bb]$, the set $\cH[\bb]$ is compact.
\end{proof}

\begin{prop}\label{prop:finite-convolution}
For each $i\in\{1,\,2\}$, let $C_i\in(0,1)$ and $\bb^{(i)}$ be a $C_i$-log-$\ba$-admissible gauge, let $\nu_i$ be a Borel probability measure on $\Lip(\Sigma_{\bF},d_{\ba})$ with $\nu_i\bigl(\cH\bigl[\bb^{(i)}\bigr]\bigr)=1$.
For each $\bw\in\Sigma_*(\bF)$, set $b_{\bw}\= b_{\bw}^{(1)}+b_{\bw}^{(2)}$.
Then $\bb \= \{b_{\bw}\}_{\bw \in \Sigma_*(\bF)}$ is $\ba$-admissible and
\begin{equation*}
    (\nu_1*\nu_2)(\cH[\bb])=1.
\end{equation*}
If $C_1+C_2<1$, then $\nu_1\ast\nu_2$ is $(C_1+C_2)$-log-$\ba$-admissible.
\end{prop}

\begin{proof}
Let $n \in \mathbb{N}$. Note that
\begin{equation*}
    \max\{b_{\bw}/a_n:\bw\in\Sigma_n(\bF)\} \le (C_1+C_2) / \log n \longrightarrow0 \quad\text{ as }n\to+\infty.
\end{equation*}
Then, Lemma~\ref{lem:hilbert-brick} implies that $\bb$ is $\ba$-admissible.
For each $i\in\{1,\,2\}$ and each $g^{(i)}\in\cH\bigl[\bb^{(i)}\bigr]$, choose a coefficient representation
\begin{equation*}
    g^{(i)}=\sum_{n=0}^{+\infty}\sum_{\bw \in\Sigma_n} g_{\bw}^{(i)}\mathbbold1_{[\bw]},\quad \text{ and } \quad
    \Absbig{g_{\bw}^{(i)}}\le b_{\bw}^{(i)}.
\end{equation*}
The coefficients of the sum satisfy $\Absbig{g_{\bw}^{(1)}+g_{\bw}^{(2)}}\le b_{\bw}$.
Thus $g_1+g_2\in\cH[\bb]$.
By definition, $\nu_1*\nu_2$ is the law of the sum of independent $\Lip(\Sigma_{\bF},d_{\ba})$-valued random variables with laws $\nu_1$ and $\nu_2$.
Since $\cH\bigl[\bb^{(1)}\bigr]+\cH\bigl[\bb^{(2)}\bigr]\subseteq\cH[\bb]$, we have $(\nu_1*\nu_2)(\cH[\bb])=1$.
If $C_1+C_2<1$, then $\bb$ satisfies \eqref{eq:logadmis} with $C\=C_1+C_2$, which proves the last assertion.
\end{proof}

\begin{prop}\label{prop:infinite-convolution}
For each $j\in\N$, let $C_j\in(0,1)$, let $\bb^{(j)}$ be a $C_j$-log-$\ba$-admissible gauge, and let $\nu_j$ be a Borel probability measure on $\Lip(\Sigma_{\bF},d_{\ba})$ with $\nu_j\bigl(\cH\bigl[\bb^{(j)}\bigr]\bigr)=1$.
Assume that
\begin{equation}\label{eq:infinite-convolution-1}
    C_*\=\sum_{j=1}^{+\infty}C_j<+\infty.
\end{equation}
For each $\bw\in\Sigma_*(\bF)$, define $b_{\bw}\=\sum_{j=1}^{+\infty}b_{\bw}^{(j)}$.
Then $\bb$ is $\ba$-admissible, and the finite convolutions $\nu_1*\cdots*\nu_M$ converge weakly to a Borel probability measure $\nu$ as $M\to + \infty$ on $\Lip(\Sigma_{\bF},d_{\ba})$, with $\nu(\cH[\bb])=1$.
For independent $\Lip(\Sigma_{\bF},d_{\ba})$-valued random variables $g_j$ with laws $\nu_j$, the series $\sum_{j=1}^{+\infty}g_j$ converges absolutely in $\Lip(\Sigma_{\bF},d_{\ba})$ almost surely, and its law is $\nu$.
If $C_*<1$, then $\nu$ is $C_*$-log-$\ba$-admissible.
\end{prop}

\begin{proof}
For each $j\in\N$, set $\cH_j\=\cH\bigl[\bb^{(j)}\bigr]$, which is compact.
We first obtain a bound for each $j\in\N$ and $g\in \cH_j$.
Since the coefficients vanish below level $N$, and since the cylinders of the same length are disjoint, we have
\begin{equation*}
    \norm{g}_\infty \le C_j\sum_{n=N}^{+\infty}\frac{a_n}{\log n} \le\frac{C_ja_N}{(1-\theta)\log N}.
\end{equation*}
For distinct $\bw,\,\bw'\in\Sigma_{\bF}$ with $d_{\ba}(\bw, \bw') = a_m$, setting $n'\=\max\{N,\,m\}$, we have
\begin{equation*}
\begin{aligned}
    \abs{g(\bw)-g(\bw')} \le2C_j\sum_{n=n'}^{+\infty}\frac{a_n}{\log n}\le\frac{2C_ja_{n'}}{(1-\theta)\log N} \le\frac{2C_j}{(1-\theta)\log N}d_{\ba}(\bw,\bw').
\end{aligned}
\end{equation*}
Set $R\=(a_N+2)/((1-\theta)\log N)$.
For each $j\in \N$ and $g\in \cH_j$, the preceding estimates give
\begin{equation}\label{eq:infinite-convolution-2}
    \norm{g}_{\Lip}\le RC_j.
\end{equation}

Equip $\cI\=\prod_{j=1}^{+\infty}\cH_j$ with the product topology.
The space $\cI$ is compact.
For each $M\in\N$, define the continuous map $\Pi_M'\: \cI\rightarrow\Lip(\Sigma_{\bF},d_{\ba}),$
\begin{equation*}
     \Pi_M'(\tg)\=\sum_{j=1}^{M}g_j.
\end{equation*}
By \eqref{eq:infinite-convolution-1} and \eqref{eq:infinite-convolution-2}, for each $\tg=\{g_j\}_{j\in\N}\in \cI$,
\begin{equation*}
    \sum_{j=1}^{+\infty}\norm{g_j}_{\Lip}\le RC_*<+\infty.
\end{equation*}
Since $\Lip(\Sigma_{\bF},d_{\ba})$ is a Banach space, the sum defines a map $\Pi'\: \cI\rightarrow\Lip(\Sigma_{\bF},d_{\ba}),$
\begin{equation*}\label{eq:infinite-convolution-3}
     \Pi'(\tg)\=\sum_{j=1}^{+\infty}g_j, \;\;\; \norm{\Pi'(\tg)-\Pi_M'(\tg)}_{\Lip} \le R\sum_{j=M+1}^{+\infty} C_j.
\end{equation*}
The right-hand side tends to zero as $M\to+\infty$ and is independent of $\{g_j\}_{j\in\N}\in \cI$.
Thus $\Pi_M'\to \Pi'$ uniformly, so $\Pi$ is continuous.
In particular, $\Pi'(\cI)$ is compact.

Regard each $\nu_j$ as a probability measure on $\cH_j$, and define
\begin{equation*}
    \nu_{\times}\=\bigotimes_{j=1}^{+\infty}\nu_j, \qquad \nu\= \Pi'_*\nu_{\times}.
\end{equation*}
Then $(\Pi_M)'_*\nu_{\times}=\nu_1\ast\cdots\ast\nu_M$.
For each bounded continuous function $\varphi\:\Lip(\Sigma_{\bF},d_{\ba})\rightarrow\R$, the dominated convergence theorem gives
\begin{equation*}
\begin{aligned}
    \int_{\Lip(\Sigma_{\bF},d_{\ba})}\!\varphi\,\mathrm{d}(\nu_1*\cdots*\nu_M) 
    &=\int_{\cI}\!\varphi(\Pi'_M(g))\,\mathrm{d}\nu_{\times}(g)
    \longrightarrow\int_{\cI}\!\varphi(\Pi'(g))\,\mathrm{d}\nu_{\times}(g) 
    =\int_{\Lip(\Sigma_{\bF},d_{\ba})}\!\varphi\,\mathrm{d}\nu.
\end{aligned}
\end{equation*}
Hence the finite convolutions converge weakly to $\nu$.
The coordinate random variables on $(\cI,\nu_{\times})$ are independent with laws $\nu_j$, and their sum has law $\nu$.
For any other independent family with the same laws, all $g_j\in \cH_j$ simultaneously with probability one, and the joint law on $\cI$ is $\nu_{\times}$.
The same norm bound gives absolute convergence, and the sum again has law $\nu$.

The total gauge $\bb$ vanishes on words of length less than $N$.
For each $\bw\in\Sigma_n(\bF)$ with $\abs{\bw}\ge N$, the series defining $b_{\bw}$ is finite and satisfies
\begin{equation*}
    b_{\bw}\le C_*a_{\abs{\bw}}/\log\abs{\bw} .
\end{equation*}
It follows that
\begin{equation*}
    \max\{b_{\bw}/a_n:\bw\in\Sigma_n(\bF)\} \le C_*/\log n \longrightarrow0 \quad\text{ as }n\to +\infty.
\end{equation*}
Lemma~\ref{lem:hilbert-brick} shows that $\cH[\bb]$ is compact and hence closed in $\Lip(\Sigma_{\bF},d_{\ba})$.
Adding finitely many coefficient representations gives $\Pi_M(\cI)\subseteq\cH[\bb]$ for each $M\in\N$.
Passing to the limit gives $\Pi(\cI)\subseteq\cH[\bb]$.
Thus $\nu(\cH[\bb])=1$.
If $C_*<1$, then $\bb$ is $C_*$-log-$\ba$-admissible, and so is $\nu$.
\end{proof}

\section{Locking criterion proof}
\label{sec:lockingproof}
Here we prove Lemma~\ref{lemma.betastar}. We follow the steps in \cite{BZ16}. We begin with Lemma~4.3 of \cite{BZ16}, which is proven identically.
\begin{lemma}    \label{lem:b1}
    Suppose $\boo \in \PO_{*}$ satisfies $\Delta(\boo) > a_n$. Then for each $s \in \N$,
    \begin{equation*}
        \{\bw \in \Sigma_{\bF} : \bw_{[1,n+s]} \in \Sigma_{n+s}(\boo) \} 
        = \bigcap_{i = 0}^s \bigl\{\bw \in \Sigma_{\bF} : \sigma^i(\bw)_{[1,n]} \in \Sigma_{n}(\boo)\bigr\}.
    \end{equation*}
    \vspace{-\baselineskip}
\end{lemma}
We continue by restating Lemma~4.4 of \cite{BZ16}:
\begin{lemma}
    Suppose $m\ge n\geq 1$ are integers, $\boo,\, \boo' \in \PO_*$ satisfy $\Delta(\boo) > a_n$, and $\Delta(\boo') >a_m$. Assume $\abs{\boo'} \ge \abs{\boo}$. Then, if $f = A_mf$, then
    \begin{equation*}
        \mpe(\boo', f) - \mpe(\boo, f) 
        \le 2(m-n+1) \cdot \norm{f}_{\infty}\cdot \mu_{\boo'}([\Sigma_n^c(\boo)]).
    \end{equation*}
    \vspace{-\baselineskip}
\end{lemma}
Then, as $m > n_{\theta}$, we have that the existence of $\nu_k$ and $u_k$ of \cite{BZ16} and the proof of Lemma~4.5 of \cite{BZ16} follow from Lemma~\ref{lemma.gapcompose}. Afterwards, we can define
\begin{equation*}
    g_k \= A_{k+1}f - A_kf.
\end{equation*}
Then, using \eqref{q.modifiedestimate}, we get the modified Lemma~4.1 of \cite{BZ16}, which gives the following locking criterion:
\begin{equation*}
    \ouv_m^{\sharp}(f) > 2\sum_{k \ge m} (k-m+1) \norm{g_k}_{\infty}.
\end{equation*}
We use \eqref{eq:fAmfclose} and the estimates on $a_k$ to bound the right-hand side:
\begin{equation*}
    \begin{aligned}
        2\sum_{k \ge m} (k-m+1) \norm{g_k}_{\infty} 
        &\le 2\sum_{k \ge m} (k-m+1)( \norm{f-A_k f}_{\infty} + \norm{f-A_{k+1}}_{\infty})\\
        & \le  2\LIP(f)\sum_{k \ge m} (k-m+1) (a_{k+1}+a_{k+2})  \\
        &\le  4\LIP(f)a_m \sum_{k \ge m} (k-m+1)\theta^{k-m+1}\\
        & \le 4\LIP(f)a_m \sum_{k\ge 1} k\theta^k \\
        &\le 4\LIP(f)a_m(1-\theta)^{-2}. 
    \end{aligned}
\end{equation*}
This finishes the proof of Lemma~\ref{lemma.betastar}.

\end{document}